\documentclass[12pt,a4paper]{article}

\usepackage[english]{babel} 
\usepackage{ucs}
\usepackage[utf8]{inputenc}

\RequirePackage{color} 

\usepackage[pdftex,colorlinks=true,urlcolor=blue,pdfstartview=FitH]{hyperref}
\usepackage[pdftex]{graphicx}
\DeclareGraphicsExtensions{.pdf, .png, .jpg, .pdftex} 

\RequirePackage{amsmath,amssymb,amsfonts,latexsym,mathrsfs}

\newcommand{\A}{\mathcal{A}}

\newcommand{\ee}{\mathbf{e}}

\newcommand{\TT}{\mathcal{T}}

\newcommand{\CC}{\mathcal{C}}

\newcommand{\dgh}{\mathtt{d}_{{\rm GH}}}

\newcommand{\dq}{d_{\mathbf{q}}} 

\newcommand{\bq}{\mathbf{q}}
 \newcommand{\bp}{\mathbf{p}}
\newcommand{\B}{\mathbb{B}}
 
\newcommand{\bm}{\mathbf{m}}
\newcommand{\bt}{\mathbf{t}}
 
\newcommand{\bl}{{\bf l}} 

\newcommand{\bv}{\mathbf{v}}
 
 \newcommand{\bS}{{\bf S}}
 \newcommand{\bQ}{{\bf Q}}
 \newcommand{\RR}{{\cal R}}

\newcommand{\sg}{\mathfrak{s}}

\newcommand{\R}{\mathbb{R}}
\newcommand{\N}{\mathbb{N}} 
\newcommand{\M}{\mathbb{M}}
\newcommand{\E}{\mathbb{E}}
 
\newcommand{\Z}{\mathbb{Z}}

\newcommand{\Q}{\mathbb{Q}} 
 
\newcommand{\FF} { {\cal F }} 
\def\build#1_#2^#3{\mathrel{ \mathop{\kern 0pt#1}\limits_{#2}^{#3}}}

\def\cq{$\hfill \square$}
\def\un{\underline}
\def\d{\mathrm{d}}
\def\n{\mathcal{N}}

\def\eps{\varepsilon}

\def\S{{\sf S}}

\def\bLM{\mathbf{LM}}

\def\bCLM{\mathbf{CLM}}
\def\LM{{\sf LM}}
\def\CLM{{\sf CLM}}
\def\ov{\overline}

\def\proof{\noindent{\bf Proof. }}

\renewcommand{\P}{\mathbb{P}}

\newcommand{\btau}{\boldsymbol{\tau}}
\newcommand{\ind}{\mathbf{1}}

\newtheorem{thm}{Theorem}
\newtheorem{lmm}[thm]{Lemma}
\newtheorem{prp}[thm]{Proposition}
\newtheorem{defn}[thm]{Definition}
\newtheorem{crl}[thm]{Corollary}

\usepackage{fancyhdr}

\usepackage[outerbars]{changebar} 

\title{The Brownian map is the scaling limit of uniform random plane
  quadrangulations}

\author{Grégory Miermont\\
  Université de Paris-Sud}

\begin{document}

\selectlanguage{english}

\maketitle

\begin{abstract}
  We prove that uniform random quadrangulations of the sphere with $n$
  faces, endowed with the usual graph distance and renormalized by
  $n^{-1/4}$, converge as $n\to\infty$ in distribution for the
  Gromov-Hausdorff topology to a limiting metric space. We validate a
  conjecture by Le~Gall, by showing that the limit is (up to a scale
  constant) the so-called {\em Brownian map}, which was introduced by
  Marckert \& Mokkadem and Le Gall as the most natural candidate for
  the scaling limit of many models of random plane maps. The proof
  relies strongly on the concept of {\em geodesic stars} in the map,
  which are configurations made of several geodesics that only share a
  common endpoint and do not meet elsewhere.
\end{abstract}

\tableofcontents

\section{Introduction}

Random plane maps and their scaling limits are a sort of
two-dimensional analog of random walks and Brownian motion, in which
one wants to approximate a random continuous surface using large
random graphs drawn on the $2$-sphere \cite{ADJ}. A {\em plane map} is
a proper embedding, without edge-crossings, of a finite connected
graph in the two-dimensional sphere. Loops or multiple edges are
allowed. We say that a map is {\em rooted} if one of its oriented
edges, called the root, is distinguished. Two (rooted) maps $\bm,\bm'$
are considered equivalent if there exists a direct homeomorphism of
the sphere that maps $\bm$ onto $\bm'$ (and maps the root of $\bm$ to
that of $\bm'$ with preserved orientations). Equivalent maps will
systematically be identified in the sequel, so that the set of maps is
a countable set with this convention.

From a combinatorial and probabilistic perspective, the maps called
{\em quadrangulations}, which will be the central object of study in
this work, are among the simplest to manipulate.  Recall that the {\em
  faces} of a map are the connected component of the complement of the
embedding.  A map is called a quadrangulation if all its faces have
degree $4$, where the degree of a face is the number of edges that are
incident to it (an edge which is incident to only one face has to be
counted twice in the computation of the degree of this face).  Let
$\bQ$ be the set of plane, rooted quadrangulations, and $\bQ_n$ be the
set of elements of $\bQ$ with $n$ faces.  Let $Q_n$ be uniformly
distributed in $\bQ_n$. We identify $Q_n$ with the finite metric space
$(V(Q_n),d_{Q_n})$, where $V(Q_n)$ is the set of vertices of $Q_n$,
and $d_{Q_n}$ is the usual graph distance on $V(Q_n)$. We see $Q_n$ as
a random variable with values in the space $\M$ of compact metric
spaces considered up to isometry. The space $\M$ is endowed with the
Gromov-Hausdorff topology \cite{burago01}: The distance between two
elements $(X,d),(X',d')$ in $\M$ is given by
$$\dgh((X,d),(X',d'))=\frac{1}{2}\inf_{\mathcal{R}}\mathrm{dis}(\mathcal{R})\,
,$$ where the infimum is taken over all {\em correspondences}
$\mathcal{R}$ between $X$ and $X'$, i.e.\ subsets of $X\times X'$
whose canonical projections $\mathcal{R}\to X,X'$ are onto, and
$\mathrm{dis}(\mathcal{R})$ is the {\em distortion} of $\mathcal{R}$,
defined by 
$$\mathrm{dis}(\mathcal{R})=\sup_{(x,x'),(y,y')\in
  \mathcal{R}}|d(x,y)-d'(x',y')|\, .$$ The space $(\M,\dgh)$ is then a
separable and complete metric space \cite{evpiwin}.


It turns out that typical graph distances in $Q_n$ are of order
$n^{1/4}$ as $n\to\infty$, as was shown in a seminal paper by
Chassaing \& Schaeffer \cite{CSise}. Since then, much attention has
been drawn by the problem of studying the scaling limit of $Q_n$,
i.e.\ to study the asymptotic properties of the space
$n^{-1/4}Q_n:=(V(Q_n),n^{-1/4}d_{Q_n})$ as $n\to\infty$. Le Gall
\cite{legall06} obtained a compactness result, namely that the laws of
$(n^{-1/4}Q_n,n\geq 1)$, form a relatively compact family of
probability distributions on $\M$. This means that along suitable
subsequences, the sequence $n^{-1/4}Q_n$ converges in distribution in
$\M$ to some limiting space. Such limiting spaces are called
``Brownian maps'', so as to recall the fact that Brownian motion
arises as the scaling limit of discrete random walks. Many properties
satisfied by {\em any} Brownian map (i.e.\ by any limit in
distribution of $n^{-1/4}Q_n$ along some subsequence) are known. In
particular, Le Gall showed that their topology is independent of the
choice of a subsequence \cite{legall06}. Then Le Gall \& Paulin
identified this topology with the topology of the 2-dimensional sphere
\cite{lgp,miermontsph}. Besides, the convergence of {\em two-point
  functions} and {\em three-point functions}, that is, of the joint
laws of rescaled distances between $2$ or $3$ randomly chosen vertices
in $V(Q_n)$, was also established respectively by Chassaing \&
Schaeffer \cite{CSise} and Bouttier \& Guitter \cite{BoGu08b}. These
convergences occur without having to extract subsequences. See
\cite{LGMi11} for a recent survey of the field.

It is thus natural to conjecture that all Brownian maps should in fact
have the same distribution, and that the use of subsequences in the
approximation by quadrangulations is superfluous. A candidate for a
limiting space (which is sometimes also called the Brownian map,
although is was not proved that this space arises as the limit of
$n^{-1/4}Q_n$ along some subsequence) was described in equivalent, but
slightly different forms, by Marckert and Mokkadem \cite{MM05} and Le
Gall \cite{legall06}. 

The main goal of this work is to prove these conjectures, namely, that
$n^{-1/4}Q_n$ converges in distribution as $n\to\infty$ to the
conjectured limit of \cite{MM05,legall06}. This unifies the several
existing definitions of Brownian map we just described, and lifts the
ambiguity that there could have been more than one limiting law along
a subsequence for $n^{-1/4}Q_n$.

In order to state our main result, let us describe the limiting
Brownian map. This space can be described from a pair of random
processes $(\ee,Z)$. Here, $\ee=(\ee_t,0\leq t\leq 1)$ is the
so-called {\em normalized Brownian excursion}. It can be seen as a
positive excursion of standard Brownian motion conditioned to have
duration $1$. The process $Z=(Z_t,0\leq t\leq 1)$ is the so-called
{\em head of the Brownian snake} driven by $\ee$: Conditionally given
$\ee$, $Z$ is a centered Gaussian process with continuous
trajectories, satisfying
$$E[|Z_s-Z_t|^2\, |\, \ee]=d_\ee(s,t)\, \qquad s,t\in [0,1]\, ,$$
where
$$d_\ee(s,t)=\ee_s+\ee_t-2\inf_{s\wedge t\leq u\leq s\vee t}\ee_u\,
.$$ The function $d_\ee$ defines a pseudo-distance on $[0,1]$: This
means that $d_\ee$ satisfies the properties of a distance, excepting
separation, so that it can hold that $d_\ee(s,t)=0$ with $s\neq t$.
We let $\TT_\ee=[0,1]/\{d_\ee=0\}$, and denote the canonical
projection by $p_\ee:[0,1]\to \TT_\ee$. It is obvious that $d_\ee$
passes to the quotient to a distance function on $\TT_\ee$, still
called $d_\ee$.  The space $(\TT_\ee,d_\ee)$ is called Aldous' {\em
  continuum random tree}. The definition of $Z$ implies that a.s., for
every $s,t$ such that $d_\ee(s,t)=0$, one has $Z_s=Z_t$. Therefore, it
is convenient to view $Z$ as a function $(Z_a,a\in \TT_\ee)$ indexed
by $\TT_\ee$. Roughly speaking, $Z$ can be viewed as the Brownian
motion indexed by the Brownian tree.

Next, we set
$$D^\circ(s,t)=Z_s+Z_t-2\max\Big(\inf_{s\leq u\leq  t}Z_u,\inf_{t\leq
  u\leq s}Z_u\Big)\,
,\qquad s,t\in [0,1]\, ,$$
where $s\leq u\leq t$ means $u\in [s,1]\cup [0,t]$ when $t<s$. Then
let, for $a,b\in \TT_\ee$, 
$$D^\circ(a,b)=\inf\{D^\circ(s,t):s,t\in
[0,1],p_\ee(s)=a,p_\ee(t)=b\}\, .$$ The function $D^\circ$ on
$[0,1]^2$ is a pseudo-distance, but $D^\circ$ on $\TT_\ee^2$ does not
satisfy the triangle inequality. This motivates writing, for $a,b\in
\TT_\ee$, 
$$D^*(a,b)=\inf\Big\{\sum_{i=1}^{k-1}D^\circ(a_i,a_{i+1}):k\geq
1,a=a_1,a_2,\ldots, a_{k-1},a_k=b\Big\}\, .$$ 
The function $D^*$ is now a pseudo-distance on $\TT_\ee$, and we
finally define
$$S=\TT_\ee/\{D^*=0\}\, ,$$
which we endow with the quotient distance, still denoted by
$D^*$. Alternatively, 
letting, for $s,t\in [0,1]$,
$$D^*(s,t)=\inf\Big\{\sum_{i=1}^{k}D^\circ(s_i,t_i):k\geq
1,s=s_1,t=t_k, d_\ee(t_i,s_{i+1})=0,1\leq i\leq k-1\Big\}\, ,$$ one
can note that $S$ can also be defined as the quotient metric space
$[0,1]/\{D^*=0\}$. The space $(S,D^*)$ is a {\em geodesic metric
  space}, meaning that for every two points $x,y\in S$, there exists
an isometry $\gamma:[0,D^*(x,y)]\to S$ such that
$\gamma(0)=x,\gamma(D^*(x,y))=y$. The function $\gamma$ is called a
{\em geodesic} from $x$ to $y$.

The main result of this paper is the following.

\begin{thm}\label{sec:introduction-4}
  As $n\to\infty$, the metric space $(V(Q_n),(8n/9)^{-1/4}d_{Q_n})$
  converges in distribution for the Gromov-Hausdorff topology on $\M$
  to the space $(S,D^*)$.
\end{thm}

The strategy of the proof is to obtain properties of the geodesic
paths that have to be satisfied in any distributional limit of
$n^{-1/4}Q_n$, and which are sufficient to give ``formulas'' for
distances in these limiting spaces that do not depend on the choice of
subsequences. An important object of study is the occurence of certain
types of {\em geodesic stars} in the Brownian map, i.e.\ of points
from which radiate several disjoint geodesic {\em arms}. We hope that
the present study will pave the way to a deeper understanding of
geodesic stars in the Brownian map. 


The rest of the paper is organised as follows. The next section
recalls an important construction from \cite{legall06} of the
potential limits in distributions of the spaces $n^{-1/4}Q_n$, which
will allow to reformulate slightly Theorem \ref{sec:introduction-4}
into the alternative Theorem \ref{sec:introduction}.  Then we show how
the proof of Theorem \ref{sec:introduction} can be obtained as a
consequence of two intermediate statements, Propositions
\ref{sec:introduction-1} and
\ref{sec:th1}. Section~\ref{sec:rough-comp-betw} is devoted to the
proof of the first proposition, while Section
\ref{sec:exceptional-nature-3-4} reduces the proof of Proposition
\ref{sec:th1} to two key ``elementary'' statements, Lemmas
\ref{sec:covering-3-star} and \ref{sec:covering-3-star-1}, which deal
with certain properties of geodesic stars in the Brownian map, with
two or three arms. The proofs of these lemmas contains most of the
novel ideas and techniques of the paper. Using a generalization of the
Cori-Vauquelin-Schaeffer bijection found in \cite{miertess}, we will
be able to translate the key statements in terms of certain
probabilities for families of labeled maps, which have a simple enough
structure so that we are able to derive their scaling limits. These
discrete and continuous structures will be described in Sections
\ref{sec:coding-labeled-maps} and \ref{sec:scal-limits-label}, while
Section \ref{sec:proof-key-lemmas} is finally devoted to the proof of
the two key lemmas.

Let us end this Introduction by mentioning that, simultaneously
to the elaboration of this work, Jean-François Le Gall independently
found a proof of Theorem \ref{sec:introduction-4}. His method is
different from ours, and we believe that both approaches present
specific interests.

\paragraph{Notation conventions. }  In this paper, we let
$V(\bm),E(\bm)$ and $F(\bm)$ be the sets of vertices, edges and faces
of the map $\bm$.  Also, we let $\vec{E}(\bm)$ be the set of oriented
edges of $\bm$, so that $\#\vec{E}(\bm)=2\#E(\bm)$.  If $e\in
\vec{E}(\bm)$, we let $e_-e_+\in V(\bm)$ be the origin and target of
$e$. The reversal of $e\in \vec{E}(\sg)$ is denoted by $\ov{e}$.

If $e\in \vec{E}(\bm),f\in F(\bm),v\in V(\bm)$, we say that $f$ and
$e$ are incident if $f$ lies {\em to the left of $e$} when following
the orientation of $e$. We say that $v$ and $e$ are incident if
$v=e_-$. If $e\in E(\bm)$ is not oriented, we say that $f$ and $e$ are
incident if $f$ is incident to $e$ or $\ov{e}$. A similar definition
holds for incidence between $v$ and $e$. Finally, we say that $f$ and
$v$ are incident if $v$ and $f$ are incident to a common $e\in
E(\bm)$. The sets of vertices, edges and oriented edges incident to a
face $f$ are denoted by $V(f),E(f),\vec{E}(f)$. We will also use the
notation $V(f\cap f'),E(f\cap f')$ for vertices and edges
simultaneously incident to the two faces $f,f'\in F(\bm)$. 

If $\bm$ is a map and $v,v'$ are vertices of $\bm$, a {\em chain} from
$v$ to $v'$ is a finite sequence $(e^{(1)},e^{(2)},\ldots,e^{(k)})$ of
oriented edges, such that $e^{(1)}_-=v,e^{(k)}_+=v'$, and
$e^{(i)}_+=e^{(i+1)}_-$ for every $i\in \{1,2,\ldots,k-1\}$. The
integer $k$ is called the length of the chain, and we allow also the
chain of length $0$ from $v$ to itself. The {\em graph distance}
$d_\bm(v,v')$ in $\bm$ described above is then the minimal $k$ such
that there exists a chain with length $k$ from $v$ to $v'$.  A chain
with minimal length is called a {\em geodesic chain}.

In this paper, we will often let
$C$ denote a positive, finite constant, whose value may vary from line
to line. Unless specified otherwise, the random variables considered
in this paper are supposed to be defined on a probability space
$(\Omega,\FF,P)$.

\section{Preliminaries}\label{sec:brownian-map}

\subsection{Extracting distributional limits from large
  quadrangulations}\label{sec:extr-distr-limits}

As mentioned in the introduction, it is known that the laws of the
rescaled quadrangulations $(8n/9)^{-1/4}Q_n$ form a compact sequence
of distributions on $\M$. Therefore, from every subsequence, it is
possible to further extract a subsequence along which
$(8n/9)^{-1/4}Q_n$ converges in distribution to a random variable
$(S',D')$ with values in $\M$. Theorem \ref{sec:introduction-4} then
simply boils down to showing that this limit has the same distribution
as $(S,D^*)$, independently on the choices of subsequences.

In order to be able to compare efficiently the spaces $(S',D')$ and
$(S,D^*)$, we perform a particular construction, due to
\cite{legall06}, for which the spaces are the same quotient space,
i.e.\ $S=S'$. This is not restrictive, in the sense that this
construction can always be performed up to (yet) further
extraction. We recall some of its important aspects.

Recall that a quadrangulation $\bq\in \bQ_n$, together with a
distinguished vertex $v_*$ can be encoded by a labeled tree with $n$
edges {\em via} the so-called Cori-Vauquelin-Schaeffer bijection,
which was introduced by Cori and Vauquelin \cite{CoVa}, and
considerably developed by Schaeffer, starting with his thesis
\cite{schaeffer98}. If $(\bt,\bl)$ is the resulting labeled tree, then
the vertices of $\bq$ distinct from $v_*$ are exactly the vertices of
$\bt$, and $\bl$ is, up to a shift by its global minimum over $\bt$,
the function giving the graph distances to $v_*$ in $\bq$. In turn,
the tree $(\bt,\bl)$ can be conveniently encoded by {\em contour and
  label functions}: Heuristically, this function returns the height
(distance to the root of $\bt$) and label of the vertex visited at
time $0\leq k\leq 2n$ when going around the tree clockwise. These
functions are extended by linear interpolation to continuous functions
on $[0,2n]$.

If $\bq=Q_n$ is a uniform random variable in $\bQ_n$ and $v_*$ is
uniform among the $n+2$ vertices of $Q_n$, then the resulting labeled
tree $(T_n,\ell_n)$ has a contour and label function $(C_n,L_n)$ such
that $C_n$ is a simple random walk in $\Z$, starting at $0$ and ending
at $0$ at time $2n$, and conditioned to stay non-negative. Letting
$u^n_i$ be the vertex of $\bt_n$ visited at step $i$ of the contour,
we let
$$D_n(i/2n,j/2n)=\Big(\frac{9}{8n}\Big)^{1/4}d_{Q_n}(u^n_i,u^n_j)
\, ,\qquad 0\leq i,j\leq 2n\, ,$$ and extend $D_n$ to a continuous
function on $[0,1]^2$ by interpolation, see \cite{legall06} for
details. Then, the distributions of the triples of processes
$$\Bigg(\bigg(\frac{C_n(2ns)}{\sqrt{2n}}\bigg)_{0\leq s\leq
  1},\bigg(\frac{L_n(2ns)}{(8n/9)^{1/4}}\bigg)_{0\leq s\leq
  1},(D_n(s,t))_{0\leq s,t\leq 1}\Bigg)\, ,\qquad n\geq 1$$ form a
relatively compact family of probability distributions. Therefore,
from every subsequence, we can further extract a certain subsequence
$(n_k)$, along which the above triples converge in distribution to a
limit $(\ee,Z,D)$, where $D$ is a random pseudo-distance on $[0,1]$, and
$(\ee,Z)$ is the head of the Brownian snake process described in the
Introduction.  We may and will assume that this convergence holds
a.s., by using the Skorokhod representation theorem, and up to
changing the underlying probability space.  Implicitly, until Section
\ref{sec:exceptional-nature-3-4}, all the integers $n$ and limits
$n\to\infty$ that are considered will be along the subsequence
$(n_k)$, or some further extraction.

The function $D$ is a class function for $\{d_\ee=0\}$, which induces
a pseudo-distance on $\TT_\ee/\{d_\ee=0\}$, still denoted by $D$ for
simplicity.  Viewing successively $D$ as a pseudo-distance on $[0,1]$
and $\TT_\ee$, we can let
$$S' =[0,1]/\{D=0\}=\TT_\ee/\{D=0\}\, , $$ 
and endow it with the distance induced by $D$, still written $D$ by a
similar abuse of notation as above for simplicity.  The space $(S',D)$
is then a random geodesic metric space.

On the other hand,
we can define $D^\circ$ and $D^*$ out of $(\ee,Z)$ as in the
Introduction, and let $S=\TT_\ee/\{D^*=0\}$. The main result of
\cite{legall06} is the following.

\begin{prp}[Theorem 3.4 in \cite{legall06}]
\noindent {\rm (i)}
The three subsets of $\TT_\ee\times \TT_\ee$
$$\{D=0\}\, ,\quad \{D^\circ=0\}\, ,\quad \{D^*=0\}$$ 
are a.s.\ the same. In particular, the quotient sets $S'$ and $S$ are
a.s.\ equal, and $(S,D),(S,D^*)$ are homeomorphic. 

\noindent{\rm (ii)}
Along the subsequence $(n_k)$, $(Q_n,(9/8n)^{1/4}d_{Q_n})$ converges
a.s.\ in the Gromov-Hausdorff sense to $(S,D)$.
\end{prp}

Using the last statement, Theorem \ref{sec:introduction-4} is now a
consequence of the following statement, which is the result that we
are going to prove in the remainder of this paper.

\begin{thm}\label{sec:introduction}
Almost-surely, it holds that $D=D^*$. 
\end{thm}

\subsection{A short review of results on $S$}\label{sec:short-review-results}

\paragraph{A word on notation. }
Since we are considering several metrics $D,D^*$ on the same set $S$,
a little care is needed when we consider balls or geodesics, as we
must mention to which metric we are referring to. For $x\in S$ and
$r\geq 0$, we let
$$B_D(x,r)=\{y\in S:D(x,y)<r\}\, ,\qquad B_{D^*}(x,r)=\{y\in
S:D^*(x,y)<r\}\, ,$$ and we call them respectively the (open) $D$-ball
and the $D^*$-ball with center $r$ and radius $r$. Similarly, a
continuous path $\gamma$ in $S$ will be called a $D$-geodesic, resp.\
a $D^*$-geodesic, if it is a geodesic path in $(S,D)$ resp.\
$(S,D^*)$. Note that since $(S,D)$ and $(S,D^*)$ are a.s.\
homeomorphic, a path in $S$ is continuous for the metric $D$ if and
only if it is continuous for the metric $D^*$. When it is unambiguous
from the context which metric we are dealing with, we sometimes omit
the mention of $D$ or $D^*$ when considering balls or geodesics.

\medskip

\paragraph{Basic properties of $D,D^\circ,D^*$. }
Recall that $p_\ee:[0,1]\to \TT_\ee$
is the canonical projection, we will also let $\bp_Z:\TT_\ee\to S$ be
the canonical projection, and $\bp=\bp_Z\circ p_\ee$.  The function
$D^\circ$ was defined on $[0,1]^2$ and $(\TT_\ee)^2$, it also induces
a function on $S^2$ by letting
$$D^\circ(x,y)=\inf\{D^\circ(a,b):a,b\in
\TT_\ee,\bp_Z(a)=x,\bp_Z(b)=y\}\, .$$ Again, $D^\circ$ does not
satisfy the triangle inequality. However, it holds that
$$D(x,y)\leq D^*(x,y)\leq
D^\circ(x,y)\, ,\qquad x,y\in S\, .$$
One of the difficulties in
handling $S$ is its definition using two successive quotients, so we
will always mention whether we are considering $D,D^*,D^\circ$ on
$[0,1],\TT_\ee$ or $S$.

We define the {\em uniform measure}
$\lambda$ on $S$ to be the push-forward of the Lebesgue measure on
$[0,1]$ by $\bp$. This measure will be an important tool to sample
points randomly on $S$.

Furthermore, by \cite{legweill}, the process $Z$ attains its overall
minimum a.s.\ at a single point $s_*\in [0,1]$, and the class
$\rho=\bp(s_*)$ is called the {\em root} of the space $(S,D)$.  One
has, by \cite{legall06},
\begin{equation}\label{eq:14}
D(\rho,x)=D^\circ(\rho,x)=D^*(\rho,x)=Z_x-\inf Z\, ,\qquad \mbox{ for every
}x\in S\, .
\end{equation}
Here we viewed $Z$ as a function on $S$, which is licit because $Z$ is
a class function for $\{D=0\}$, coming from the fact that $D(s,t)\geq
|Z_s-Z_t|$ for every $s,t\in [0,1]$. The latter property can be easily
deduced by passing to the limit from the discrete counterpart
$D_n(i/2n,j/2n)\geq (9/8n)^{1/4}|L_n(i)-L_n(j)|$, which is a
consequence of standard properties of the Cori-Vauquelin-Schaeffer
bijection.

\paragraph{Geodesics from the root in $S$. }
Note that the definition of the function $D^\circ$ on $[0,1]^2$ is
analogous to that of $d_\ee$, using $Z$ instead of $\ee$. Similarly to
$\TT_\ee$, we can consider yet another quotient
$\TT_Z=[0,1]/\{D^\circ=0\}$, and endow it with the distance induced by
$D^\circ$. The resulting space is a random {\em $\R$-tree}, that is a
geodesic metric space in which any two points are joined by a unique
continuous injective path up to reparameterization. This comes from
general results on encodings of $\R$-trees by continuous functions,
see \cite{duqlegprep} for instance.  The class of $s_*$, that we still
call $\rho$, is distinguished as the root of $\TT_Z$, and any point in
this space (say encoded by the time $s\in [0,1]$) is joined to $\rho$
by a unique geodesic. A formulation of the main result of
\cite{legall08} is that this path projects into $(S,D)$ as a geodesic
$\gamma^{(s)}$ from $\rho$ to $\bp(s)$, and that {\em any} $D$-geodesic
from $\rho$ can be described in this way. In
particular, this implies that $D^\circ(x,y)=D(x,y)=Z_y-Z_x$ whenever
$x$ lies on a $D$-geodesic from $\rho$ to $y$. 

This implies the following improvement of \eqref{eq:14}.  In a metric
space $(X,d)$, we say that $(x,y,z)$ are aligned if
$d(x,y)+d(y,z)=d(x,z)$: Note that this notion of alignment depends on
the order in which $x,y,z$ are listed (in fact, of the middle term
only). In a geodesic metric space, this is equivalent to saying that
$y$ lies on a geodesic from $x$ to $z$.
\begin{lmm}\label{sec:confl-geod--1}
  Almost surely, for every $x,y\in S$, it holds that $(\rho,x,y)$ are
  aligned in $(S,D)$ if and only if they are aligned in $(S,D^*)$, and
  in this case it holds that
$$D(x,y)=D^\circ(x,y)=D^*(x,y)$$
\end{lmm}

To prove this lemma, assume that $(\rho,x,y)$ are aligned in
$(S,D)$. We already saw that this implies that $D(x,y)=D^\circ(x,y)$,
so that necessarily $D^*(x,y)=D(x,y)$ as well, since $D\leq D^*\leq
D^\circ$. By \eqref{eq:14} this implies that $(\rho,x,y)$ are aligned
in $(S,D^*)$, as well as the last conclusion of the lemma. Conversely,
if $(\rho,x,y)$ are aligned in $(S,D^*)$, then by \eqref{eq:14}, the
triangle inequality, and the fact that $D\leq D^*$,
$$D^*(\rho,y)=D(\rho,y)\leq D(\rho,x)+D(x,y)\leq
D^*(\rho,x)+D^*(x,y)=D^*(\rho,y)\, ,$$ meaning that there must be
equality throughout. Hence $(\rho,x,y)$ are aligned in $(S,D)$, and we
conclude as before.

Another important consequence of this description of geodesics is that
the geodesics $\gamma^{(s)},\gamma^{(t)}\in [0,1]$ are bound to
merge into a single path in a neighborhood of $\rho$, a phenomenon
called {\em confluence of geodesics}. This particularizes to the
following statement.

\begin{lmm}\label{sec:confl-geod-}
  Let $s,t\in [0,1]$. Then the images of $\gamma^{(s)},\gamma^{(t)}$
  coincide in the complement of $B_D(\bp(s),D^\circ(s,t))$ (or of
  $B_D(\bp(t),D^\circ(s,t))$).
\end{lmm}

To prove this lemma, it suffices to note that $D^\circ(s,t)$ is
the length of the path obtained by following $\gamma^{(s)}$ back from
its endpoint $\bp(s)$ until it coalesces with $\gamma^{(t)}$ in the
tree $\TT_Z$, and then following $\gamma^{(t)}$ up to $\bp(t)$. Note
that the same is true with $D^*$ instead of $D$ because of
Lemma \ref{sec:confl-geod--1}.

\subsection{Plan of the proof}\label{sec:plan-proof}

In this section, we decompose the proof of Theorem
\ref{sec:introduction} into several intermediate statements. The two
main ones (Propositions \ref{sec:introduction-1} and \ref{sec:th1})
will be proved in Sections \ref{sec:rough-comp-betw} and
\ref{sec:exceptional-nature-3-4}. 

The first step is to show that the distances $D$ and $D^*$ are almost
equivalent distances.

\begin{prp}\label{sec:introduction-1}
Let $\alpha\in (0,1)$ be fixed. Then almost-surely, there exists a
(random) $\eps_1>0$ such that for every $x,y\in S$ with $D(x,y)\leq
\eps_1$, one has $D^*(x,y)\leq D(x,y)^\alpha$.
\end{prp}

The second step is based on a study of particular points of the space
$(S,D)$, from which emanate {\em stars} made of several disjoint
geodesic paths, which we also call {\em arms} by analogy with the
so-called ``arm events'' of percolation. 

\begin{defn}\label{sec:geodesic-arms-1}
  Let $(X,d)$ be a geodesic metric space, and $x_1,\ldots, x_k,x$ be
  $k+1$ points in $X$. We say that $x$ is a $k$-star point with
  respect to $x_1,\ldots,x_k$ if for every geodesic paths (arms)
  $\gamma_1,\ldots,\gamma_k$ from $x$ to $x_1,\ldots,x_k$
  respectively, it holds that for every $i,j\in \{1,2,\ldots,k\}$ with
  $i\neq j$, the paths $\gamma_i,\gamma_j$ intersect only at $x$.  We
  let $\mathcal{G}(X;x_1,\ldots,x_k)$ be the set of points $x\in X$
  that are $k$-star points with respect to $x_1,\ldots,x_k$.
\end{defn}

Conditionally on $(S,D)$, let $x_1,x_2$ be random points of $S$ with
distribution $\lambda$.  These points can be constructed by
considering two independent uniform random variables $U_1,U_2$ in
$[0,1]$, independent of $(\ee,Z,D)$, and then setting $x_1=\bp(U_1)$
and $x_2=\bp(U_2)$. These random variables always exist up to
enlarging the underlying probability space if necessary.  Then by
\cite[Corollary 8.3 (i)]{legall08} (see also \cite{miertess}), with
probability $1$, there is a unique $D$-geodesic from $x_1$ to $x_2$,
which we call $\gamma$. By the same result, the geodesics from $\rho$
to $x_1$ and $x_2$ are also unique, we call them $\gamma_1$ and
$\gamma_2$.  Moreover, $\rho$ is a.s.\ not on $\gamma$, because
$\gamma_{1},\gamma_{2}$ share a common initial segment (this is the
confluence property mentioned earlier). So trivially
$D(x_1,\rho)+D(x_2,\rho)>D(x_1,x_2)$ meaning that $(x_1,\rho,x_2)$ are
not aligned.

We let
$$\Gamma=\gamma([0,D(x_1,x_2)])\cap {\cal G}(S;x_1,x_2,\rho)\, .$$
Equivalently, with probability $1$, it holds that $y\in \Gamma$ if and
only if any geodesic from $y$ to $\rho$ does not intersect $\gamma$
except at $y$ itself: Note that the a.s.\ unique geodesic from $y$ to
$x_1$ is the segment of $\gamma$ that lies between $y$ and $x_1$, for
otherwise, there would be several distinct geodesics from $x_1$ to
$x_2$.

\begin{prp}\label{sec:th1}
  There exists $\delta\in (0,1)$ for which the following property is
  satisfied almost-surely: There exists a (random) $\eps_2>0$ such
  that for every $\eps\in(0,\eps_2)$, the set $\Gamma$ can be covered
  with at most $\eps^{-(1-\delta)}$ balls of radius $\eps$ in $(S,D)$.
\end{prp}

Note that this implies a (quite weak) property of the Hausdorff
dimension of $\Gamma$.

\begin{crl}\label{sec:plan-proof-1}
  The Hausdorff dimension of $\Gamma$ in $(S,D)$ is a.s.\ strictly
  less than $1$.
\end{crl}

We do not know the exact value of this dimension. The largest constant
$\delta$ that we can obtain following the approach of this paper is
not much larger than $0.00025$, giving an upper bound of $0.99975$ for
the Hausdorff dimension of $\Gamma$.

 With Propositions \ref{sec:introduction-1}
and \ref{sec:th1} at hand, proving Theorem \ref{sec:introduction}
takes only a couple of elementary steps, which we now perform.

\begin{lmm}\label{sec:introduction-2}
Let $(s,t)$ be a non-empty subinterval of $[0,D(x_1,x_2)]$ such that
$\gamma(v)\notin \Gamma$ for every $v\in (s,t)$. Then,
there exists a unique $u\in [s,t]$ such that 
\begin{itemize}
\item $(\rho, \gamma(s),\gamma(u))$ are aligned, and
\item $(\rho, \gamma(t),\gamma(u))$ are aligned.
\end{itemize}
\end{lmm}

\noindent{\bf Proof. } Fix $v\in (s,t)$. Since $\gamma(v)\notin
\Gamma$, there exists a geodesic from $\gamma(v)$ to $\rho$ that
intersects the image of $\gamma$ at some point $\gamma(v')$, with
$v'\neq v$. In particular, the points $(\rho,\gamma(v'),\gamma(v))$
are aligned. Let us assume first that $v'<v$, and let
$$w=\inf\{v''\in [s,v]: (\rho,\gamma(v''),
\gamma(v))\mbox{ are aligned}\}\, .$$ Then $w\in [s,v)$ since $v'$ is
in the above set. Let us show that $w=s$. If if was true that $w>s$,
then $\gamma(w)$ would not be an element of $\Gamma$, and some
geodesic from $\gamma(w)$ to $\rho$ would thus have to intersect the
image of $\gamma$ at some point $\gamma(w')$ with $w'\neq w$. But we
cannot have $w'<w$ by the minimality property of $w$, and we cannot
have $w'>w$ since otherwise, both $(\rho,\gamma(w),\gamma(w'))$ on the
one hand and, $(\rho,\gamma(w'),\gamma(w))$ on the other hand, would
be aligned. This shows that $w=s$, and this implies that $(\rho,
\gamma(s),\gamma(v))$ are aligned. By a similar reasoning, if it holds
that $v'>v$, then $(\rho,\gamma(t),\gamma(v))$ are aligned.

At this point, we can thus conclude that for every $v\in (s,t)$,
either $(\rho, \gamma(s),\gamma(v))$ are aligned or
$(\rho, \gamma(t),\gamma(v))$ are aligned. The
conclusion follows by defining 
$$u=\sup\{u'\in [s,t]: (\rho, \gamma(s),\gamma(u'))
\mbox{ are aligned}\}\, ,$$ as can be easily checked.  \cq

\bigskip

Let $\delta$ be as in Proposition \ref{sec:th1}.  Let $\eps>0$ be
chosen small enough so that $\Gamma$ is covered by balls
$B_D(x_{(i)},\eps),1\leq i\leq K$, for some
$x_{(1)},x_{(2)},\ldots,x_{(K)}$ with $K=\lfloor
\eps^{-(1-\delta)}\rfloor$. Without loss of generality, we can assume
that $\{x_1,x_2\}\subset\{x_{(1)}\ldots,x_{(K)}\}$, up to increasing
$K$ by $2$, or leaving $K$ unchanged at the cost of taking smaller
$\delta$ and $\eps$. We can also assume that all the balls
$B_D(x_{(i)},\eps)$ have a non-empty intersection with the image of
$\gamma$, by discarding all the balls for which this property does not
hold. Let
$$r_i=\inf\{t\geq 0:\gamma(t)\in B_D(x_{(i)},\eps)\}\, ,\qquad
r'_i=\sup\{t\leq D(x_1,x_2):\gamma(t)\in B_D(x_{(i)},\eps)\}\, ,$$ so
that $r_i<r'_i$ for $i\in \{1,2,\ldots,K\}$, since $B_D(x_{(i)},\eps)$
is open, and let
\begin{equation}\label{eq:5}
A=\bigcup_{i=1}^K[r_i,r'_i]\, ,
\end{equation}
so that $\Gamma\subset \gamma(A)$.

\begin{lmm}\label{sec:cand-brown-map-1}
  For every $i\in \{1,2,\ldots,K\}$ and $r\in [r_i,r'_i]$, one has
  $D(\gamma(r),x_{(i)})\leq 2\eps$.
\end{lmm}

\noindent{\bf Proof. } If this was not the case, then since $\gamma$
is a $D$-geodesic path that passes through $\gamma(r_i),\gamma(r)$ and
$\gamma(r_i')$ in this order, we would have
$D(\gamma(r_i),\gamma(r_i'))=
D(\gamma(r_i),\gamma(r))+D(\gamma(r),\gamma(r_i'))> 2\eps$. But on
the other hand,
$$D(\gamma(r_i),\gamma(r_i'))\leq
D(\gamma(r_i),x_{(i)})+D(x_{(i)},\gamma(r_i'))<2\eps\, ,$$ a
contradiction. \cq

\bigskip

Let $I$ be the set of indices $j\in \{1,2,\ldots,K\}$ such that
$[r_j,r'_j]$ is maximal for the inclusion order in the family
$\{[r_j,r'_j],1\leq j\leq K\}$. Note that the set $A$ of (\ref{eq:5})
is also equal to $A=\bigcup_{i\in I}[r_i,r'_i]$.  This set can be
uniquely rewritten in the form $A=\bigcup_{i=0}^{K'-1}[t_i,s_{i+1}]$
where $K'\leq \#I$, and
$$0=t_0<s_1<t_1<s_2<\ldots<s_{K'-1}<t_{K'-1}<s_{K'}=D(x_1,x_2)\, .$$ Here,
the fact that $t_0=0$ and $s_{K'}=D(x_1,x_2)$ comes from the
assumption that $x_1,x_2$ are in the set $\{x_{(i)},1\leq i\leq K\}$.
We let $x_{(i)}=\gamma(s_i)$ for $1\leq i\leq K'$ and
$y_{(i)}=\gamma(t_i)$ for $0\leq i\leq K'-1$.

\begin{lmm}\label{sec:cand-brown-map}
Almost-surely, it holds that for every $\eps$ small enough, 
$$\sum_{i=0}^{K'-1}D^*(y_{(i)},x_{(i+1)})\leq 4K\eps^{(2-\delta)/2}\, .$$
\end{lmm}

\noindent{\bf Proof. } Consider one of the connected components
$[t_i,s_{i+1}]$ of $A$. Let 
$$J_i=\{j\in I:[r_j,r'_j]\subset [t_i,s_{i+1}]\}\, ,$$ so that
$\sum_{0\leq i\leq K'-1}\#J_i=\#I\leq K$. The lemma will thus follow
if we can show that $D^*(y_{(i)},x_{(i+1)})\leq 2\#J_i\eps^{(2-\delta)/2}$ for
every $\eps$ small enough. By the definition of $I$, it holds that if
$[r_j,r_j']$ and $[r_k,r'_k]$ with $j,k\in I$ have non-empty
intersection, then these two intervals necessarily overlap, meaning
that $r_j\leq r_k\leq r'_j\leq r'_k$ or vice-versa. Let us re-order
the $r_j,j\in J_i$ as $r_{j_k},1\leq k\leq \#J_i$ in non-decreasing
order. Then $\gamma(r_{j_1})=y_{(i)}$ and
$\gamma(r_{j_{\#J_i}}')=x_{(i+1)}$, and
\begin{equation}\label{eq:4}
D^*(y_{(i)},x_{(i+1)})\leq \sum_{k=1}^{\#J_i-1}D^*(\gamma(r_{j_k}),
\gamma(r_{j_{k+1}}))+D^*(\gamma(r_{j_{\#J_i}}),x_{(i+1)})\, .
\end{equation}
By Lemma \ref{sec:cand-brown-map-1} and the overlapping property of
the intervals $[r_j,r'_j]$, we have that
$$D(\gamma(r_{j_k}),\gamma(r_{j_{k+1}}))\leq 4\eps\, ,\qquad
D(\gamma(r_{j_{\#J_i}}),x_{(i+1)})\leq 4\eps\, .$$ Now apply
Proposition \ref{sec:introduction-1} with $\alpha=(2-\delta)/2$ to
obtain from (\ref{eq:4}) that a.s.\ for every $\eps$ small enough,
$$D^*(y_{(i)},x_{(i+1)})\leq \#J_i(4\eps)^{(2-\delta)/2}\leq
4\#J_i\eps^{(2-\delta)/2}\, ,$$ which concludes the proof. \cq

\bigskip

We can now finish the proof of Theorem \ref{sec:introduction}. The
complement of the set $A$ in $[0,D(x_1,x_2)]$ is the union of the
intervals $(s_i,t_i)$ for $1\leq i\leq K'-1$. Now, for every such $i$,
the image of the interval $(s_i,t_i)$ by $\gamma$ does not intersect
$\Gamma$, since $\Gamma \subset\gamma(A)$ and $\gamma$ is injective.
By Lemma \ref{sec:introduction-2}, for every $i\in
\{1,2,\ldots,K'-1\}$, we can find $u_i\in [s_i,t_i]$ such that $(\rho,
\gamma(s_i),\gamma(u_i))$ are aligned, as well as
$(\rho,\gamma(t_i),\gamma(u_i))$.  Letting
$x_{(i)}=\gamma(s_i),y_{(i)}=\gamma(t_i)$ and $z_{(i)}=\gamma(u_i)$,
by Lemma \ref{sec:confl-geod--1},
$$D^*(x_{(i)},z_{(i)})=D(x_{(i)},z_{(i)})\, ,\qquad D^*(y_{(i)},z_{(i)})=D(y_{(i)},z_{(i)})\, ,\qquad
1\leq i\leq K'-1\, .$$ By the triangle inequality, and the fact that
$\gamma$ is a $D$-geodesic, we have a.s.\ 
\begin{eqnarray*}
D^*(x_1,x_2)&\leq &\sum_{i=1}^{K'-1}(D^*(x_{(i)},z_{(i)})+D^*(z_{(i)},y_{(i)}))+
\sum_{i=0}^{K'-1}D^*(y_{(i)},x_{(i+1)})\\ &= &
\sum_{i=1}^{K'-1}(D(x_{(i)},z_{(i)})+D(z_{(i)},y_{(i)}))+\sum_{i=0}^{K'-1}D^*(y_{(i)},x_{(i+1)})\\ &\leq
&D(x_1,x_2)+\sum_{i=0}^{K'-1}D^*(y_{(i)},x_{(i+1)})\\ &\leq
&D(x_1,x_2)+4K\eps^{(2-\delta)/2}\, ,
\end{eqnarray*}
where we used Lemma \ref{sec:cand-brown-map} at the last step,
assuming $\eps$ is small enough. Since $K\leq \eps^{-(1-\delta)}$,
this is enough to get $D^*(x_1,x_2)\leq D(x_1,x_2)$ by letting
$\eps\to 0$. Since $D\leq D^*$ we get $D(x_1,x_2)=D^*(x_1,x_2)$.

Note that the previous statement holds a.s.\ for
$\lambda\otimes\lambda$-almost every $x_1,x_2\in S$. This means that
if $x_1,x_2,\ldots$ is an i.i.d.\ sequence of $\lambda$-distributed
random variables (this can be achieved by taking an i.i.d.\ sequence
of uniform random variables on $[0,1]$, independent of the Brownian
map, and taking their images under $\bp$), then almost-surely one has
$D^*(x_i,x_j)=D(x_i,x_j)$ for every $i,j\geq 1$. But since the set
$\{x_i,i\geq 1\}$ is a.s.\ dense in $(S,D)$ (or in $(S,D^*)$), the
measure $\lambda$ being of full support, we obtain by a density
argument that a.s., for every $x,y\in S$, $D^*(x,y)=D(x,y)$. This ends
the proof of Theorem \ref{sec:introduction}.

\section{Rough comparison between $D$ and $D^*$}\label{sec:rough-comp-betw}

The goal of this Section is to prove Proposition
\ref{sec:introduction-1}. We start with an elementary statement in
metric spaces.

\begin{lmm}\label{sec:rough-comp-betw-1}
Let $(X,d)$ be an arcwise connected metric space, and let $x,y$ be two
distinct points in $X$. Let $\gamma$ be a continuous path with
extremities $x$ and $y$. Then for every $\eta> 0$, there exist at
least $K= \lfloor d(x,y)/(2\eta)\rfloor+1$ points $y_1,\ldots,y_K$ in
the image of $\gamma$ such that $d(y_i,y_j)\geq 2\eta$ for every
$i,j\in \{1,2,\ldots,K\}$ with $i\neq j$.
\end{lmm}

\noindent{\bf Proof. } Let us assume without loss of generality that
$\gamma$ is parameterized by $[0,1]$ and that
$\gamma(0)=x,\gamma(1)=y$. Also, assume that $d(x,y)\geq 2\eta$, since
the statement is trivial otherwise. In this case we have $K\geq 2$.

  Set $s_0=0$, and by induction, let
$$s_{i+1}=\sup\{t\leq 1:d(\gamma(t),\gamma(s_i))<2\eta\}\, ,\qquad
  i\geq 0\, .$$ The sequence $(s_i,i\geq 0)$ is non-decreasing, and it
  holds that $d(\gamma(s_i),\gamma(s_{i+1}))\leq 2\eta$ for every
  $i\geq 0$.  Let $i\in \{0,1,\ldots,K-2\}$. Then, since
  $x=\gamma(0)=\gamma(s_0)$,
$$d(x,\gamma(s_i))\leq \sum_{j=0}^{i-1}d(\gamma(s_j),\gamma(s_{j+1}))
\leq 2\eta i \leq 2\eta(K-2)\leq d(x,y)-2\eta\, ,$$ which implies, by
the triangle inequality,
$$d(\gamma(s_i),y)\geq d(x,y)-d(x,\gamma(s_i))\geq
d(x,y)-(d(x,y)-2\eta)=2\eta\, .$$ From this, and the definition of
$(s_i,i\geq 0)$ it follows that $d(\gamma(s_i),\gamma(s_j))\geq 2\eta$
for every $i\in \{0,1,\ldots,K-2\}$ and $j>i$. This yields the wanted
result, setting $y_i=\gamma(s_{i-1})$ for $1\leq i\leq K$.  \cq

\bigskip

Our next tool is a uniform estimate for the volume of $D$-balls in
$S$.

\begin{lmm}\label{sec:rough-comp-betw-2}
Let $\eta\in (0,1)$ be fixed. Then almost-surely, there exists a
(random) $C\in (0,\infty)$ such that for every $r\geq 0$ and every
$x\in S$, one has
$$\lambda(B_D(x,r))\leq Cr^{4-\eta} \, .$$
\end{lmm}

This is an immediate consequence of a result by Le Gall
\cite[Corollary 6.2]{legall08}, who proves the stronger fact that the
optimal random ``constant'' $C$ of the statement has moments of all
orders.  In the remainder of this section, we will use let $C,c$ for
such almost-surely finite positive random variables. As long as no
extra property besides the a.s.\ finiteness of these random variables
is required, we keep on calling them $C,c$ even though they might
differ from statement to statement or from line to line, just as if
they were universal constants.

The following statement is a rough uniform lower estimate for
$D^*$-balls.

\begin{lmm}\label{sec:rough-comp-betw-3}
Let $\eta\in (0,1)$ be fixed. Then almost-surely, there exists a
(random) $c\in (0,\infty)$ and $r_0>0$ such that for every $r\in
[0,r_0] $ and every $x\in S$, one has
$$\lambda(B_{D^*}(x,r))\geq c r^{4+\eta} \, .$$
\end{lmm}

\noindent{\bf Proof. } We use the fact that $B_{D^\circ}(x,r)
\subseteq B_{D^*}(x,r)$ for every $x\in S$ and $r\geq 0$, where
$B_{D^\circ}(x,r)=\{y\in S:D^\circ(x,y)<r\}$. We recall that $D^\circ$
does not satisfy the triangle inequality, which requires a little
extra care when manipulating ``balls'' of the form $B_{D^\circ}(x,r)$.

By definition, $D^\circ(x,y)=\inf_{s,t}D^\circ(s,t)$ where the infimum
is over $s,t\in [0,1]$ such that $\bp(s)=x,\bp(t)=y$. Consequently, for
every $s\in [0,1]$ such that $\bp(s)=x$, 
$$\bp(\{t\in [0,1]:D^\circ(s,t)<r\})\subseteq B_{D^\circ}(x,r)\, ,$$
which implies, by definition of $\lambda$, 
$$
\lambda(B_{D^*}(x,r))\geq \lambda(B_{D^\circ}(x,r)) \geq {\rm
  Leb}(\{t\in [0,1]:D^\circ(s,t)<r\})\, .
$$ Consequently, for every $r>0$,
\begin{equation}\label{eq:1}
\lambda(B_{D^*}(x,r))\geq {\rm Leb}(\{t\in [0,1]:D^\circ(s,t)\leq
r/2\})\, .
\end{equation}
We use the fact that $Z$ is a.s.\ Hölder-continuous with exponent
$1/(4+\eta)$, which implies that a.s.\ there exists a random $C\in
(0,\infty)$, such that for every $h\geq 0$, $\omega(Z,h)\leq
Ch^{1/(4+\eta)}$, where
$$\omega(Z,h)=\sup \{|Z_t-Z_s|:s,t\in [0,1], |t-s|\leq h\}$$ is the
oscillation of $Z$. Since
$$D^\circ(s,t)\leq Z_s+Z_t-2\min_{u\in [s\wedge t,s\vee t]}Z_u\leq
2\omega(Z,|t-s|)\, ,$$ we obtain that for every $s\in [0,1]$, for
every $h>0$ and $t\in [(s-h)\vee 0,(s+h)\wedge 1]$,
$$D^\circ(s,t)\leq 2Ch^{1/(4+\eta)}\, .$$ Letting $h=(r
/(4C))^{4+\eta}$ yields
$${\rm Leb}(\{t\in [0,1]:D^\circ(s,t)\leq r/2\})\geq (2h)\wedge 1\,
.$$ This holds for every $s\in[0,1],r\geq 0$ so that, by (\ref{eq:1}),
$$\lambda(B_{D^*}(x,r))\geq
\Big(\frac{2}{(4C)^{4+\eta}}r^{4+\eta}\Big)\wedge 1\, .$$ This
yields the wanted result with $r_0=4C/2^{1/(4+\eta)}$ and
$c=2/(4C)^{4+\eta}$.  \cq

\bigskip

We are now able to prove Proposition \ref{sec:introduction-1}, and we
argue by contradiction. Let $\alpha\in (0,1)$ be fixed, and assume
that the statement of the proposition does not hold. This implies that
with positive probability, one can find two sequences $(x_n,n\geq 0)$
and $(y_n,n\geq 0)$ of points in $S$ such that $D(x_n,y_n)$ converges
to $0$, and $D^*(x_n,y_n)>D(x_n,y_n)^\alpha$ for every $n\geq 0$. From
now on until the end of the proof, we restrict ourselves to this event
of positive probability, and almost-sure statements will implicitly be
in restriction to this event.

Let $\gamma_n$ be a $D$-geodesic path from $x_n$ to $y_n$. Let
$V_\beta^D(\gamma_n)$ be the $D,\beta$-thikening of the image of
$\gamma_n$:
$$V_\beta^D(\gamma_n)=\{x\in S:\exists t\in [0,D(x_n,y_n)],
D(\gamma_n(t),x)<\beta\}\, .$$
Then $V_\beta^D(\gamma_n)$ is contained in a union of at most $\lfloor
D(x_n,y_n)/(2\beta)\rfloor+1$ $D$-balls of radius $2\beta$: Simply take
centers $y_n$ and $\gamma_n(2\beta k)$ for $1\leq k\leq \lfloor
D(x_n,y_n)/(2\beta)\rfloor$. Consequently, for every $\beta>0$, 
$$\lambda(V_\beta^D(\gamma_n))\leq
\Big(\frac{D(x_n,y_n)}{2\beta}+1\Big)\sup_{x\in
  S}\lambda(B_D(x,2\beta))\, .$$ By applying Lemma
\ref{sec:rough-comp-betw-2}, for any $\eta\in (0,1)$, whose value
will be fixed later on, we obtain a.s.\ the existence of $C\in
(0,\infty)$ such that for every $n\geq 0$ and $\beta\in
[0,D(x_n,y_n)]$,
\begin{equation}\label{eq:2}
\lambda(V_\beta^D(\gamma_n))\leq C\beta^{3-\eta}D(x_n,y_n)\, .
\end{equation}

Let $V_\beta^{D^*}(\gamma_n)$ be the $D^*,\beta$-thikening of
$\gamma_n$, defined as $V_\beta^D(\gamma_n)$ above but with $D^*$
instead of $D$.  The spaces $(S,D^*)$ and $(S,D)$ being homeomorphic,
we obtain that $(S,D^*)$ is arcwise connected and $\gamma_n$ is a
continuous path in this space. Therefore Lemma
\ref{sec:rough-comp-betw} applies: For every $\beta>0$ we can find
points $y_1,\ldots,y_K$, $K=\lfloor D^*(x_n,y_n)/(2\beta)\rfloor+1$,
such that $D^*(y_i,y_j)\geq 2\beta$ for every $i,j\in
\{1,2,\ldots,K\}$ with $i\neq j$. From this, it follows that the balls
$B_{D^*}(y_i,\beta),1\leq i\leq K$ are pairwise disjoint, and they are
all included in $V_\beta^{D^*}(\gamma_n)$. Hence,
$$\lambda(V_\beta^{D^*}(\gamma_n))\geq
\sum_{i=1}^K\lambda(B_{D^*}(y_i,\beta))\geq K\inf_{x\in
  S}\lambda(B_{D^*}(x,\beta))\geq \frac{D^*(x_n,y_n)}{2\beta}\inf_{x\in
  S}\lambda(B_{D^*}(x,\beta))\, .$$ For the same $\eta\in (0,1)$ as
before, by using Lemma \ref{sec:rough-comp-betw-3} and the definition
of $x_n,y_n$, we conclude that a.s.\ there exists $c\in (0,\infty)$
and $r_0>0$ such that for every $\beta\in[0,r_0]$,
\begin{equation}\label{eq:3}
\lambda(V_\beta^{D^*}(\gamma_n))\geq c\beta^{3+\eta}D^*(x_n,y_n)
\geq c\beta^{3+\eta}D(x_n,y_n)^\alpha\, .
\end{equation}
But since $D\leq D^*$, we have that $V_\beta^{D^*}(\gamma_n)\subseteq
V_\beta^{D}(\gamma_n)$, so that (\ref{eq:2}) and (\ref{eq:3}) entail
that for every $\beta\in [0, D(x_n,y_n)\wedge r_0 ]$,
$$\beta^{2\eta}\leq CD(x_n,y_n)^{1-\alpha}\, ,$$ for some
a.s.\ finite $C>0$. But taking $\eta=(1-\alpha)/4$ and then
$\beta=D(x_n,y_n)\wedge r_0$, we obtain since $D(x_n,y_n)\to 0$ that
$D(x_n,y_n)^{(1-\alpha)/2}=O(D(x_n,y_n)^{1-\alpha})$, which is a
contradiction. This ends the proof of
Proposition \ref{sec:introduction-1}. \cq

\section{Covering $3$-star points on typical
  geodesics}\label{sec:exceptional-nature-3-4}

We now embark in our main task, which is to prove Proposition
\ref{sec:th1}. 

\subsection{Entropy number estimates}\label{sec:entr-numb-estim}

We use the same notation as in Section \ref{sec:plan-proof}. In this
section, we are going to fix two small parameters $\delta,\beta\in
(0,1/2)$, which will be tuned later on: The final value of $\delta$
will be the one that appears in Proposition \ref{sec:th1}.

We want to estimate the number of $D$-balls of radius $\eps$ needed to
cover the set $\Gamma$.  Since we are interested in bounding this
number by $\eps^{-(1-\delta)}$ we can consider only the points of
$\Gamma$ that lie at distance at least $8\eps^{1-\beta}$ from
$x_1,x_2$ and $\rho$. Since $\Gamma$ is included in the image of
$\gamma$, the remaining part of $\Gamma$ can certainly be covered with
at most $32\eps^{-\beta}$ balls of radius $\eps$. Since we chose
$\beta,\delta<1/2$, we have $\beta<1-\delta$, so this extra number of
balls will be negligible compared to $\eps^{-(1-\delta)}$.

So for $\eps>0$, let $\n_\Gamma(\eps)$ be the minimal $n\geq
1$ such that there exist $n$ points $x_{(1)},\ldots,x_{(n)}\in S$ such
that
$$\Gamma\setminus
(B_D(x_1,8\eps^{1-\beta})\cup B_D(x_2,8\eps^{1-\beta})\cup
B_D(\rho,8\eps^{1-\beta}))\subset\bigcup_{i=1}^{n} B_D(x_{(i)},\eps)\,
.$$ We call $\Gamma_\eps$ the set on the left-hand side.  We first
need a simple control on $D^\circ$.

\begin{lmm}\label{sec:exceptional-nature-3-6}
  Fix $\delta>0$.  Almost surely, there exists a (random)
  $\eps_3(\delta)\in(0,1)$ such that for every $t\in [0,1]$,
  $\eps\in(0,\eps_3(\delta))$ and $s\in [(t-\eps^{4+\delta})\vee
  0,(t+\eps^{4+\delta})\wedge 1]$, it holds that $D^\circ(s,t)\leq
  \eps/2$.
\end{lmm}

\proof This is an elementary consequence of the fact that $Z$ is a.s.\
Hölder-continuous with any exponent $\alpha\in (0,1/4)$, and the
definition of $D^\circ$. From this, we get that $D^\circ(s,t)\leq
2K_\alpha|t-s|^\alpha$, where $K_\alpha$ is the $\alpha$-Hölder norm
of $Z$. If $|t-s|\leq \eps^{4+\delta}$, then $D^\circ(s,t)\leq
2K_{\alpha}\eps^{(4+\delta)\alpha}$, so it suffices to choose $\alpha\in
(1/(4+\delta),1/4)$, and then
$\eps_3(\delta)=(4K_\alpha)^{-((4+\delta)\alpha-1)^{-1}}\wedge (1/2)$. \cq

\medskip

For every $y\in S$, let $t\in [0,1]$ such that $\bp(t)=y$. We
let 
$$F(y,\eps)=\bp([(t-\eps^{4+\delta})\vee 0,(t+\eps^{4+\delta})\wedge 1])\, ,\qquad
\eps>0\, .$$ Note that in general, $F(y,\eps)$ does depend on the
choice of $t$, so we let this choice be arbitrary, for instance, $t$
can be the smallest possible in $[0,1]$. Lemma
\ref{sec:exceptional-nature-3-6} and the fact that $D\leq D^\circ$
entails that $F(y,\eps)$ is included in $B_D(y,\eps/2)$ for every
$\eps\in (0,\eps_3(\delta))$, but $F(y,\eps)$ is not necessarily a
neighborhood of $y$ in $S$. We can imagine it as having the shape of a
fan with apex at $y$.

\begin{lmm}\label{sec:exceptional-nature-3-7}
  Let $\delta>0$ and $\eps\in (0,\eps_3(\delta))$ be as in Lemma
  \ref{sec:exceptional-nature-3-6}. Let
  $t_{(1)},t_{(2)},\ldots,t_{(N)}$ be elements in $[0,1]$ such that
  the intervals $[(t_{(i)}-\eps^{4+\delta}/2)\vee
  0,(t_{(i)}+\eps^{4+\delta}/2)\wedge 1],1\leq i\leq N$ cover
  $[0,1]$. Then, letting $x_{(i)}=\bp(t_{(i)})$,
$$\n_\Gamma(\eps)\leq \sum_{i=1}^{N}\ind_{\{x_{(i)}\in \bigcup_{y\in
    \Gamma_\eps}F(y,\eps)\}}\, .$$
\end{lmm}

\proof Let $y\in \Gamma_\eps$ be given. By assumption, for every $t\in
[0,1]$, the interval $[(t-\eps^{4+\delta})\vee 0,(t+\eps^{4+\delta})\wedge 1]$ contains
at least one point $t_{(i)}$. By definition, this means that
$F(y,\eps)$ contains at least one of the points $x_{(i)},1\leq i\leq
N$. By the assumption on $\eps$, we obtain that $y\in
B_D(x_{(i)},\eps/2)$. Therefore, the union of balls
$B_D(x_{(i)},\eps)$, where $x_{(i)}$ is contained in
$F(y,\eps)$ for some $y\in \Gamma_\eps$, cover $\Gamma_\eps$.  \cq

\medskip

Let $(U_0^{(i)},i\geq 1)$
be an i.i.d.\ sequence of uniform random variables on $[0,1]$,
independent of $(\ee,Z,D)$ and of $U_1,U_2$. We let
$x_0^{(i)}=\bp(U_0^{(i)})$.  For $\eps>0$ let $N_\eps=\lfloor
\eps^{-4-2\delta}\rfloor$. The probability of the event $\mathcal{B}_\eps$
that the intervals $[(U_0^{(i)}-\eps^{4+\delta}/2)\vee
0,(U_0^{(i)}+\eps^{4+\delta}/2)\wedge 1],1\leq i\leq N_\eps$ do not cover
$[0,1]$ is less than the probability that there exists a $j\leq
2\eps^{-4-\delta}$ such that $[j\eps^{4+\delta}/2,((j+1)\eps^{4+\delta}/2)\wedge 1]$
does not contain any of the $U_0^{(i)},1\leq i\leq N_\eps$. This has a
probability at most $2\eps^{-4-\delta}e^{-\eps^{-\delta}/2}$, which decays
faster than any power of $\eps$ as $\eps\to 0$. Using Lemma
\ref{sec:exceptional-nature-3-7}, we get the existence of a finite
constant $C$ such that for every $\eps>0$,
\begin{eqnarray*}
\lefteqn{P(\n_\Gamma(\eps)\geq \eps^{-(1-\delta)}, \eps\leq
  \eps_3(\delta) )}\\
 &\leq &
P(\mathcal{B}_\eps)+
P\bigg(\sum_{i=1}^{N_\eps}\ind_{\{x_0^{(i)}\in \bigcup_{y\in
    \Gamma}F(y,\eps)\}}\geq \eps^{-(1-\delta)}, \eps\leq \eps_3(\delta)
\bigg)\\
&\leq & C\eps^4+\eps^{1-\delta}\eps^{-4-2\delta}P\bigg(x_0\in
\bigcup_{y\in \Gamma_\eps}F(y,\eps),\eps\leq \eps_3(\delta)\bigg)\, ,
\end{eqnarray*}
where we used the Markov inequality in the second step, and let
$x_0=x_0^{(1)}$. 

The uniqueness properties of geodesics in the Brownian map already
mentioned before imply that a.s.\ there is a unique geodesic
$\gamma_0$ from $x_0$ to $\rho$, and a unique geodesic between $x_i$
and $x_j$ for $i,j\in \{0,1,2\}$. We let $\gamma_{01},\gamma_{02}$ be
the geodesics from $x_0$ to $x_1$ and $x_2$, and recall that $\gamma$
is the unique geodesic between $x_1$ and $x_2$.

Since $F(y,\eps)\subseteq B_D(y,\eps/2)$ whenever $\eps\leq
\eps_3(\delta)$, and $\Gamma_\eps\subset \gamma$, we get that if $x_0\in
\bigcup_{y\in \Gamma_\eps}F(y,\eps)$ then $B_D(x_0,\eps/2)$ intersects
$\gamma$. Moreover, recall from Lemma \ref{sec:confl-geod-} that
$D^\circ$ is a measure of how quickly two geodesics in the Brownian
map coalesce. By definition, we have $D^\circ(x_0,y)\leq \eps/2$
whenever $x_0\in F(y,\eps)$, so that outside of $B_D(x_0,\eps/2)$,
the image of $\gamma_0$ is included in some geodesic $\gamma'$ going
from $y$ to $\rho$. Since $y\in \Gamma$, by definition, $\gamma'$ does
not intersect $\gamma_{12}$ except at the point $y$ itself.

We conclude that the event $\A_0(\eps)=\{x_0\in \bigcup_{y\in
  \Gamma_\eps}F(y,\eps)\mbox{ and }\eps\in (0,\eps_3(\delta))\}$
implies that either of the events $\A_1(\eps,\beta)$ or
$\A_2(2\eps^{1-\beta})$ occurs, where
\begin{itemize}
\item $\A_1(\eps,\beta)$ is the event that $D(x_0,x_1)\wedge
  D(x_0,x_2)\geq 7\eps^{1-\beta}$, that $\gamma$ intersects
  $B_D(x_0,\eps/2)$, and there exists a point of $\gamma_{01}\cup
  \gamma_{02}$ out of $B_D(x_0,2\eps^{(1-\beta)})$ not belonging to
  $\gamma$,
\item $\A_2(\eps)$ is the event that $D(x_0,x_1)\wedge
  D(x_0,x_2)\wedge D(x_0,\rho)\geq 7\eps/2$, that $\gamma$ intersects
  $B_D(x_0,\eps)$, and that the geodesics
  $\gamma_0,\gamma_{01},\gamma_{02}$ do not intersect outside of
  $B_D(x_0,\eps)$.
\end{itemize}
Indeed, on $\A_0(\eps)$, we first note that $D(x_0,\Gamma_\eps)
<\eps$, and since 
$$D(x_1,\Gamma_\eps)\wedge D(x_2,\Gamma_\eps)\wedge D(\rho,\Gamma_\eps)\geq
8\eps^{1-\beta}$$ we obtain 
$$D(x_0,x_1)\wedge D(x_0,x_2)\wedge D(x_0,\rho)\geq
7\eps^{1-\beta}=(7/2)\cdot 2\eps^{1-\beta}\, .$$ Then on
$\A_0(\eps)\setminus \A_1(\eps,\beta)$, outside of
$B_D(x_0,2\eps^{1-\beta})$, the geodesics $\gamma_{01}$ and
$\gamma_{02}$ are included in $\gamma$, and $\gamma_0$ is included in
the geodesic $\gamma'$ discussed above, so $\A_2(2\eps^{1-\beta})$
occurs.

The key lemmas give an estimation of the probabilities for
these events. 

\begin{lmm}\label{sec:covering-3-star}
  For every $\beta\in(0,1)$, there exists $C\in
  (0,\infty)$ such that for every $\eps>0$,
$$P(\A_1(\eps,\beta))\leq C\eps^{3+\beta}\, .$$
\end{lmm}

\begin{lmm}\label{sec:covering-3-star-1}
  There exist $\chi\in(0,1)$ and $C\in (0,\infty)$ such that for every
  $\eps>0$,
$$P(\A_2(\eps))\leq C\eps^{3+\chi}\, .$$
\end{lmm}

Taking these lemmas for granted, we can now conclude the proof of
Proposition \ref{sec:th1}. The constant $\chi$ of Lemma
\ref{sec:covering-3-star-1} will finally allow us to tune the
parameters $\delta,\beta$, and we choose first $\beta$ so that
$(1-\beta)(3+\chi)>3$, and then $\delta>0$ so that
$$\delta < \frac{\beta}{3} \wedge
\frac{(1-\beta)(3+\chi)-3}{3}\, .$$ From our discussion, we have, for
$\eps_3>0$ fixed and $\eps\in (0,\eps_3)$
\begin{eqnarray}
  \lefteqn{P(\n_\Gamma(\eps)\geq \eps^{-(1-\delta)},\eps_3(\delta)\geq
    \eps_3)}\nonumber\\
  & \leq &
  C\eps^4+\eps^{-3-3\delta}(P(\A_1(\eps,\beta))+
  P(\A_2(2\eps^{1-\beta}))\nonumber\\
  &\leq &
  C(\eps^4+\eps^{\beta-3\delta}+\eps^{(1-\beta)(3+\chi)-3-3\delta})\,
  .\label{eq:6}
\end{eqnarray}
By our choice of $\delta,\beta$, the exponents in \eqref{eq:6} are
strictly positive. This gives the existence of some $\psi>0$ such that
for every $\eps_3>0$, there exists a $C>0$ such that $\eps\leq \eps_3$
implies
$$P(\n_\Gamma(\eps)\geq \eps^{-(1-\delta)},\eps_3(\delta)\geq
\eps_3)\leq C\eps^\psi\, .$$ Applying this first to $\eps$ of the form
$2^{-k},k\geq 0$ and using the Borel-Cantelli Lemma and the
monotonicity of $\n_\Gamma(\eps)$, we see that a.s., on the event
$\{\eps_3(\delta)\geq \eps_3\}$ for every $\eps>0$ small enough,
$\n_\Gamma(\eps)<\eps^{-(1-\delta)}$. Since $\eps_3(\delta)>0$ a.s.,
we obtain the same result without the condition $\eps_3(\delta)\geq
\eps_3$. This proves Proposition \ref{sec:th1}, and it remains to
prove Lemmas \ref{sec:covering-3-star} and
\ref{sec:covering-3-star-1}.

\subsection{Back to geodesic stars in discrete
  maps}\label{sec:back-geodesic-stars}

Our strategy to prove Lemmas \ref{sec:covering-3-star} and
\ref{sec:covering-3-star-1} is to relate them back to asymptotic
properties of random quadrangulations. In turn, these properties can
be obtained by using a bijection between quadrangulations and certain
maps with a simpler structure, for which the scaling limits can be
derived (and do not depend on the subsequence $(n_k)$ used to define
the space $(S,D)$).

We start by reformulating slightly the statements of Lemmas
\ref{sec:covering-3-star} and \ref{sec:covering-3-star-1}, in a way
that is more symmetric in the points $\rho,x_1,x_2,x_0$ that are
involved. For this, we use the invariance under re-rooting of the
Brownian map \cite[Theorem 8.1]{legall08} stating that $\rho$ has the
same role as a uniformly chosen point in $S$ according to the
distribution $\lambda$. So we let $x_0,x_1,x_2,x_3,x_4,\ldots$ be a
sequence of independent such points (from now on, $x_3$ will perform
the role of $\rho$, which will never be mentioned again), and let
$\gamma_{ij}$ be the a.s.\ unique geodesic from $x_i$ to $x_j$, for
$i<j$.

Both events $\A_1(\eps,\beta),\A_2(\eps)$ that are involved in Lemmas
\ref{sec:covering-3-star} and \ref{sec:covering-3-star-1} deal with
properties of ``geodesic $\eps$-stars'' in random maps, in which the
different arms of the geodesic stars separate quickly, say after a
distance at most $\eps$, rather than being necessarily pairwise
disjoint. Indeed, the event $\A_1(\eps,\beta)$ states in particular
that the random point $x_0$ lies at distance $\eps$ from a certain
point $y$ of $\gamma_{12}$, and this point can be seen as a $2$-star
point from which emanate the segments of $\gamma_{12}$ from $y$ to
$x_1,x_2$. This does imply that the geodesics
$\gamma_{01},\gamma_{02}$ are disjoint outside of the ball
$B_D(x_0,\eps)$, as is easily checked. The similar property for the
geodesics $\gamma_{03},\gamma_{01},\gamma_{02}$ under the event
$\A_2(\eps)$ is part of the definition of the latter. Therefore, we
need to estimate the probability of events of the following form:
\begin{equation}\label{eq:7}
\mathcal{G}(\eps,k)=\{\forall\, i,j\in \{1,2,\ldots,k\}
\mbox{ with }i<j, \, 
\gamma_{0i} \mbox{ is disjoint from }\gamma_{0j}\mbox{ outside
}B_D(x_0,\eps)\}\, .
\end{equation}
More precisely, we define discrete analogs for the events
$\A_1(\eps,\beta)$ and $\A_2(\eps)$. Let 
$\mathcal{A}_1^{(n)}(\eps,\beta)$ be the event that
  \begin{itemize}
  \item any
  geodesic chain from $v_1$ to $v_2$ in $Q_n$ intersects
  $B_{d_{Q_n}}(v_0,\eps (8n/9)^{1/4})$, 
\item it either holds that any geodesic chain from $v_0$ to $v_1$,
  visits at least one vertex $v$ with $d_{Q_n}(v,v_0)>\eps^{1-\beta}
  (8n/9)^{1/4}$, such that $(v_1,v,v_2)$ are not aligned, {\bf or that} any
  geodesic chain from $v_0$ to $v_2$, visits at least one vertex $v$
  with $d_{Q_n}(v,v_0)>\eps^{1-\beta} (8n/9)^{1/4}$, such that
  $(v_1,v,v_2)$ are not aligned
\item the vertices $v_0,v_1,v_2$, taken in any order, are not aligned
  in $Q_n$, and $d_{Q_n}(v_0,v_1)\wedge d_{Q_n}(v_0,v_2)\geq
  3\eps^{1-\beta}(8n/9)^{1/4}$
  \end{itemize}
  and
$\mathcal{A}_2^{(n)}(\eps)$ be the event that
  \begin{itemize}
  \item any geodesic chain from $v_1$ to $v_2$ in $Q_n$ intersects
    $B_{d_{Q_n}}(v_0,\eps (8n/9)^{1/4})$, 
  \item no two geodesic chains respectively from $v_0$ to $v_i$ and
    from $v_0$ to $v_j$ share a common vertex outside
    $B_{d_{Q_n}}(v_0,\eps(8n/9)^{1/4})$, for every $i\neq j$ in
    $\{1,2,3\}$.
  \item any three vertices among $v_0,v_1,v_2,v_3$, taken in any
    order, are not aligned in $Q_n$, and $d_{Q_n}(v_0,v_1)\wedge
    d_{Q_n}(v_0,v_2)\wedge d_{Q_n}(v_0,v_3)\geq 3\eps(8n/9)^{1/4}$
  \end{itemize}

\begin{prp}\label{sec:geod-stars-discr}
Let $Q_n$ be a uniform quadrangulation in $\bQ_n$, and conditionally
given $Q_n$, let $v_0,v_1,v_2,v_3$ be uniformly chosen points in
$V(Q_n)$.  Then 
$$P(\mathcal{A}_1(\eps,\beta))\leq
\limsup_{n\to\infty}P(\mathcal{A}^{(n)}_1(\eps,\beta))\, ,$$
and
$$P(\mathcal{A}_2(\eps))\leq
\limsup_{n\to\infty}P(\mathcal{A}^{(n)}_2(\eps))\, .$$
\end{prp}

\proof We rely on results of \cite{legall08}, see also
\cite{miertess}, stating that the marked quadrangulations
$(V(Q_n),(8n/9)^{-1/4}d_{Q_n},(v_0,v_1,\ldots,v_k))$ converge in
distribution along $(n_k)$ for the so-called $(k+1)$-pointed
Gromov-Hausdorff topology to $(S,D,(x_0,x_1,\ldots,x_k))$.  Assuming,
by using the Skorokhod representation theorem that this convergence
holds almost-surely, this means that for every $\eta>0$, and for every
$n$ large enough, it is possible to find a correspondence
$\mathcal{R}_n$ between $V(Q_n)$ and $S$, such that $(v_i,x_i)\in
\mathcal{R}_n$ for every $i\in \{0,1,\ldots,k\}$, and such that
\begin{equation}\label{eq:15}
\sup_{(v,x),(v',x')\in
  \mathcal{R}_n}\bigg|\Big(\frac{9}{8n}\Big)^{1/4}d_{Q_n}(v,v')-D(x,x')\bigg|\leq
\eta\, .
\end{equation}
%
Let us now assume that $\mathcal{A}_1(\eps,\beta)$ holds, and apply
the preceding observations for $k=2$.  Assume by contradiction that
with positive probability, along some (random) subsequence, there
exists a geodesic chain $\gamma_{(n)}$ in $Q_n$ from $v_1$ to $v_2$
such that no vertex of this chain lies at distance less than
$\eps(8n/9)^{1/4}$ from $v_0$.  Now choose $v^q_{(n)}$ on
$\gamma_{(n)}$, in such a way that
$(8n/9)^{-1/4}d_{Q_n}(v_0,v^q_{(n)})$ converges to some $q\in
(0,D(x_1,x_2))\cap \Q$. This entails in particular that
$(v_1,v^q_{(n)},v_2)$ are aligned (along the extraction considered
above). Then, let $x^q_{(n)}$ be such that $(v^q_{(n)},x^q_{(n)})\in
\mathcal{R}_n$. By diagonal extraction, we can assume that for every
$q$, $x^q_{(n)}$ converges to some $x^q\in S$, and using \eqref{eq:15}
entails both that $(x_1,x^q,x_2)$ are aligned, with $D(x_1,x^q)=q$,
and $D(x^q,x^{q'})=|q'-q|$. One concludes that the points $x^q$ are
dense in the image of a geodesic from $x_1$ to $x_2$, but this
geodesic is a.s.\ unique and has to be $\gamma$. By assumption,
$d_{Q_n}(v_0,v^q_{(n)})\geq \eps(8n/9)^{1/4}$ for every rational $q$
and $n$ chosen along the same extraction, so that \eqref{eq:15}
entails that $D(x_0,x^q)\geq \eps$ for every $q$, hence that $\gamma$
does not intersect $B_D(x_0,\eps/2)$, a contradiction.

Similarly, assume by contradiction that with positive probability, for
infinitely many values of $n$, there exists a geodesic chain
$\gamma_{(n)}$ from $v_0$ to $v_1$ such that every $v$ on this
geodesic with $d_{Q_n}(v,v_0)>\eps^{1-\beta}(8n/9)^{1/4}$, satisfies
also that $(v_1,v,v_2)$ are aligned (and similarly with the roles of
$v_1$ and $v_2$ interchanged). Similarly to the above, choose
$v^q_{(n)}$ on $\gamma_{(n)}$, in such a way that
$(8n/9)^{-1/4}d_{Q_n}(v_0,v^q_{(n)})$ converges to some $q\in
(\eps^{1-\beta},D(x_0,x_1))\cap \Q$. This entails in particular that
$(v_1,v^q_{(n)},v_2)$ are aligned (along some extraction). By using
the correspondence $\mathcal{R}_n$ and diagonal extraction, this
allows to construct a portion of geodesic in $(S,D)$ from $x_0$ to
$x_1$ lying outside of $B_D(x_0,\eps^{1-\beta})$, which visits only
points that are aligned with $x_1$ and $x_2$. The uniqueness of
geodesics allows to conclude that all points $x$ on $\gamma_{01}$
outside $B_D(x_0,\eps^{1-\beta})$ are in $\gamma$. The same holds for
$\gamma_{02}$ instead of $\gamma_{01}$ by the same argument, so
$\mathcal{A}_1(\eps,\beta)$ does not occur.

Next, we know that a.s.\ $x_0,x_1,x_2$ have no alignement relations,
and by \eqref{eq:15} this is also the case of $v_0,v_1,v_2$ for every
large $n$. A last use of \eqref{eq:15} shows that
$d_{Q_n}(v_0,v_1)\wedge d_{Q_n}(v_0,v_2)\geq
3\eps^{1-\beta}(8n/9)^{1/4}$ for $n$ large, from the fact that
$D(x_0,x_1)\wedge D(x_0,x_2)\geq 7\eps^{1-\beta}$ on the event
$\mathcal{A}_1(\eps,\beta)$.

Putting things together, we have obtained the claim on
$\mathcal{A}_1(\eps,\beta)$.  The statement concerning
$\mathcal{A}_2(\eps)$ is similar and left to the reader.  \cq

\section{Coding by labeled maps}\label{sec:coding-labeled-maps}

Our main tool for studying geodesic $\eps$-stars with $k$ arms is a
bijection \cite{miertess} between {\em multi-pointed delayed
  quadrangulations} and a class of labeled maps, which extends the
celebrated Cori-Vauquelin-Schaeffer bijection. The multi-pointed
bijection was used in \cite{miertess} to prove a uniqueness result for
typical geodesics that is related to the result of \cite{legall08}
that we already used in the present work. It was also used in
\cite{BoGu08b} to obtain the explicit form of the joint law of
distances between three randomly chosen vertices in the Brownian map
$(S,D)$. The way in which we use the multi-pointed bijection is in
fact very much inspired from the approach of \cite{BoGu08b}.

\subsection{The multi-pointed bijection}\label{sec:multi-point-scha}

\subsubsection{Basic properties}\label{sec:basic-properties}

Let $\bq\in \bQ$, and $\bv=(v_0,v_1,\ldots,v_k)$ be $k+1$ vertices of
$\bq$. Let also $\btau=(\tau_0,\tau_1,\ldots,\tau_k)$ be {\em delays} between
the points $v_i,0\leq i\leq k$, i.e.\ relative integers such that for
every $i,j\in \{0,1,\ldots,k\}$ with $i\neq j$,
\begin{eqnarray}
|\tau_i-\tau_j|<\dq(v_i,v_j)\, ,\label{eq:8}\\
 \dq(v_i,v_j)+\tau_i-\tau_j\in 2\N\, .\label{eq:9}
\end{eqnarray}
Such vertices and delays exist as soon as $\dq(v_i,v_j)\geq 2$ for
every $i\neq j$ in $\{0,1,\ldots, k\}$. We let $\bQ^{(k+1)}$ be the
set of triples $(\bq,\bv,\btau)$ as described, and we let
$\bQ^{(k+1)}_n$ be the subset of those triples such that $\bq$ has $n$
faces.  

On the other hand, a {\em labeled map} with $k+1$ faces is a pair
$(\bm,\bl)$ such that $\bm$ is a rooted map with $k+1$ faces, named
$f_0,f_1,\ldots,f_k$, while $\bl:V(\bm)\to \Z_+$ is a labeling function
such that $|\bl(u)-\bl(v)|\leq 1$ for every $u,v$ linked by an edge of
$\bm$.  If $\bm$ has $n$ edges, then $\bm$ has $n-k+1$ vertices by
Euler's formula.  We should also mention that the function $\bl$ and
the delays $\btau$ are defined up to a common additive constant, but
we are always going to consider particular representatives in the
sequel.

The bijection of \cite{miertess} associates with $(\bq,\bv,\btau)\in
\bQ^{(k+1)}_n$ a labeled map $(\bm,\bl)$ with $n$ edges, denoted by
$\Phi^{(k+1)}(\bq,\bv,\btau)$, in such a way that
$V(\bm)=V(\bq)\setminus\{v_0,v_1,\ldots,v_k\}$: This identification
will be implicit from now on. Moreover, the function $\bl$ satisfies
\begin{equation}\label{eq:10}
\bl(v)=\min_{0\leq i\leq k}(\dq(v,v_i)+\tau_i)\, ,\qquad v\in
V(\bm)\, ,
\end{equation}
where on the right-hand side one should understand that the vertex $v$
is a vertex of $V(\bq)$. The function $\bl$ is extended to $V(\bq)$ in
the obvious way, by letting $\bl(v_i)=\tau_i,0\leq i\leq k$. This
indeed extends \eqref{eq:10} by using \eqref{eq:8}, and we see in
passing that $\tau_i+1$ is the minimal label $\bl(v)$ for all vertices
$v$ incident to $f_i$. The interpretation of
the labels is the following. Imagine that $v_i$ is a source of liquid,
which starts to flow at time $\tau_i$. The liquid then spreads in the
quadrangulation, taking one unit of time to traverse an edge. The
different liquids are not miscible, so that they end up entering in
conflict and becoming jammed. The vertices $v$ such that
$$\bl(v)=\dq(v,v_i)+\tau_i<\min_{j\in \{0,1,\ldots,k\}\setminus\{i\}}(\dq(v,v_j)+\tau_j)$$
should be understood as the set of vertices that have only been attained by
the liquid starting from $v_i$. 

The case where there are ties is a little more elaborate as we have to
give {\em priority rules} to liquids at first encounter. The property
that we will need is the following. The label function $\bl$ is such
that $|\bl(u)-\bl(v)|=1$ for every adjacent $u,v\in V(\bq)$, so that
there is a natural orientation of the edges of $\bq$, making them
point toward the vertex of lesser label. We let
$\vec{E}^{\bv,\btau}(\bq)$ be the set of such oriented edges.  Maximal
oriented chains made of edges in $\vec{E}^{\bv,\btau}(\bq)$ are then
geodesic chains, that end at one of the vertices $v_0,v_1,\ldots,v_k$.
For every oriented edge $e\in \vec{E}^{\bv,\btau}(\bq)$, we consider
the oriented chain starting from $e$, and turning to the left as much
as possible at every step. For $i\in \{0,1,\ldots,k\}$, the set
$\vec{E}_i^{\bv,\btau}(\bq)$ then denotes the set of $e\in
\vec{E}^{\bv,\btau}(\bq)$ for which this leftmost oriented path ends
at $v_i$. One should see $\vec{E}_i^{\bv,\btau}(\bq)$ as the set of
edges that are traversed by the liquid emanating from source $v_i$.

\subsubsection{The reverse construction}\label{sec:reverse-construction}

We will not specify how $(\bm,\bl)$ is constructed from an element
$(\bq,\bv,\btau)\in \bQ^{(k+1)}$, but it is important for our purposes
to describe how one goes back from the labeled map $(\bm,\bl)$ to the
original map with $k+1$ vertices and delays. 
We first set a couple of
extra notions.
%

For each $i\in \{0,1,\ldots,k\}$, we can arrange the oriented edges of
$\vec{E}(f_i)$ cyclically in the so-called {\em facial order}: Since
$f_i$ is located to the left of the incident edges, we can view the
faces as polygons bounded by the incident edges, oriented so that they
turn around the face counterclockwise, and this order is the facial
order. If $e,e'$ are distinct oriented edges incident to the same
face, we let $[e,e']$ be the set of oriented edges appearing when
going from $e$ to $e'$ in facial order, and we let
$[e,e]=\{e\}$. Likewise, the oriented edges incident to a given vertex
$v$ are cyclically ordered in counterclockwise order when turning
around $v$. The {\em corner} incident to the oriented edge $e$ is a
small angular sector with apex $e_-$, that is delimited by $e$ and the
edge that follows around $v$: These sectors should be simply connected
and chosen small enough so that they are paiwise disjoint. We will
often assimilate $e$ with its incident corner. The label of a corner
is going to be the label of the incident vertex. In particular, we
will always adopt the notation $\bl(e)=\bl(e_-)$.

The converse construction from $(\bm,\bl)$ to $(\bq,\bv,\btau)$ goes
as follows.  Inside the face $f_i$ of $\bm$, let us first add an extra
vertex $v_i$, with label
\begin{equation}\label{eq:11}
\bl(v_i)=\tau_i=\min_{v\in V(f_i)}\bl(v)-1\, ,
\end{equation}
consistently with \eqref{eq:10} and the discussion below. We view
$v_i$ as being incident to a single corner $c_i$.  For every $e\in
\vec{E}(\bm)$, we let $f_i$ be the face of $\bm$ incident
to $e$, and define the {\em successor} of $e$ as the first corner $e'$
following $e$ in the facial order, such that $\bl(e')=\bl(e)-1$, we
let $e'=s(e)$. If there are no such $e'$, we let $s(e)=c_i$. The
corners $c_i$ themselves have no successors.

For every $e\in \vec{E}(\bm)$, we draw an {\em arc} between the corner
incident to $e$ and the corner incident to $s(e)$. It is possible to
do so in such a way that the arcs do not intersect, nor cross an edge
of $\bm$. The graph with vertex set $V(\bm)\cap\{v_0,v_1,\ldots,v_k\}$
and edge-set the set of arcs (so that the edges of $\bm$ are
excluded), is then a quadrangulation $\bq$, with distinguished
vertices $v_0,v_1,\ldots,v_k$ and delays $\tau_0,\tau_1,\ldots,\tau_k$
defined by \eqref{eq:11}. More precisely, if $e$ is incident to $f_i$,
then the arc from $e$ to $s(e)$ is an oriented edge in
$E_i^{\bv,\btau}(\bq)$, and every such oriented edge can be obtained
in this way: In particular, the chain $e,s(e),s(s(e)),\ldots$ from
$e_-$ to $v_i$ is the leftmost geodesic chain described when defining
the sets $E_i^{\bv,\btau}(\bq)$. 

To be complete, we should describe how the graph made of the arcs is
rooted, but we omit the exact construction as it is not going to play
an important role here. What is relevant is that for a given map
$(\bm,\bl)$, there are two possible rooting conventions for
$\bq$. Therefore, the mapping $\Phi^{(k+1)}$ associating a labeled map
$(\bm,\bl)$ with a delayed quadrangulation $(\bq,\bv,\btau)$ is
two-to-one. Consequently, the mapping $\Phi^{(k+1)}$ pushes forward
the counting measure on $\bQ^{(k+1)}_n$ to {\em twice} the counting
measure on labeled maps with $n$ edges, as far as we are interested in
events that do not depend on the root of $\bq$.

\subsection{Geodesic $r$-stars in
  quadrangulations}\label{sec:geodesic-r-stars}

We want to apply the previous considerations to the estimation of the
probabilities of the events
$\mathcal{A}_1^{(n)}(\eps,\beta),\mathcal{A}_2^{(n)}(\eps)$ of Section
\ref{sec:back-geodesic-stars}. To this end, we will have to specify
the appropriate discrete counterpart to the event
$\mathcal{G}(\eps,k)$ of \eqref{eq:7}. Contrary to the continuous
case, in quadrangulation there are typically many geodesic chains
between two vertices. In uniform quadrangulations with $n$ faces, the
geodesic chains between two typical vertices will however form a thin
pencil (of width $o(n^{1/4})$), that will degenerate to a single
geodesic path in the limit.

Let $k,r>0$ be integers.  We denote by $G(r,k)$ set of pairs
$(\bq,\bv)$ with $\bq\in \bQ$, $\bv=(v_0,v_1,\ldots,v_k)\in
V(\bq)^{k+1}$, and such that
\begin{itemize}
\item if $(v_0,v,v_i)$ are aligned for some $i\in \{1,2,\ldots,k\}$
  and $d_\bq(v,v_0)\geq r$, then $(v_0,v,v_j)$ are not aligned for
  every $j\neq i$,
\item no three distinct vertices in $\{v_0,v_1,\ldots,v_k\}$, taken in
  any order, are aligned in $\bq$, and $\min\{\dq(v_0,v_i),1\leq i\leq
  k\}\geq 3r$,
\end{itemize}
Let $(\bq,\bv)\in G(r,k)$, and fix $r'\in \{r+1,r+2,\ldots,2r\}$. We
let $\btau^{(r')}=(\tau^{(r')}_0,\tau^{(r')}_1,\ldots,\tau^{(r')}_k)$
be defined by
\begin{equation}\label{eq:12}
\left\{\begin{array}{l}\tau_0^{(r')}=-r'\\
\tau_i^{(r')}=-\dq(v_0,v_i)+r'\, ,\qquad 1\leq i\leq k\, .
\end{array}\right.
\end{equation}

\begin{lmm}\label{sec:geodesic-r-stars-1}
  If $(\bq,\bv)\in G(r,k)$ and $r'\in \{r+1,\ldots,2r\}$, then
  $(\bq,\bv,\btau^{(r')})\in\bQ^{(k+1)}$.
\end{lmm}

\proof Let us verify \eqref{eq:8}. We write $\btau$ instead of
$\btau^{(r')}$ for simplicity. We have
$\tau_0-\tau_i=\dq(v_0,v_i)-2r'$, and by the assumption that
$\dq(v_0,v_i)\geq 3r> r'$, we immediately get
$|\tau_0-\tau_i|<\dq(v_0,v_i)$. Next, if $i,j\in \{1,\ldots,k\}$ are
distinct, we have
$|\tau_i-\tau_j|=|\dq(v_0,v_j)-\dq(v_0,v_i)|<\dq(v_i,v_j)$, since
$(v_i,v_j,v_0)$ are not aligned, and neither are $(v_j,v_i,v_0)$.

We now check \eqref{eq:9}. First note that for every $i\in
\{1,\ldots,k\}$, $\dq(v_0,v_j)+\tau_i-\tau_0=2r'$ which is even. Next,
for $i,j\in \{1,\ldots,k\}$ distinct, consider the mapping $h:v\in
V(\bq)\mapsto \dq(v_i,v_j)-\dq(v,v_i)+\dq(v,v_j)$. We have $h(v_j)=0$,
which is even. Moreover, since $\bq$ is a bipartite graph, we have
$h(u)-h(v)\in \{-2,0,2\}$ if $u$ and $v$ are adjacent vertices. Since
$\bq$ is a connected graph, we conclude that $h$ takes all its values
in $2\Z$, so $h(v_0)$ is even, and this is \eqref{eq:9}.  \cq

\medskip

Under the hypotheses of Lemma \ref{sec:geodesic-r-stars-1}, let
$(\bm,\bl)=\Phi^{(k+1)}(\bq,\bv,\btau^{(r')})$, where $\Phi^{(k+1)}$
is the mapping described in Section \ref{sec:multi-point-scha}. The
general properties of this mapping entail that
\begin{equation}\label{eq:13}
  \min_{v\in V(f_0)}\bl(v)=-r'+1 > -2r
\end{equation}
and
$$\min_{v\in
  f_i}\bl(v)=-\dq(v_0,v_i)+r'-1<0\quad \mbox{ for every }\quad i\in
\{1,\ldots,k\}\, .$$ We now state a key combinatorial lemma.  Let
$\bLM^{(k+1)}$ set of labeled maps $(\bm,\bl)$ with $k+1$ faces such
that for every $i\in \{1,\ldots,k\}$, the face $f_0$ and $f_i$ have at
least one common incident vertex, and $\min_{v\in V(f_0\cap
  f_i)}\bl(v)=0$.

\begin{lmm}\label{sec:geodesic-r-stars-2}
  Let $(\bq,\bv)\in G(r,k)$, and $r'\in \{r+1,r+2,\ldots,2r\}$. Then
  the labeled map $(\bm,\bl)=\Phi^{(k+1)}(\bq,\bv,\btau^{(r')})$
  belongs to $\bLM^{(k+1)}$.
\end{lmm}

\proof Consider a geodesic chain
$\gamma_i=(e_1,\ldots,e_{\dq(v_0,v_i)})$ from $v_0$ to $v_i$, and let
$e=e_{r'+1},e'=\ov{e}_{r'}$ and $v=e_-=e'_-$ be the vertex visited by
this geodesic at distance $r'$ from $v_0$. Then it holds that
$\dq(v,v_i)=\dq(v_0,v_i)-r'$, so that 
$$\dq(v,v_0)+\tau_0=0=\dq(v,v_i)+\tau_i\, .$$
Let us show that $v\in V(f_0\cap f_i)$. Since $(\bq,\bv)\in G(r,k)$, we
know that $v$ is not on any geodesic from $v_0$ to $v_j$, for $j\in
\{1,\ldots,k\}\setminus\{i\}$. Therefore, 
$$\dq(v,v_j)+\tau_j=\dq(v,v_j)-\dq(v_0,v_j)+r'>-\dq(v_0,v)+r'=0\, .$$
From this, we conclude that $\bl(v)=0$, and that $v$ can be incident
only to $f_0$ or $f_i$. It is then obvious, since $\gamma_i$ is a
geodesic chain, that $e\in \vec{E}_i^{\bv,\btau}(\bq)$ and $e'\in
\vec{E}_0^{\bv,\btau}(\bq)$. From the reverse construction, we see
that $e$ and $e'$ are arcs drawn from two corners of the same vertex,
that are incident to $f_i$ and $f_0$ respectively. \cq

\medskip

In fact, the proof shows that all the geodesic chains from $v_0$ to
$v_i$ in $\bq$ visit one of the vertices of label $0$ in $V(f_0\cap
f_i)$. This will be useful in the sequel. 

We let $\bLM^{(k+1)}_n$ be the subset of elements with $n$ edges, and
$\LM^{(k+1)}_n$ be the counting measure on $\bLM^{(k+1)}_n$, so its
total mass is $\#\bLM^{(k+1)}_n$. We want to consider the asymptotic
behavior of this measure as $n\to\infty$, and for this, we need to
express the elements of $\bLM^{(k+1)}$ in a form that is appropriate
to take scaling limits.

%

\subsection{Decomposition of labeled
  maps in $\bLM^{(k+1)}$}\label{sec:decomp-label-maps}

It is a standard technique both in enumerative combinatorics and in
probability theory to decompose maps in simpler objects: Namely, a
homotopy type, or {\em scheme}, which is a map of fixed size, and a
labeled forest indexed by the edges of the scheme. See
\cite{okounkov00,ChMaSc,miertess,bettinelli10} for instance. Due to
the presence of a positivity constraint on the labels of vertices
incident to $f_0$ and $f_i$, this decomposition will be more elaborate
than in these references, it is linked in particular to the one
described in \cite{BoGu08b}.

\subsubsection{Schemes}\label{sec:schemes}

From this point on, the notation $k$ will always stand for an integer
$k\geq 2$.  We call {\em pre-scheme} with $k+1$ faces, an unrooted map
$\sg_0$ with $k+1$ faces named $f_0,f_1,\ldots,f_k$, in which every
vertex has degree at least $3$, and such that for every $i\in
\{1,2,\ldots,k\}$, the set $V(f_0\cap f_i)$ of vertices incident to
both $f_0$ and $f_i$ is not empty.

It is easy to see, by applying the
Euler formula, that there are only a finite number of pre-schemes with
$k+1$ faces: Indeed, it has at most $3k-3$ edges and $2k-2$ vertices,
with equality if and only if all vertices have degree exactly $3$, in
which case we say that $\sg_0$ is {\em dominant}, following
\cite{ChMaSc}. 

\begin{figure}
\begin{center}
\includegraphics{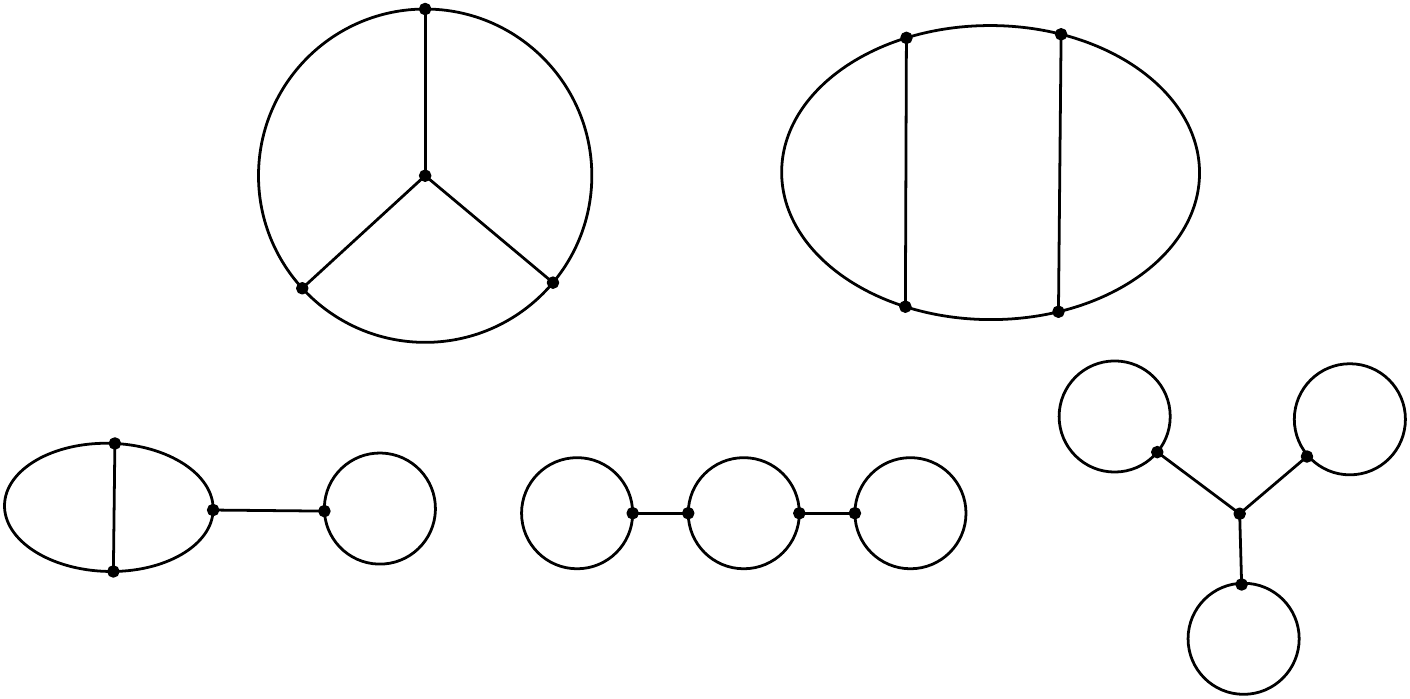}
\end{center}
\caption{The five dominant pre-schemes with $4$ faces, where $f_0$ is
 the unbounded face, and after forgetting the names $f_1,f_2,f_3$ of
 the other faces (there are sixteen dominant pre-schemes with $4$
 faces). In the first two cases, the boundary of the exterior
 face is a Jordan curve, while in the last three it is not simple.}
\label{fig:shapes}
\end{figure}

A {\em scheme} with $k+1$ faces is an unrooted map that can be
obtained from a pre-scheme $\sg_0$ in the following way. For every
edge of $\sg_0$ that is incident to $f_0$ and some face $f_i$ with
$i\in \{1,2,\ldots,k\}$, we allow the possibility to split it into two
edges, incident to a common, distinguished new vertex of degree $2$
called a {\em null-vertex}. Likewise, some of the vertices of
$V(f_0\cap f_i)$ with $i\in \{1,2,\ldots,k\}$ are allowed to be
distinguished as null-vertices. These operations should be performed
in such a way that every face $f_i$, for $i\in \{1,2,\ldots,k\}$, is
incident to at least one null-vertex.  Furthermore, a null vertex of
degree $2$ is not allowed to be adjacent to any other null vertex (of
any degree).

In summary, a scheme is a map $\sg$ with $k+1$ faces
labeled $f_0,f_1,\ldots,f_k$, such that
\begin{itemize}
\item For every $i\in \{1,2,\ldots,k\}$, the set $V(f_0\cap f_i)$ is
  not empty,
\item the vertices of $\sg$ have degrees greater than or equal to $2$,
\item vertices of degree $2$ are all in $\bigcup_{i=1}^k V(f_0\cap
  f_i)$, and no two vertices of degree $2$ are adjacent to each other,
\item every vertex of degree $2$, plus a subset of the other vertices
  in $\bigcup_{i=1}^k V(f_0\cap f_i)$, are distinguished as
  null-vertices, in such a way that every face $f_i$ with $i\in
  \{1,2,\ldots,k\}$ is incident to a null-vertex, and no degree-$2$
  vertex is adjacent to another null vertex.
\end{itemize}

Since the number of pre-schemes with $k+1$ faces is finite, the number
of schemes is also finite. Indeed, passing from a pre-scheme to a
scheme boils down to specifying a certain subset of edges and vertices
of this pre-scheme. We say that the scheme is {\em dominant} if it can
be obtained from a dominant pre-scheme, and if is has exactly $k$
null vertices, which are all of degree $2$. Since a dominant
pre-scheme has $3k-3$ edges and $2k-2$ vertices, we see that a
dominant scheme has $4k-3$ edges and $3k-2$ vertices.  We let
$\bS^{(k+1)}$ be the set of schemes with $k+1$ faces, and
$\bS_d^{(k+1)}$ the subset of dominant ones.

The vertices of a
scheme can be partitioned into three subsets: 
\begin{itemize}
\item 
The set $V_N(\sg)$ of
null vertices, 
\item
The set $V_I(\sg)$ of vertices that are incident to
$f_0$ and to some other face among $f_1,\ldots,f_k$, but which are not
in $V_N(\sg)$, and 
\item
The set $V_O(\sg)$ of all other
vertices. 
\end{itemize}
Similarly, the edges of $\sg$ can be partitioned into 
\begin{itemize}
\item 
the
set $E_N(\sg)$ of edges incident to $f_0$ and to some other face among
$f_1,\ldots,f_k$, and having at least one extremity in $V_N(\sg)$, 
\item
the
set $E_I(\sg)$ of edges that are incident to $f_0$ and some other face
in $f_1,\ldots,f_k$, but that are not in $E_N(\sg)$, and
\item the set
$E_O(\sg)$ of all other edges.
\end{itemize}

\begin{figure}
\begin{center}
\includegraphics{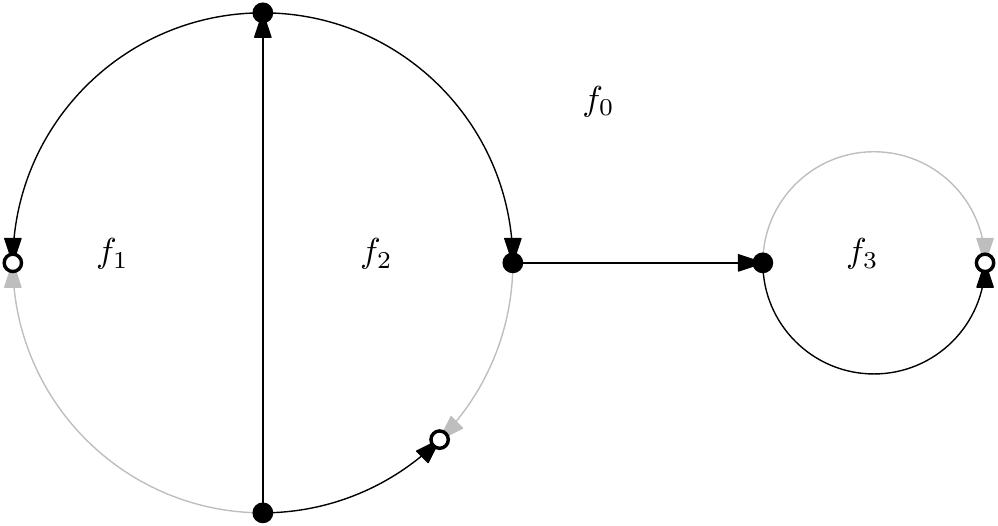}
\end{center}
\caption{A scheme in $\bS^{(3)}_d$, where we indicated thin edges in light
  gray, and specified the orientation conventions. Here the cardinalities of
  $V_N(\sg),V_I(\sg),V_O(\sg)$ are respectively $3,4,0$, and the
  cardinalities of $E_N(\sg),E_I(\sg),E_O(\sg)$ are respectively
  $6,1,2$. }
\label{fig:scheme1}
\end{figure}

It will be convenient to adopt once and for all an orientation
convention valid for every scheme, meaning that every element of
$E(\sg)$ comes with a privileged orientation. We add the constraint that
an edge in $E_N(\sg)$ is always oriented towards a vertex of
$V_N(\sg)$, with priority given to those with degree $2$: That is, if
an edge of $E_N(\sg)$ links two vertices of $V_N(\sg)$, the
orientation is arbitrary if both vertices have degree at least $3$,
and points towards to unique incident vertex of degree $2$
otherwise. The other orientations are arbitrary, as in Figure
\ref{fig:scheme1}. We let $\check{E}(\sg)$ be the orientation
convention of $\sg$.

For every null vertex $v$ of degree $2$, there is only one $e\in
\check{E}(\sg)$ satisfying both $e_+=v$, and $e\in \vec{E}(f_0)$ (by
definition, the face incident to $\ov{e}$ is then some other face
among $f_1,\ldots,f_k$).  The corresponding (non-oriented) edge is
distinguished as a {\em thin edge}.  Similarly, any edge of $E_N(\sg)$
incident to at least one null vertex of degree at least $3$ is counted
as a thin edge. We let $E_T(\sg)$ be the set of thin edges. Dominant
schemes are the ones having exactly $k$ null-vertices, which are all
of degree $2$: These also have $k$ thin edges.

\subsubsection{Labelings and
  edge-lengths}\label{sec:label-edge-lengths}

An {\em admissible labeling} for a (planted) scheme $\sg$ is a family $(\ell_v,v\in
V(\sg))\in \Z^{V(\sg)}$ such that
\begin{enumerate}
\item $\ell_v=0$ for every $v\in V_N(\sg)$,
\item $\ell_v> 0$ for every $v\in V_I(\sg)$,
\end{enumerate}
A {\em family of edge-lengths} for a scheme $\sg$ is a family
$(r_e,e\in E(\sg))\in \N^{E(\sg)}$ of positive integers indexed by the
edges of $\sg$. A family of edge-lengths can be
naturally seen as being indexed by oriented edges rather
than edges, by setting $r_e=r_{\ov{e}}$ to be equal to the edge-length
of the edge with orientations $e,\ov{e}$.

\subsubsection{Walk networks}\label{sec:walk-networks-snakes}

We call {\em Motzkin walk}\footnote{This is not a really standard
  denomination in combinatorics, where Motzkin paths usually denote
  paths that are non-negative besides the properties we
  require} a finite sequence $(M(0),M(1),\ldots,M(r))$ with values in
$\Z$, where the integer $r\geq 0$ is the duration of the walk, and
$$M(i)-M(i-1)\in \{-1,0,1\}\, ,\qquad 1\leq
  i\leq r\, .$$

Given a scheme with admissible labeling $(\ell_v,v\in V(\sg))$ and
edge-lengths $(r_e,e\in \vec{E}(\sg))$, a {\em compatible walk
  network} is a family $(M_e,e\in \vec{E}(\sg))$ of Motzkin walks
indexed by $\vec{E}(\sg)$, where for every $e\in
\vec{E}(\sg)$,
\begin{enumerate}
\item
 one has $M_e(i)=M_{\ov{e}}(r_e-i)$ for every $0\leq i\leq r_e$, 
\item
one has $M_e(0)=\ell_{e_-}$, 
\item if $e\in E_N(\sg)\cup E_I(\sg)$, then $M_e$ takes only
  non-negative values. If moreover $e$ is not a thin edge, then all
  the values taken by $M_e$ are positive, except $M_e(r_e)$ (when $e$
  is canonically oriented towards the only null vertex is is incident
  to).
\end{enumerate}
The first condition says that the walks can really be defined as being
labeled by edges of $\sg$ rather than oriented edges: The family
$(M_e,e\in \vec{E}(\sg))$ is indeed entirely determined by $(M_e,e\in
\check{E}(\sg))$, where $\check{E}(\sg)$ is the orientation convention
on $E(\sg)$. Also, note that 1.\ and 2.\ together imply that
$M_e(r_e)=\ell_{e_+}$, so we see that $M_e(r_e)=0$ whenever $e\in
E_N(\sg)$ (with orientation pointing to a vertex of $V_N(\sg)$, which
is the canonical orientation choice we made).  The distinction arising
in 3.\ between thin edges and non-thin edges in $E_N(\sg)$ is slightly
annoying, but unavoidable as far as exact counting is involved. Such
distinctions will disappear in the scaling limits studied in Section
\ref{sec:scal-limits-label}.

\subsubsection{Forests and discrete snakes}

Our last ingredient is the notion of {\em plane forest}. We will not
be too formal here, and refer the reader to
\cite{miertess,bettinelli10} for more details. A plane tree is a
rooted plane map with one face, possibly reduced to a single vertex,
and a plane forest is a finite sequence of plane trees
$(\bt_1,\bt_2,\ldots,\bt_r)$. We view a forest itself as a plane map,
by adding an oriented edge from the root vertex of $\bt_i$ to the root
vertex of $\bt_{i+1}$ for $1\leq i\leq r-1$, and adding another such
edge from the root vertex of $\bt_r$ to an extra vertex. These special
oriented edges are called {\em floor edges}, and their incident
vertices are the $r+1$ floor vertices. The vertex map, made of a
single vertex and no edge, is considered as a forest with no tree.

With a scheme $\sg$, an admissible labeling $(\ell_v,v\in V(\sg))$ and
a compatible walk network $(M_e,e\in \vec{E}(\sg))$, a {\em compatible
  labeled forest} is the data, for every $e\in \vec{E}(\sg)$, of a
plane forest $F_e$ with $r_e$ trees, and with an integer-valued
labeling function $(L_e(u),u\in V(F_e))$ such that $L_e(u)=M_e(i)$ if
$u$ is the $i+1$-th floor vertex of $F_e$, for $0\leq i\leq r_e$ and
$L_e(u)-L_e(v)\in \{-1,0,1\}$ whenever $u$ and $v$ are adjacent
vertices in the same tree of $F_e$, or adjacent floor vertices.

In order to shorten the notation, we can encode labeled forests in
discrete processes called {\em discrete snakes}. If $\bt$ is a rooted
plane tree with $n$ edges, we can consider the facial ordering
$(e^{(0)},e^{(1)},\ldots,e^{(2n-1)})$ of oriented edges starting from
its root edge, and let 
$$C_\bt(i)=d_{\bt}(e^{(0)}_-,e^{(i)}_-)\, ,\qquad 0\leq i\leq 2n-1\,
,$$ and then $C_{\bt}(2n)=0$ and $C_{\bt}(2n+1)=-1$. The sequence
$(C_{\bt}(i),0\leq i\leq 2n+1)$ is called the contour sequence of
$\bt$, and we turn it into a continuous function defined on the time
interval $[0,2n+1]$ by linearly interpolating between values taken at
the integers. Roughly speaking, for $0\leq s\leq 2n$, $C_\bt(s)$ is
the distance from the root of the tree at time $s$ of a particle going
around the tree at unit speed, starting from the root.

If $F$ is a plane forest with trees $\bt_1,\ldots,\bt_r$, the contour
sequence $C_F$ is just the concatenation of
$r+C_{\bt_1},r-1+C_{\bt_2},\ldots,1+C_{\bt_r}$, starting at $r$ and
finishing at $0$ at time $r+2n$, where $n$ is the total number of
edges in the forest distinct from floor edges. Note the fact that the
sequence visits $r-i$ for the first time when it starts exploring the
$i+1$-th tree. For simplicity, we still denote the sequence by $F$. If
$L$ is a labeling function on $F$, the label process is defined by
letting $L_F(i)$ be the label of the corner explored at the $i$-th step
of the exploration. 
Both processes $C_F,L_F$ are extended by linear interpolation between
integer times.

\begin{figure}[htb!]
\begin{center}
\includegraphics[scale=.8]{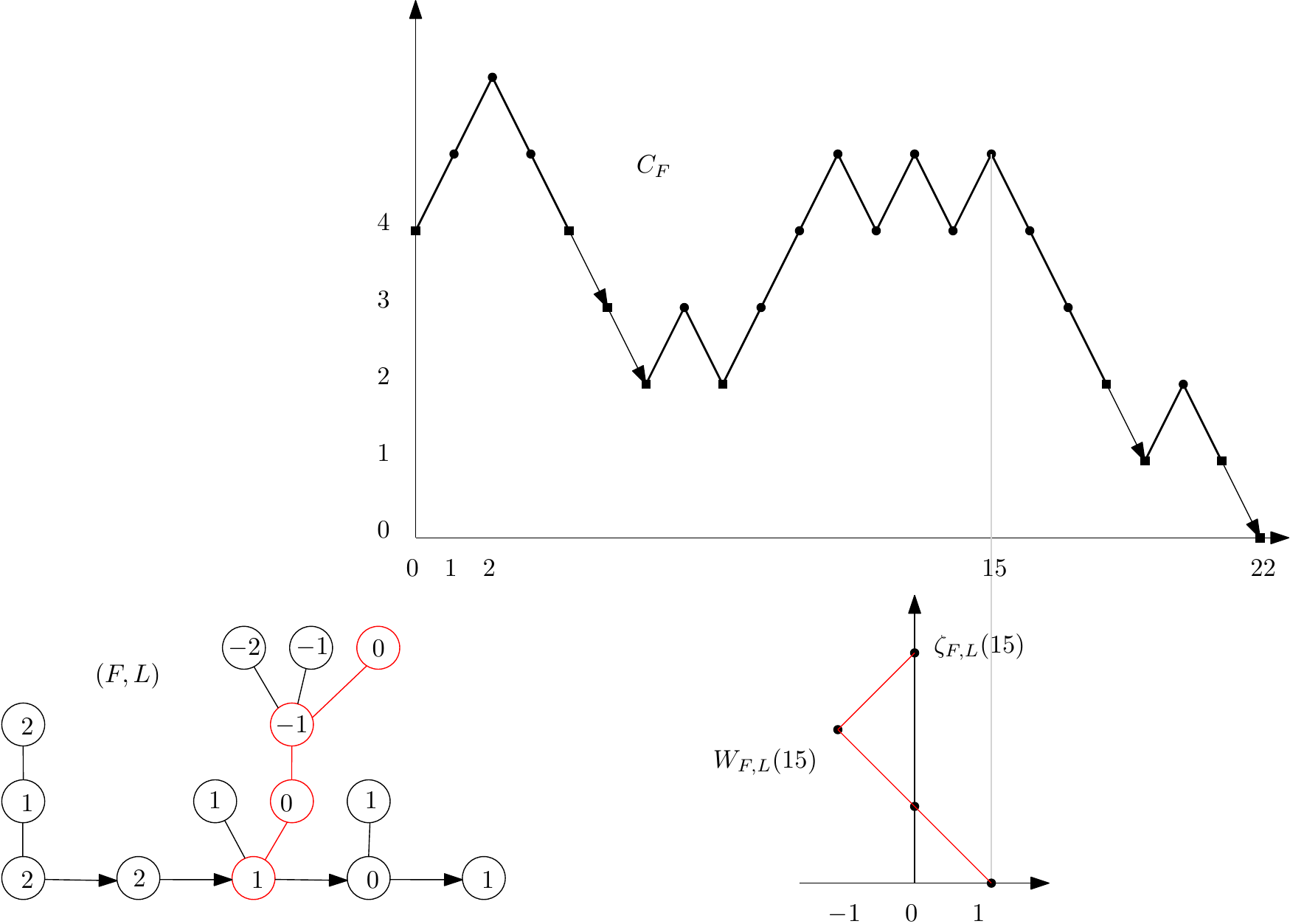}
\end{center}
\caption{A labeled forest with $4$ trees and $22$ oriented edges, 
its contour process, and the associated
  discrete snake evaluated at time $15$}
\end{figure}

The information carried by $(F,L)$ can be summarized into one unique
path-valued process
$$(W_{F,L}(i), 0\leq i\leq 2n+r)\, ,\quad \mbox{ with }\qquad W_{F,L}(i)=(W_{F,L}(i,j),0\leq
j\leq \zeta_F(i))\, ,$$ where $\zeta_F(i)$ is the distance to the
floor of the vertex $u^{(i)}=e^{(i)}_-$ of $F$ visited at the $i$-th
step in contour order, and $W_{F,L}(i,j)$ is the label $L(u)$ of the
ancestor $u$ of $u^{(i)}$ at distance $j$ from the floor. So, for
every $i\in \{0,1,\ldots,2n+r\}$, $W_{F,L}(i)$ is a finite sequence
with length $\zeta_F(i)$, and this sequence is a Motzkin walk. The
initial value $W_{F,L}(0)$ is the Motzkin walk
$(M(0),M(1),\ldots,M(r))$ given by the labels of the floor vertices of
$F$. Finally, we extend $W_{F,L}$ to a process in $\CC(\CC(\R))$ by
interpolation: $W_{F,L}(i)$ is now really a path in $\CC(\R)$ obtained
by interpolating linearly between integer times $(W_{F,L}(i,j),0\leq
j\leq \zeta_F(i))$, and for $s\in [i,i+1]$ we simply let $W_{F,L}(s)$
be the path $(W_{F,L}(i+1,t),0\leq t\leq \zeta_F(i)+s-i)$ if
$\zeta_F(i)<\zeta_F(i+1)$, and $(W_{F,L}(i,t),0\leq t\leq
\zeta_F(i)-s+i)$ if $\zeta_F(i)>\zeta_F(i+1)$.

We define a {\it discrete snake} to be a process $W_{F,L}$ obtained in
this way from some labeled forest $(F,L)$. From $W_{F,L}$, it is
possible to recover $(F,L)$. In particular, if
$\un{\zeta}_F(i)=\min_{0\leq i'\leq i}\zeta_F(i')$, then for $i\in
\{0,1,\ldots,2n+r-1\}$, $r+1-\un{\zeta}_F(i)$ is the label of the tree
of $F$ visited at time $i$ in the contour process. In particular, we
have
\begin{align*}
C_F(i)=r+1+\zeta_F(i)-\un{\zeta}_F(i)\, ,&\qquad &0\leq i\leq 2n+r\, ,\\
L_F(i)=W_{F,L}(i,\zeta_F(i))\, ,&\qquad& 0\leq i\leq 2n+r\, ,\\
M(j)= W_{F,L}(\inf\{i\geq 0:\zeta_F(i)=r-j\},r-j)\, ,&\qquad &0\leq
j\leq r\, .
\end{align*}

From this, the data of a scheme $\sg$, an admissible labeling
$(\ell_v,v\in V(\sg))$, admissible edge-lengths $(r_e,e\in E(\sg))$, a
compatible walk network $(M_e,e\in E(\sg))$, and compatible labeled
forests $((F_e,L_e),e\in \vec{E}(\sg))$ can be summed up in the family
$(\sg,(W_e,e\in \vec{E}(\sg)))$ where $W_e$ is the discrete snake
associated with $(F_e,L_e)$. We call a family $(W_e,e\in
\vec{E}(\sg))$ obtained in this way an {\em admissible family of
  discrete snakes} on the scheme $\sg$.

\subsubsection{The reconstruction}

Let us reconstruct an element of $\bLM^{(k+1)}$, starting from a
scheme $\sg$, and an admissible family of discrete snakes $(W_e,e\in
\vec{E}(\sg))$. The latter defines labeling, edge-lengths, a walk
network and a family of labeled forests, and we keep the same notation
as before.

First, we label every vertex of $\sg$ according to $(\ell_v,v\in
V(\sg))$. Second, we replace every edge $e$ of $\sg$ with a chain of
$r_e$ edges. Since $r_e=r_{\ov{e}}$, this is really an operation on
edges rather than oriented edges. 
Then, the vertices inside each of these chains are
labeled according to $M_e(0),M_e(1),\ldots,M_e(r_e)$.  Since
$M_e(0)=\ell_{e_-}$ and $M_e(r_e)=\ell_{e+}$, these labelings are
consistent with the labeling at the vertices of $\sg$.  At this stage,
we get a map with $k+1$ faces and labeled vertices of degree at least
$2$.

Then, we graft the labeled forest $(F_e,L_e)$ in such a way that the
$r_e$ floor edges of $F_e$ coincide with the oriented edges of the
chain with $r_e$ edges corresponding to $e$, and are all incident to
the same face as $e$ (i.e.\ the face located to the left of $e$). Note
that at this step the construction depends strongly on the orientation
of $e$, and that $F_e$ and $F_{\ov{e}}$ can be very different, despite
having the same number of trees, and the fact that the labels of the
floor vertices are given by $M_e$ and $M_{\ov{e}}$ respectively.

This yields a labeled map $(\bm,\bl)\in \bLM^{(k+1)}$. To be completely
accurate, we need to specify the root of $\bm$. To this end, 
we mark one of the oriented edges of
the forest $F_{e^*}$ (comprising the $r$ {\em oriented} floor
edges) for some $e_*\in \vec{E}(\sg)$. 
This edge specifies the root edge of $\bm$. Note that the
number of oriented edges of $\bm$ equals
\begin{equation}\label{eq:18}
\#\vec{E}(\bm)=\sum_{e\in \vec{E}(\sg)}\#\vec{E}(F_e)\, .
\end{equation}

This construction can be inverted: Starting from a labeled map
$(\bm,\bl)\in \bLM^{(k+1)}$, we can erase the degree-$1$ vertices
inductively, hence removing families of labeled forests grafted on
maximal chains, joining two vertices with degrees $\geq 3$ by passing
only through degree-$2$ vertices. Any given such chain is incident to
one or two faces only. The resulting map has only vertices of degrees
at least $2$, we call it $\bm'$. In turn, the degree-$2$ vertices of
$\bm'$ can be deleted, ending with a pre-scheme $\sg_0$, each edge of
which corresponds to a maximal chain of vertices with degree $2$ in
$\bm'$. Its faces are the same as $\bm,\bm'$ and inherit their names
$f_0,\ldots,f_k$. If a maximal chain in $\bm'$ is incident to $f_0$
and $f_i$ for some $i\in \{1,2,\ldots,k\}$, and if the label of one of
the degree-$2$ vertices of this chain is $0$ while the labels of both
extremities are strictly positive, then we add an extra degree-$2$
vertex to the corresponding edge of $\sg_0$, that is distinguished as
a null-vertex. Likewise, a vertex of label $0$ with degree at least
$3$ in $\bm'$ that is incident to $f_0$ and some other face is
distinguished as a null-vertex. This family of extra null vertices
turns $\sg_0$ into a scheme $\sg$, since by definition of
$\bLM^{(k+1)}$, every face $f_i$ with index $i\in \{1,2,\ldots,k\}$
has at least an incident vertex with label $0$ that is also incident
to $f_0$.

It remains to construct walks indexed by the edges of the scheme.
Consider a maximal chain as in the previous paragraph, corresponding
to an edge $e\in E(\sg)$. If $e$ is not incident to a null-vertex of
degree $2$, then the labels of the successive vertices of the chain
define a walk $M_e$ with positive duration, equal to the number of
edges in the chain. By default, the chain can be oriented according to
the orientation convention of $e$ in $\check{E}(\sg)$, which specifies
the order in which we should take the labels to define $M_e$ (changing
the orientation would define $M_{\ov{e}}$, which amounts to property
1.\ in the definition of compatible walk networks).  On the other
hand, if $e$ is incident to a null vertex of $\sg$ with degree $2$,
then $e$ is obtained by splitting an edge $e''$ of $\sg_0$ in two
sub-edges, $e,e'$. For definiteness, we will assume that $e$ is the
thin edge, meaning that when $e$ is oriented so that it points towards
the null vertex with degree $2$, then $f_0$ lies to its left. The
canonical orientation for $e'$ makes it point to the null vertex of
degree $2$ as well, as is now customary.  The chain of $\bm'$ that
corresponds to the edge $e''$, when given the same orientation as $e$,
has vertices labeled $l_0,l_1,\ldots,l_r$, in such a way that $\{i\in
\{1,\ldots,r-1\}:k_i=0\}$ is not empty.  Let $T=\max\{i\in
\{1,\ldots,r-1\}:l_i=0\}$, and set
$$M_e=(l_0,\ldots,l_T)\, ,\qquad
M_{e'}=(l_r,l_{r-1},\ldots,l_T)\, .$$ This defines two walks with
positive durations ending at $0$, the second one with only positive
values except at the ending point, and the first one takes
non-negative values (the starting value being positive).  We end up
with a scheme, admissible labelings and edge-lengths, and compatible
families of walks and forests.

To sum up our study, we have the following result. 

\begin{prp}\label{sec:decomp-label-maps-1}
The data of
\begin{itemize}
\item a scheme $\sg$,
\item an admissible family of
  discrete snakes $(W_e,e\in \vec{E}(\sg))$, 
\item an extra
  distinguished oriented edge in $F_{e_*}$, for some $e_*\in \vec{E}(\sg)$
\end{itemize}
determines a unique element in $\bLM^{(k+1)}$, and every such element
can be uniquely determined in this way.
\end{prp}

The only part that requires further justification is the word
``uniquely'' in the last statement: We have to verify that two
different elements $(\sg,(W_e,e\in \vec{E}(\sg))),(\sg',(W_e',e\in
\vec{E}(\sg)))$ cannot give rise to the same element in
$\bLM^{(k+1)}$. This is due to the fact that schemes have a trivial
automorphism group, i.e.\ every map automorphism of $\sg$ that
preserves the labeled faces is the identity automorphism. To see this,
note that there are at least two distinct indices $i,j\in
\{0,1,\ldots,k\}$ such that $f_i$ and $f_j$ are incident to a common
edge. If there is a unique oriented edge $e$ incident to $f_i$ with
$\ov{e}$ incident to $f_j$, then any map automorphism preserving the
labeled faces should preserve this edge, hence all the edges by
standard properties of maps. On the other hand, if there are several
edges incident to both $f_i$ and $f_j$, then there are multiple edges
between the vertices corresponding to $f_i$ and $f_j$ in the dual
graph of $\sg$, and by the Jordan Curve Theorem these edges split the
sphere into a collection $G_1,G_2,\ldots,G_p$ of $2$-gons. Since
$k\geq 2$, a scheme has at least three faces, so there exists $r\in
\{0,1,2,\ldots,k\}\setminus\{i,j\}$, and we can assume that the vertex
corresponding to $f_r$ in the dual graph of $\sg$ lies in $G_1$, up to
renumbering. But then, a graph automorphism preserving the labeled
faces, which boils down to a graph automorphism of the dual preserving
the labeled vertices, should preserve $G_1$, and in particular, it
should fix its two boundary edges, oriented from $f_i$ to
$f_j$. Therefore, this automorphism has to be the identity.

Note that this argument is only valid in planar geometry, and does not
apply to surfaces of positive genus.

\subsection{Using planted schemes to keep track of the
  root}\label{sec:keeping-track-root}

In this section we present a variant of the preceding description,
that allows to keep track of the root of the labeled map by an
operation on the scheme called {\em planting}. A {\em planted scheme}
satisfies the same definition as scheme, except that we allow exactly
one exceptional vertex $v_{**}$ with degree $1$. In the canonical
orientation, the edge incident to the only degree-one vertex always
points towards this vertex. We define the sets
$V_N(\sg),V_I(\sg),V_O(\sg)$ and $E_N(\sg),E_I(\sg),E_O(\sg)$ in the
same way as we did for schemes. In particular, one will note that
$v_{**}$ is always in $V_O(\sg)$, since it can only be incident to a
single face, while the only edge $e_{**}$ incident to $v_{**}$ is in $E_O(\sg)$
for the same reason. We let $\dot{\bS}^{(k+1)}$ be the set of planted
schemes with $k+1$ faces, and $\dot{\bS}_d^{(k+1)}$ the set of
dominant planted schemes, i.e.\ those having exactly $k$
null-vertices, all with degree $2$, all other vertices except $v_{**}$
being of degree $3$. Due to the distinguished nature of the edge
incident to $v_{**}$, planted schemes always have a trivial automorphism
group (not only in planar geometry).

The definition of admissible labelings, edge-lengths, walk networks
and discrete snakes associated with a planted scheme is the same as in
the non-planted case, with only one exception, stating that the path
$M_{e_{**}}$ is allowed to have duration $r_{e_{**}}$ equal to $0$, as
opposed to all others having positive durations, and the number of
trees in the labeled forests encoded by $W_{e_{**}},W_{\ov{e}_{**}}$
are both equal to $r_{e_{**}}+1$. Moreover, if $e'$ is the edge of
$\vec{E(\sg)}$ that comes after $e_{**}$ in clockwise order, then we
lower the number of trees of the forest $F_{e'}$ to $r_{e'}-1$ instead
of $r_{e'}$.

Let us explain how to construct a pair $(\sg,(W_e)_{e\in
  \vec{E}(\sg)})$, with $\sg$ a planted scheme and $(W_e)_{e\in
  E(\sg)}$ an admissible family of discrete snakes, starting from a
rooted labeled map $(\bm,\bl)$. The construction is very similar as in
the preceding section, so we will omit some details, and focus on how
we obtain the extra edge $e_{**}$, and the associated snakes
$W_{e_{**}},W_{\ov{e}_{**}}$.

Let $e_*$ be the root edge of $\bm$ and $v_*=(e_*)_-$ be its
origin. Recall that $\bm'$ is the map obtained from $\bm$ after
inductively removing inductively all the vertices of degree $1$, and
that $\bm$ can be obtained by grafting a tree component at each corner
of $\bm'$. The root $e_*$ belongs to one of these trees, let us call
it $\bt_*$, and we consider the chain
$v^{(0)}=v_*,v^{(1)},\ldots,v^{(r_{e_{**}})}$ down from $v_*$ to the
root of this tree (it can occur that $r_{e_{**}}$, meaning that $v_*$
belongs to $V(\bm')$). This chain is the one that gives rise to the
extra edge $e_{**}$ (i.e.\ even if it has zero length), and the labels
are $\ell_{(e_{**})_-}=\bl(v^{(r_{e_{**}})})$ and
$\ell_{(e_{**})_+}=\bl(v^{(0)})$. Likewise, the labels along the
chain define the path $M_{e_{**}}$.

\begin{figure}[htb!]
\begin{center}
\includegraphics{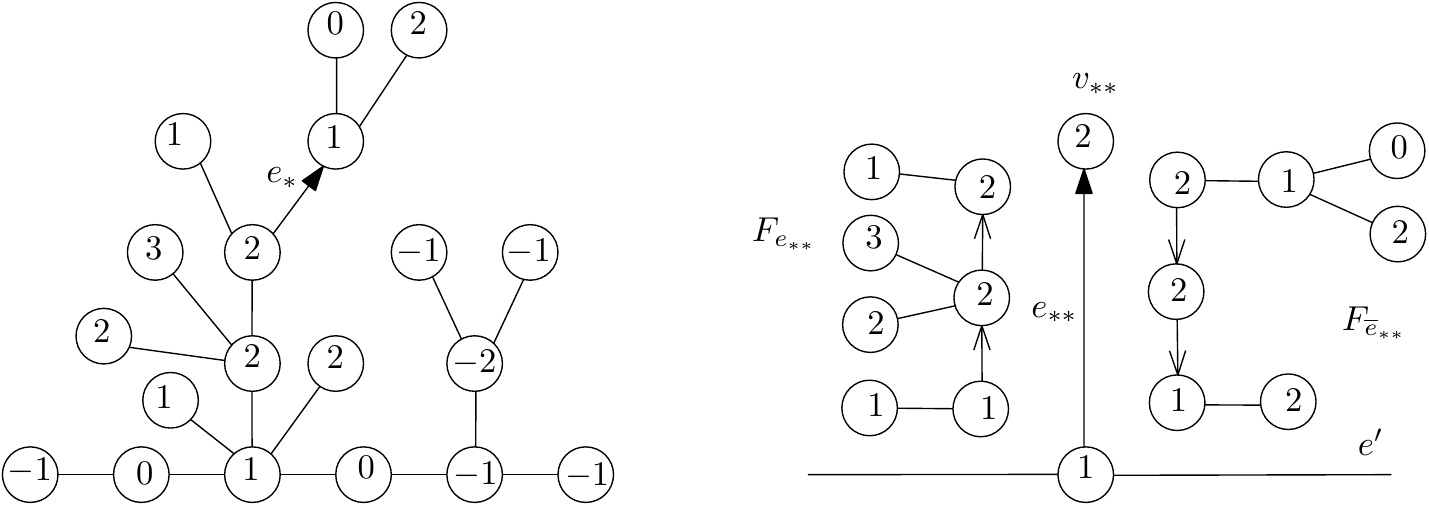}
\end{center}
\caption{The planting convention: On the left, a small
  portion of $\bm$ around its root $e_*$ is represented, in which the bottom vertices and
  edges all belong to the map $\bm'$. The root $e_*$ in $\bm$
  determines the location of the special edge $e_{**}$, and determines
  two forests with a positive number of trees, but no terminating
  floor edge. Observe that $e_*$ could well be an oriented edge in
  $\vec{E}(\bm')$, which occurs precisely when $F_{\ov{e}_{**}}$ is
  made of a single tree ($r_{e_{**}}=0$) made of a single vertex. The
  forest associated with $e'$ has one tree less than the forest
  associated with $\ov{e}'$, this tree being ``stolen'' as the last
  tree of the forest $F_{\ov{e}_{**}}$. }
\label{fig:plantedscheme}
\end{figure}

The tree $\bt_*$ and the distinguished chain between the root corner
of $\bt_*$ and $e_*$ can be seen in turn as being made of two forests
$F_{e_{**}},F_{\ov{e}_{**}}$ with the same (positive) number of trees,
equal to $r_{e_{**}}$, and in which we decide to forget the last floor
edge, which does not play a role here. The last tree explored in
$F_{\ov{e}_{**}}$ is then ``stolen'' from the forest $F_{e'}$, as
explained in Figure \ref{fig:plantedscheme}. We leave the last details
of the construction as an exercise to the reader.

\section{Scaling limits of labeled maps}\label{sec:scal-limits-label}

The description of labeled maps from Proposition
\ref{sec:decomp-label-maps-1} is particularly appropriate when one is
interested in taking scaling limits. We first introduce the proper
notion of ``continuum labeled map''. 

\subsection{Continuum measures on labeled
  maps}\label{sec:cont-meas-label}

Let $(E,d)$ be a metric space. We let $\CC(E)$ be the set of
$E$-valued continuous paths, i.e.\ of continuous functions
$f:[0,\zeta]\to \R$, for some $\zeta=\zeta(f)\geq 0$ called the {\em
  duration} of $f$. This space is endowed with the distance
$$\mathrm{dist}(f,g)=\sup_{t\geq 0}d\big(f(t\wedge \zeta(f)),g(t\wedge
\zeta(g))\big)+|\zeta(f)-\zeta(g)|\, ,$$ that makes it a Polish space
if $E$ is itself Polish. Sometimes, we will also have to consider
continuous functions $f:\R_+\to \R$ with infinite duration, so we let
$\ov{\CC}(E)$ be the set of continuous functions with finite or
infinite duration.

We let $\bCLM^{(k+1)}$ be the set of pairs of the form
$(\sg,(W_e)_{e\in \vec{E}(\sg)})$, where $\sg$ is a scheme with $k+1$
faces, and for every $e\in \vec{E}(\sg)$, $W_e$ is an element of
$\CC(\CC(\R))$, meaning that for every $s\in [0,\zeta(W_e)]$, $W_e(s)$
is a function in $\CC(\R)$ that can be written $(W_e(s,t),0\leq t\leq
\zeta(W_e(s)))$. The space $\bCLM^{(k+1)}$ is a Polish space, a
complete metric being for instance the one letting $(\sg,(W_e)_{e\in
  \vec{E}(\sg)})$ and $(\sg',(W'_e)_{e\in \vec{E}(\sg')})$ be at
  distance $1$ if $\sg\neq \sg'$, and at distance $\max_{e\in
    \vec{E}(\sg)}\mathrm{dist}(W_e,W'_e)$ if $\sg=\sg'$.

In order to make clear distinctions between quantities like
$\zeta(W_e),\zeta(W_e(s))$, we will adopt the following notation in
the sequel. For $e\in \vec{E}(\sg)$, we let 
$$M_e=W_e(0)\in \CC(\R)\, ,\quad r_e=\zeta(M_e)\, ,\quad 
\sigma_e=\zeta(W_e)\, .$$

The measure that will play a central role is a {\em continuum measure}
$\CLM^{(k+1)}$ on labeled maps. To define it, we first need to
describe the continuum analogs of the walks and discrete snakes
considered in the previous section.

\subsubsection{Bridge measures}\label{sec:bridge-measures}

In the sequel, we let $(X_t,t\in [0,\zeta(X)])$ be the canonical
process on the space $\ov{\CC}(\R)$ of $\R$-valued continuous functions
with finite or infinite duration. The infimum process of $X$ is the
process $\un{X}$ defined by 
$$\un{X}_t=\inf_{0\leq s\leq t}X_s\, ,\qquad 0\leq t\leq \zeta(f)\,
.$$ For every $x\in \R$, we also let $T_x=\inf\{t\geq 0:X_t=x\}\in
[0,\infty]$ be the first hitting time of $x$ by $X$. 

The building blocks of $\CLM^{(k+1)}$ are going to be several paths
measures.  We let $\P_x$ be the law of standard $1$-dimensional
Brownian motion started from $x$, and for $y<x$, we let
$\E_x^{(y,\infty)}$ be the law of $(X_t,0\leq t\leq T_y)$ under
$\P_x$, that is, Browian motion killed at first exit time of
$(y,\infty)$.

Next, let $\P_{x\to y}^t$ be the law of the $1$-dimensional Brownian
bridge from $x$ to $y$ with duration $t$:
$$\P_{x\to y}^t(\cdot)=\P_x( (X_s)_{0\leq s\leq t}\in \cdot \, |\,
X_t=y)\, .$$ A slick way to properly define this singular conditioned
measure is to let $\P_{x\to y}^t$ be the law of $(X_s+(y-X_t)s/t,0\leq
s\leq t)$ under $\P_x$, see \cite[Chapter I.3]{revyor}. We also let
$\ov{\P}_x^{t}$ be the law of the {\em first-passage
  Brownian bridge} from $x$ to $0$ with duration $t$, defined formally
by 
$$\ov{\P}_x^{t}(\cdot)=\P_x((X_s,0\leq s\leq t)\in \cdot\, |\,
T_0=t)\, .$$ Regular versions for this singular conditioning can be
obtained using space-time Doob $h$-transforms. We refer to
\cite{bettinelli10} for a recent and quite complete treatment of
first-passage bridges and limit theorems for their discrete
versions, which will be helpful to us later on.

For $t>0$ let
\begin{eqnarray*}
p_t(x,y)&=&\frac{\P_x(X_t\in \d y)}{\d y}=\frac{1}{\sqrt{2\pi
    t}}e^{-\frac{(x-y)^2}{2t}} \qquad 
x,y\in \R\, ,\\
\ov{p}_t(x,0)&=&\frac{\P_x(T_0\in \d t)}{\d t}=\frac{x}{\sqrt{2\pi t^3}}e^{-\frac{x^2}{2t}}
\qquad x>0\, .
\end{eqnarray*}
The {\em bridge measure} from $x$ to $y$ is then the $\sigma$-finite
measure $\B_{x\to y}$ defined by 
$$\B_{x\to y}(\d X)=\int_0^\infty \d t\, p_t(x,y)\P_{x\to y}^t(\d
X)\, .$$ We also let $\B_{x\to y}^+$ be the restriction of $\B_{x\to
  y}$ to the set $\CC^+(\R)=\{f\in \CC(\R):f\geq 0\}$.  The reflexion
principle entails that $\P_{x\to y}^t(\CC^+(\R))=1-e^{-2xy/t}$, and it
is a simple exercise to show that for $x,y>0$, letting
$p_t^+(x,y)=p_t(x,y)-p_t(x,-y)=\P_x(X_t\in \d y,\un{X}_t\geq 0)/\d y$,
\begin{equation}
  \label{eq:20}
  \B^+_{x\to y}(\d X)=\int_0^\infty\d t\, p_t^+(x,y)\P_{x\to y}^t(\d
  X\, |\, \CC^+(\R))\, ,
\end{equation}

We now state two path decomposition formulas for the measures
$\B_{x\to y}$ that will be useful for our later purposes.  If $\mu$ is
a measure on $\CC$, we let $\widehat{\mu}$ be the image measure of
$\mu$ under the time-reversal operation $f\mapsto \widehat{f}$.  If
$\mu$ and $\mu'$ are two probability measures on $\CC$ such that
$\mu(X(\zeta(X))=z)=\mu'(X(0)=z)=1$ for some $z\in \R$, then we let
$\mu\bowtie\mu'$ be the image measure of $\mu\otimes \mu'$ under the
concatenation operation $(f,g)\mapsto f\bowtie g$ from $\CC(\R)^2$ to
$\CC(\R)$, where
$$(f\bowtie g) (t)=\left\{\begin{array}{cl}
f(t) & \mbox{ if } 0\leq t\leq \zeta(f)\\
g(t-\zeta(f)) & \mbox{ if } \zeta(f)\leq t\leq \zeta(f)+\zeta(g)\, .
  \end{array}\right.$$
The {\em agreement formula} of \cite[Corollary 3]{PiYo96} states that
\begin{equation}\label{eq:27}
  \mathbb{B}_{x\to y}(\d X)=\int_{-\infty}^{x\wedge y}\d z\, 
  \big(\E_x^{(z,\infty)}\bowtie \widehat{\E}_y^{(z,\infty)}\big)(\d X)\,
  . 
\end{equation}
In particular,
$$\mathbb{B}^+_{x\to y}(\d X)=\int_0^{x\wedge y}\d z\, 
\big(\E_x^{(z,\infty)}\bowtie \widehat{\E}_y^{(z,\infty)}\big)(\d X)\,
.$$ This decomposition should be seen as one of the many versions of
Williams' decompositions formulas for Brownian paths:
here, the variable $z$ plays the role of the minimum of the generic
path $X$.

A second useful path decomposition states that $x,y\in \R$,
\begin{equation}
  \label{eq:28}
\int_{\R}\d z\, (\mathbb{B}_{x\to z}\bowtie\mathbb{B}_{z\to y})(\d X)=
\zeta(X)\mathbb{B}_{x\to y}(\d X)\, . 
\end{equation}
To see this, write
$$\zeta(X)\mathbb{B}_{x\to y}(\d X)=\int_0^\infty\d r\, 
p_r(x,y)\int_0^r\d s\, \P_{x\to y}^r(\d X)\, ,$$ then note that for every
$r\geq s\geq 0$,
$$\P_{x\to y}^r(X_s\in \d z)=\frac{p_s(x,z)p_{r-s}(z,y)}{p_r(x,y)}\d z\,
,$$ and use the fact that given $X_s=z$, the paths $(X_u,0\leq u\leq
s)$ and $(X_{s+u},0\leq u\leq r-s)$ under $\P_{x\to y}^r$ are
independent Brownian bridges, by the Markov property.

\subsubsection{Brownian snakes}\label{sec:snakes}

Let $W$ be the canonical process on $\ov{\CC}(\CC(\R))$. That is to
say, for every $s\in [0,\zeta(W)]$ (or just $s\geq 0$ if
$\zeta(W)=\infty$), $W(s)$ is an element of $\CC(\R)$. For simplicity,
we let $\zeta_s=\zeta(W(s))$ be the duration of this path, and let
$W(s,t)=W(s)(t)$ for $0\leq t\leq \zeta_s$. The process
$(\zeta_s,0\leq s\leq \zeta(W))$ is called the {\em lifetime process}
of $W$. In order to clearly distinguish the duration of $W$ with that
of $W(s)$ for a given $s$, we will rather denote the duration
$\zeta(W)$ of $W$ by the letter $\sigma(W)$.

Let us describe the law of Le Gall's Brownian snake (see
\cite{legall99} for an introduction to the subject). Conditionally
given the lifetime process $(\zeta_s)$, the process $W$ under the
Brownian snake distribution is a non-homogeneous Markov process with
the following transition kernel. Given $W(s)=(w(t),0\leq t\leq
\zeta_t)$, the law of the path $W(s+s')$ is that of the path
$(w'(t),0\leq t\leq \zeta_{s+s'})$ defined by
$$w'(t)=\left\{\begin{array}{cl}
W(s,t) & \mbox{ if } 0\leq  t\leq \check{\zeta}_{s,s+s'}\\
W(s, \check{\zeta}(s,s+s'))+B_{t- \check{\zeta}(s,s+s')} & \mbox{ if } \check{\zeta}_{s,s+s'}\leq
t\leq \zeta_{s+s'}
\end{array}\right.\, ,
$$
where $(B_t,t\geq 0)$ is a standard Brownian motion independent of
$\zeta$, and 
$$\check{\zeta}_{s,s+s'}=\inf_{s\leq u\leq s+s'}\zeta_u\, .$$
We let $\Q_w$ be the law of the process $W$ started from the path
$w\in \CC$, and driven by a Brownian motion started from $\zeta(w)$
and killed at first hitting of $0$ (i.e.\ a process with law
$\E_{\zeta(w)}^{(0,\infty)}$).

For our purposes, the key property of Brownian snake will be the
following representation using Poisson random measures, which can be
found for instance in \cite{legall99,legweill}. Namely, let $w\in
\CC(\R)$.  Recall that $(\zeta_s,0\leq s\leq \sigma(W))$ is the lifetime
process driving the canonical process $W$, and that under $\Q_w$ we
have $W(0)=w$ and $\zeta_0=\zeta(w)$. For $0\leq r\leq \zeta_0$,
define $\Sigma_r=\inf\{s\geq 0:\zeta_s=\zeta_0-r\}$, so in particular
$\sigma=\Sigma_{\zeta_0}$, and the process $(\Sigma_r,0\leq r\leq
\zeta_0)$ is non-decreasing and right-continuous.  For every $r\in
[0,\zeta_0]$ such that $\Sigma_r>\Sigma_{r-}$, let
\begin{equation}\label{eq:17}
W^{(r)}(s)=(W(\Sigma_{r-}+s,t+r),0\leq t\leq \zeta_{\Sigma_{r-}+s}-r)\,
,\qquad 0\leq s\leq \Sigma_r-\Sigma_{r-}\, .
\end{equation}
In this way, every
$W^{(r)}$ is an element of $\CC(\CC(\R))$.  Then, under $\Q_w$, we
consider the measure
$$\mathcal{M}_w=\sum_{0\leq r\leq \zeta(w)}
\delta_{(r,W^{(r)})}\ind_{\{T_r>T_{r-}\}}\, .$$ Lemma 5 in
\cite[Chapter V]{legall99} states that under $\Q_w$, the measure
$\mathcal{M}_w$ is a Poisson random measure on $[0,\zeta(w)]\times
\CC$, with intensity measure given by
\begin{equation}
  \label{eq:21}
2\, \d r\, \ind_{[0,\zeta(w)]}(r)\, \N_{w(r)}(\d W)\, .
\end{equation}
Here, the measure $\N_x$ is called the {\em Itô excursion measure} of
the Brownian snake started at $x$. It is the $\sigma$-finite ``law''
of the Brownian snake started from the (trivial) path $w$ with
$w(0)=x$ and $\zeta(w)=0$, and driven by a trajectory which is a
Brownian excursion under the Itô measure of the positive excursions of
Brownian motion $n(\d \zeta)$. See \cite[Chapter XII]{revyor} for the
properties of the excursion measure $n$, which is denoted by $n_+$ in
this reference. It is important that we fix the normalization of this
measure so that the factor $2$ in \eqref{eq:21} makes sense, and we
choose it so that $n(\sup \zeta>x)=(2x)^{-1}$ for every $x>0$.

\subsubsection{The measure $\CLM^{(k+1)}$}\label{sec:measure-clmk+1}

For every dominant scheme $\sg\in \bS_d^{(k+1)}$, we let $\lambda_\sg$
be the measure on $\R^{V(\sg)}$ defined by 
$$\lambda_\sg(\d (\ell_v)_{v\in V(\sg)})=\prod_{v\in
  V_N(\sg)}\delta_0(\d \ell_v)\prod_{v\in V_I(\sg)}\d
\ell_v\ind_{\{\ell_v>0\}}\prod_{v\in V_O(\sg)}\d \ell_v\, ,$$ called
the Lebesgue measure on admissible labelings of $\sg$.  We define the
continuum measure $\CLM^{(k+1)}$ on $\bCLM^{(k+1)}$ by
\begin{eqnarray}
  \CLM^{(k+1)}(\d(\sg,(W_e)_{e\in \vec{E}(\sg)}))
  &=&\S_d^{(k+1)}(\d
  \sg)\int_{\R^{V(\sg)}}\lambda_\sg(\d (\ell_v)_{v\in V(\sg)}) \label{eq:19}\\
  & & \times \prod_{e\in
    E_N(\sg)}\int_{\CC^+(\R)}\E^{(0,\infty)}_{\ell_{e_-}}(\d M_e)\Q_{\widehat{M}_e}(\d
  W_e)\Q_{M_e}(\d W_{\ov{e}})\nonumber\\
  & & \times \prod_{e\in
    E_I(\sg)}\int_{\CC^+(\R)}\B^+_{\ell_{e_-}\to \ell_{e_+}}(\d M_e)\Q_{\widehat{M}_e}(\d
  W_e)\Q_{M_e}(\d W_{\ov{e}})\nonumber\\
  & & \times \prod_{e\in E_O(\sg)}\int_{\CC(\R)}\B_{\ell_{e_-}\to
    \ell_{e_+}}(\d M_e)
  \Q_{\widehat{M}_e}(\d
  W_e)\Q_{M_e}(\d W_{\ov{e}})\, .\nonumber
\end{eqnarray}
where $\S_d^{(k+1)}$ is the counting measure on $\bS_d^{(k+1)}$, and
$\widehat{f}(t)=f(\zeta(f)-t),0\leq t\leq \zeta(f)$ is the
time-reversed function obtained from $f\in \CC(\R)$. Note that, since
$E_N(\sg),E_I(\sg)$ and $E_O(\sg)$ partition $E(\sg)$, and since we
have fixed an orientation convention for the edges in these sets, the
oriented edges $e,\ov{e}$ exhaust $\vec{E}(\sg)$ when $e$ varies along
$E_N(\sg),E_I(\sg),E_O(\sg)$.

We also define a measure $\CLM^{(k+1)}_1$, which is rougly speaking a
conditioned version of $\CLM^{(k+1)}$ given $\sum_{e\in
  \vec{E}(\sg)}\sigma_e=1$. Contrary to $\CLM^{(k+1)}$, this is a
probability measure on $\bCLM^{(k+1)}$. Its description is more
elaborate than $\CLM^{(k+1)}$, because it does not have a simple
product structure. It is more appropriate to start by defining the
trace of $\CLM^{(k+1)}_1$ on $(\sg,(\ell_v),(r_e))$, i.e.\ its
push-forward by the mapping $(\sg,(W_e))\mapsto (\sg,(\ell_v),(r_e))$,
where $r_e=\zeta(W_e(0))$ and $\ell_v=W_e(0,0)$ whenever $v=e_-$. This trace is
\begin{equation}\label{eq:38}
\frac{1}{\Upsilon^{(k+1)}}\S_d^{(k+1)}(\d\sg) \, 
\lambda_\sg((\d \ell_v)_{v\in V(\sg)}) \Big(\prod_{e\in E(\sg)}\d r_e\, 
p^{(e)}_{r_e}(\ell_{e_-},\ell_{e_+}) \Big) \ov{p}_1(2r,0)\, ,
\end{equation}
where
$\Upsilon^{(k+1)}\in (0,\infty)$ is the normalizing constant making it a
probability distribution, and $p^{(e)}_r(x,y)$ is defined by 
\begin{itemize}
\item $p^{(e)}_r(x,y)=p_r(x,y)$ if $e\in E_O(\sg)$, 
\item $p^{(e)}_r(x,y)=p^+_r(x,y)$ if $e\in E_I(\sg)$, and 
\item $p^{(e)}_r(x,0)=\ov{p}_r(x,0)$ if $e\in E_N(\sg)$ (entailing
  automatically $\ell_{e_+}=0$).
\end{itemize}
Finally, in the last displayed expression, we let $r=\sum_{e\in
  E(\sg)}r_e$.  Then, conditionally given $(\sg,(\ell_v),(r_e))$, the
processes $(M_e,e\in \check{E}(\sg))$ are independent, and
respectively
\begin{itemize}
\item $M_e$ has law $\P^{r_e}_{\ell_{e_-}\to \ell_{e_+}}$ if $e\in
  E_O(\sg)$, 
\item $M_e$ has law $\P^{r_e}_{\ell_{e_-}\to \ell_{e_+}}(\cdot\, |\, \CC^+(\RR))$ if $e\in
  E_I(\sg)$
\item $M_e$ has law $\ov{\P}^{r_e}_{\ell_{e_-}}$ if $e\in
  E_N(\sg)$. 
\end{itemize}
As usual, the processes $(M_{\ov{e}},e\in \check{E}(\sg))$ are defined
by time-reversal: $M_{\ov{e}}=\widehat{M}_e$.  Finally, conditionally
given $(M_e,e\in E(\sg))$, the processes $(W_e,W_{\ov{e}},e\in
\check{E}(\sg))$ are independent Brownian snakes respectively started
from $\widehat{M}_e,M_e$, conditioned on the event that the sum
$\sigma$ of their durations is equal to $1$. 

This singular conditioning is obtained in the following way.  First
consider a Brownian snake $(W^\circ(s),0\leq s\leq 1)$ with lifetime
process given by a first-passage bridge $(\zeta(s),0\leq s\leq 1)$
from $2r=2\sum_{e\in \vec{E}(\sg)}r_e$ to $0$ with duration $1$, and
such that $W^\circ(0)$ is the constant path $0$ with duration $2r$.
Let $e_1,\ldots,e_{4k-3}$ be an arbitrary enumeration of
$\vec{E}(\sg)$, and let $r_i=r_{e_1}+\ldots+r_{e_i}$ for every $i\in
\{0,1,\ldots,4k-3\}$, with the convention $r_0=0$. Then, let
$\kappa_i=\inf\{s\geq 0:\zeta(s)=2r-r_i\}$, and for $1\leq i\leq
4k-3$,
\begin{equation}\label{eq:46}
W^\circ_{e_i}(s,t)=W^\circ(s+\kappa_{i-1},t+2r-r_i)\, ,
\qquad 0\leq s\leq \kappa_i-\kappa_{i-1}\, ,\quad 0\leq t\leq
\zeta(s+\kappa_{i-1})-\zeta(\kappa_i)\, .
\end{equation}
%
The processes $(W^\circ_e,e\in \vec{E}(\sg))$ are then independent
Brownian snakes started from the constant trajectories $0$ with
respective durations $r_e$, and conditioned on having total duration
$\sigma=\sum_e\sigma_e$ equal to $1$.  We finally let $\zeta_e^\circ$
be the lifetime process of $W^\circ_e$, and
\begin{equation}
  \label{eq:47}
  W_e(s,t)=\widehat{M}_e(t)+W^\circ_e(s,t)\, ,\qquad 0\leq s\leq
  \sigma_e\, ,\quad 0\leq t\leq \inf_{0\leq u\leq s}\zeta^\circ_e(u)\, ,
\end{equation}
that is, we change the initial value of $W^\circ_e$ to
$\widehat{M}_e$.  The different processes $(W_e,e\in \vec{E}(\sg))$
are then conditionally independent snakes started respectively from
$\widehat{M}_e$, conditioned on $\sigma=1$. Hence, the snakes
$(W_e,e\in \vec{E}(\sg))$ can be obtained by ``cutting into bits of
initial lengths $r_e$'' a snake with lifetime process having law
$\E_{2r}^{(0,\infty)}[\cdot\, |\, T_0=1]$ and started from the
constant zero trajectory, to which we superimpose (independently) the
initial trajectories $\widehat{M}_e$.

The relation between $\CLM^{(k+1)}$ and $\CLM^{(k+1)}_1$ goes as
follows. For every $c>0$, define a scaling operation $\Psi^{(k+1)}_c$ on
$\bCLM^{(k+1)}$, sending $(\sg,(W_e)_{e\in \vec{E}(\sg)})$ to
$(\sg,(W_e^{[c]})_{e\in \vec{E}(\sg)})$, where
$$W_e^{[c]}(s,t)=c^{1/4}W_e(s/c,t/c^{1/2})
\, ,\qquad 0\leq s\leq c \sigma_e\, ,\quad 0\leq t\leq
c^{1/2}\zeta_e(s/c) \, .$$ Then, we let $\CLM^{(k+1)}_c$ be the image
measure of $\CLM^{(k+1)}_1$ by $\Psi^{(k+1)}_c$. Sometimes we will abuse
notation and write $\Psi^{(k+1)}_c((W_e)_{e\in \vec{E}(\sg)})$ instead of
$(W^{[c]}_e)_{e\in \vec{E}(\sg)}$.

\begin{prp}\label{sec:cont-meas-label-1}
It holds that 
$$\CLM^{(k+1)}=\Upsilon^{(k+1)}\int_0^\infty\d \sigma\, \sigma^{k-9/4}\,
\CLM^{(k+1)}_\sigma\, .$$
\end{prp}

\proof The idea is to disintegrate formula \eqref{eq:19} with respect
to $\sigma=\sum_{e\in \vec{E}(\sg)}\sigma_e$. Note that conditionally
given $\sg,(M_e,e\in E(\sg))$, the snakes $(W_e)_{e\in \vec{E}(\sg)}$
are independent Brownian snakes started respectively from the
trajectories $M_{\ov{e}}$, so that their lifetime processes are
independent processes with respective laws $\E_{r_e}^{(0,\infty)}(\d
\zeta_e)$. In particular, the lifetime of $W_e$ is an independent
variable with same law as $T_{-r_e}$, the first hitting time of $-r_e$
under $\P_0$.  From this, we see that (still conditionally given
$\sg,(M_e,e\in E(\sg))$), the total lifetime $\sigma$ has same
distribution as $T_{-2r}$ under $\P_0$, where $r=\sum_{e\in
  E(\sg)}r_e$, and this has a distribution given by
$\ov{p}_\sigma(2r,0)\d \sigma$ (in several places in this proof, we
will not differentiate the random variable $\sigma$ from 
its generic value, and will do a similar abuse of notation for the
elements $(\ell_v)$ and $(r_e)$, not to introduce new
notation). Consequently, we obtain that the trace of $\CLM^{(k+1)}$
on $(\sg,(\ell_v),(r_e),\sigma)$ equals
\begin{eqnarray*}\S_d^{(k+1)}(\d \sg)\lambda_{\sg}(\d(\ell_v)_{v\in
    V(\sg)}) &&\!\!\!\!\!\prod_{e\in
    E_N(\sg)}\int_{\CC^+(\R)}\E^{(0,\infty)}_{\ell_{e_-}}(\zeta(M_e)\in
  \d r_e)\\
&&   \times \prod_{e\in
    E_I(\sg)}\int_{\CC^+(\R)}\B^+_{\ell_{e_-}\to
    \ell_{e_+}}(\zeta(M_e)\in \d r_e)\\
  && \times \prod_{e\in E_O(\sg)}\int_{\CC(\R)}\B_{\ell_{e_-}\to
    \ell_{e_+}}(\zeta(M_e)\in \d r_e)\\
&& \times
\ov{p}_\sigma(2r,0)\d \sigma\\
=\S_d^{(k+1)}(\d \sg)\lambda_{\sg}(\d(\ell_v)_{v\in
    V(\sg)}) &&\!\!\!\!\!\Big( \prod_{e\in E(\sg)}\d r_e\,
  p^{(e)}_{r_e}(\ell_{e_-},\ell_{e_+})\Big)\d \sigma\, \ov{p}_{\sigma}(2r,0)\, .
\end{eqnarray*}
Conditionally given $(\sg,(\ell_v),(r_e),\sigma)$, the path $M_e$ has
law $\P^{r_e}_{\ell_{e_-}\to \ell_{e_+}}, \P^{r_e}_{\ell_{e_-}\to
  \ell_{e_+}}(\cdot\, |\, \CC^+(\R))$ or $\ov{\P}^{r_e}_{\ell_{e_-}}$
according to whether $e$ belongs to $E_O(\sg),E_I(\sg)$ or $E_N(\sg)$,
and these paths are independent. Finally, the paths $W_e$ are
independent Brownian snakes respectively starting from $M_{\ov{e}}$,
conditioned on the sum of their durations being $\sigma$. These paths
can be defined by applying the scaling operator $\Psi^{(k+1)}_\sigma$
to a family of independent Brownian snakes started from the paths
$(\sigma^{-1/4}M_e(\sigma^{1/2}t),0\leq t\leq \sigma^{-1/2}r_e)$, and
conditioned on the sum of their durations being $1$.

Let us change variables $r_e'=\sigma^{-1/2}r_e$ and
$\ell_v'=\sigma^{-1/4}\ell_v$ and let $r'=\sum_{e\in E(\sg)}r'_e$. We
have obtained that $\CLM^{(k+1)}$ is alternatively described by a
trace on $(\sg,(\ell_v),(r_e),\sigma)$ equal to
$$\S_d^{(k+1)}(\d \sg)
\lambda_{\sg}(\d(\sigma^{1/4}\ell'_v)_{v\in V(\sg)})\Big( \prod_{e\in
  E(\sg)}\sigma^{1/2}\d r'_e\,
p^{(e)}_{\sigma^{1/2}r'_e}(\sigma^{1/4}\ell'_{e_-},\sigma^{1/4}\ell'_{e_+})\Big)\d
\sigma\, \ov{p}_{\sigma}(2r'\sigma^{1/2},0)\, ,$$ and conditionally
given these quantities, the processes $(M_e)$ are chosen as in the
previous paragraph, while the snakes $(W_e)$ are the image under
$\Psi^{(k+1)}_{\sigma}$ of independent snakes respectively started from
$M_{\ov{e}}$, conditioned on the sum of their durations being $1$.  We
then use the fact that for every $r,c>0$ and every $x,y\in \R$ such
that the following expressions make sense,
$$p_{cr}(c^{1/2}x,c^{1/2}y)=
c^{-1/2} p_r(x,y)\, ,
\qquad
p^+_{cr}(c^{1/2}x,c^{1/2}y)=
c^{-1/2} p^+_r(x,y)\, , 
$$
and 
$$
\ov{p}_{cr}(c^{1/2}x,0)=c^{-1} p^+_r(x,0)\, ,$$ which is a simple
consequence of Brownian scaling. Finally, for dominant schemes we have
$3k-2$ vertices, among which $k$ are in $V_N(\sg)$ do not contribute
in $\lambda_{\sg}(\d(\ell_v))$, and $4k-3$ edges among which $2k$ are
in $E_N(\sg)$, yielding that the trace of $\CLM^{(k+1)}$ on
$(\sg,(\ell_v),(r_e),\sigma)$ is
$$\sigma^{\frac{2(k-1)}{4}+\frac{2k-3}{4}-1}\, \d
\sigma\, \S_d^{(k+1)}(\d \sg) \lambda_{\sg}(\d (\ell_v)_{v\in
  V(\sg)})\Big( \prod_{e\in E(\sg)}\d r_e\,
p^{(e)}_{r_e}(\ell_{e_-},\ell_{e_+})\Big)\, \ov{p}_1(2r,0)\, ,$$ and
we recognize $\Upsilon^{(k+1)}\sigma^{k-9/4}\d \sigma$ times the trace
of $\CLM^{(k+1)}_1$ on $(\sg,(\ell_v),(r_e))$. Since
$\CLM^{(k+1)}_\sigma$ is the image measure of $\CLM_1^{(k+1)}$ by
$\Psi^{(k+1)}_\sigma$, we have finally obtained that the trace of
$\CLM^{(k+1)}$ on the variable $\sigma$ is
$\Upsilon^{(k+1)}\sigma^{k-9/4}\d \sigma$, and that
$\CLM_{\sigma}^{(k+1)},\sigma>0$ are conditional measures of
$\CLM^{(k+1)}$ given $\sigma$, as wanted.  \cq

\subsection{Limit theorems}\label{sec:limit-theorems}

For every $n\geq 1$, we define a new scaling operator on $\bCLM^{(k+1)}$
by $\psi^{(k+1)}_n(\sg,(W_e)_{e\in \vec{E}(\sg)})=(\sg,(W^{\{n\}}_e)_{e\in
  \vec{E}(\sg)})$, where
$$W^{\{n\}}_e(s,t)=\Big(\frac{9}{8n}\Big)^{1/4}W_e(2ns,\sqrt{2n}t)\,
,\qquad 0\leq s\leq \frac{\sigma_e}{2n}\, ,\quad 0\leq t\leq
\frac{1}{\sqrt{2n}}\zeta_e(2ns)\, .$$ So this is almost the same as
$\Psi_{1/2n}$, except that we further multiply the labels in $W_e$ by
$\sqrt{3/2}=(9/4)^{1/4}$. Recall that $\LM^{(k+1)}_n$ is the counting
measure over $\bLM^{(k+1)}_n$, we view it also as a measure on
$\bCLM^{(k+1)}$ by performing the decomposition of Proposition
\ref{sec:decomp-label-maps-1}. Similar abuse of notation will be used
in the sequel. 

\begin{prp}\label{sec:limit-theorems-1}
  We have the following weak convergence of finite measures on
  $\bCLM^{(k+1)}$:
$$\Big(\frac{9}{2}\Big)^{1/4}\frac{
\psi_{n}^{(k+1)}{}_*\LM_n^{(k+1)}}{6^k\cdot12^nn^{k-5/4} }\build\longrightarrow_{n\to\infty}^{}\Upsilon_{k+1}\CLM^{(k+1)}_1\,
  .$$
\end{prp}

This is proved in a very similar way to
\cite{miertess,bettinelli10}, but one has to pay extra care in
manipulating elements of $\bLM^{(k+1)}$,
because of the required positivity of $\ell_v$ and $M_e$ when $v\in
V_I(\sg)$ and $e\in E_N(\sg)\cup E_I(\sg)$.
We start with a preliminary observation that justifies our definition
of {\em dominant} schemes. 

\begin{lmm}
  \label{sec:limit-theorems-2} It holds
  that $\LM^{(k+1)}_n(\{(\sg,(W_e)_{e\in \vec{E}(\sg)}):\sg\in
  \bS^{(k+1)}_d\})=O(12^nn^{k-5/4})$, while 
$\LM^{(k+1)}_n(\{(\sg,(W_e)_{e\in \vec{E}(\sg)}):\sg\notin
  \bS^{(k+1)}_d\})=O(12^nn^{k-3/2})$
\end{lmm}

\proof
Let $\sg$ be an element of $\bS^{(k+1)}$, that is a scheme with $k+1$
faces. Then 
\begin{equation}\label{eq:36}
\LM_n^{(k+1)}(\{(\sg,(W_e)_{e\in
  \vec{E}(\sg)}):(W_e)_{e\in \vec{E}(\sg)} \mbox{ a compatible family
  of discrete snakes}\})
\end{equation}
is just the number of elements $(\bm,\bl)$
in $\bLM^{(k+1)}$ that induce the scheme $\sg$ via the construction of
Proposition \ref{sec:decomp-label-maps-1}, and which have $n$ edges in
total. By the discussion of Section \ref{sec:decomp-label-maps}, an
element $(\sg,(W_e)_{e\in \vec{E}(\sg)})$ can be viewed as a walk
network $(M_e,e\in E(\sg))$ and a family of labeled forests
$(F_e,L_e)_{e\in \vec{E}(\sg)}$ compatible with this walk
network. Now, once the walks $(M_e,e\in E(\sg))$ are determined, we
know that for $e\in \vec{E}(\sg)$, the forest $F_e$ has $r_e$ trees,
where $r_e$ is the duration of $M_e$. The labels of the floor vertices
of $F_e$ are $M_e(0),M_e(1),\ldots,M_e(r_e)$. For a given $F_e$, there
are exactly $3^{n_e}$ possible labelings compatible with the labeling
of the floor vertices, where $n_e$ is the number of edges in the
forest that are distinct from the floor edges, coming from the fact
that the label difference along an edge of each tree belongs to
$\{-1,0,1\}$. Therefore, if we let $r=\sum_{e\in E(\sg)}r_e$, then by
concatenating the forests $F_e,e\in \vec{E}(\sg)$ in some given order
that is fixed by convention, we obtain a forest with $2r$ trees and a
total of $n-r$ edges (distinct from the floor edges) by
\eqref{eq:18}. Once the labeling of the roots is fixed, there are
$3^{n-r}$ different labelings for these forests. This construction can
be easily inverted: Starting from a forest with $2r$ trees and $n-r$
edges (distinct from the floor edges), we can reconstruct
$(F_e,L_e)_{e\in \vec{E}(\sg)}$.

From this, and a classical counting result for plane forests
\cite{pitmancsp02}, we deduce that for a given walk network $(M_e,e\in
\check{E}(\sg))$, there are exactly
$$3^{n-r}\frac{2r}{2n}\binom{2n}{n-r}$$ 
labeled forests $(F_e,L_e)_{e\in \vec{E}(\sg)}$ compatible with this
walk network. Since, in the decomposition of Proposition
\ref{sec:decomp-label-maps-1}, we still have to select one of the
oriented edges in the forests (including floor edges), we will obtain
an extra factor of $2n$ in the end. At this point, we have obtained
that \eqref{eq:36} equals
$$2\sum_{(M_e,e\in
  E(\sg))}3^{n-r}r\binom{2n}{n-r}\, ,$$ where the sum is over all walk
networks compatible with the scheme $\sg$.

Let $\mathcal{W}(a,b;r)$ be the number of Motzkin
walks with duration $r$, starting at $a$ and ending at $b$. Likewise,
we let $\mathcal{W}^+(a,b;r)$ be the number of such paths that are
strictly positive, except perhaps at their endpoints (so that
$\#\mathcal{W}^+(0,0;1)=1$ for instance).  
Then the formula for \eqref{eq:36} becomes 
\begin{equation}
\label{eq:45}
2\sum_{(\ell_v)}\sum_{(r_e)}3^{n-r}r\binom{2n}{n-r}\prod_{e\in
  \check{E}(\sg)}\mathcal{W}^{(e)}(\ell_{e_-},\ell_{e_+};r_e)\, ,
\end{equation}
where the first sum is over all admissible labelings for $\sg$, the
second sum is over all edge-lengths on $\sg$, and the superscript
$(e)$ accounts for the constraint on the path $M_e$, namely,
\begin{itemize}
\item 
 $\mathcal{W}^{(e)}(a,b;r)=\mathcal{W}(a,b;r)$ if $e\in E_O(\sg)$, 
\item
$\mathcal{W}^{(e)}(a,b;r)=\mathcal{W}^+(a,b;r)$ if $e\in
E_I(\sg)$, 
\item 
$\mathcal{W}^{(e)}(a,0;r)=\mathcal{W}^+(a,0;r)$ if $e\in E_N(\sg)
\setminus E_T(\sg)$
\item
$\mathcal{W}^{(e)}(a,0;r)=\mathcal{W}^+(a+1,0;r+1)$ if $e\in
E_T(\sg)$. 
\end{itemize}
To explain the last point, recall that the walks $M_e$ indexed by the
thin edges are non-negative and finish at $0$, so we can turn it into
a walk taking positive values except at the last point, by translating labels by $1$
and adding an extra ``virtual'' $-1$ step in the end. 

Note that $q_r(a,b):=3^{-r}\mathcal{W}(a,b;r)$ (resp.\
$q_r^+(a,b):=3^{-r}\mathcal{W}^+(a,b;r)$) is the probability that a uniform
Motzkin walk with $r$ steps started at $a$ finishes at $b$ (resp.\
without taking non-negative values, except possibly at the start and
end point). 

At this point, we will need the following 
consequence of the local limit theorem of \cite{petrov75}, stating
that for every $a,b\in \Z$ and $p\geq 1$, 
there exists a finite constant
$C>0$ depending only on $p$ such that for every $r\geq 1$, 
\begin{equation}
  \label{eq:39}
  \sqrt{r}q_r(a,b)\leq \frac{C}{1+\Big|\frac{b-a}{\sqrt{r}}\Big|^p}\, .
\end{equation}
Likewise, by viewing $4^{-n}\binom{2n}{n-r}$ as the probability that a
simple random walk attains $2r$ in $2n$ steps, a similar use of the
local limit theorem allows to show
\begin{equation}
  \label{eq:42}
  \frac{\sqrt{2n}}{4^n}\binom{2n}{n-r}\leq
  \frac{C}{1+\Big(\frac{r}{\sqrt{n}}\Big)^p}\, .
\end{equation}
On the other hand, the reflection principle entails that for $a,b,r>0$, 
\begin{equation}
  \label{eq:40}
  q_r^+(a,b)=q_r(a,b)-q_r(a,-b)\, ,
\end{equation}
and the cyclic lemma \cite{pitmancsp02} entails that for $a,r>0$, 
\begin{equation}
  \label{eq:41}
  q_r^+(a,0)=\frac{a}{r}q_r(a,0)\, ,\quad \mbox{ and }\qquad
  q^+_r(0,0)=
\frac{1}{3}\ind_{\{r=1\}}+\frac{1}{r-1}q_{r-1}(1,0)\ind_{\{r>1\}}\, .
\end{equation}
These considerations imply that for a fixed $p$, \eqref{eq:36} is
bounded from above by
\begin{eqnarray}
  \lefteqn{C\, 
    12^n\sum_{(\ell_v)}\sum_{(r_e)}\frac{\frac{r}{\sqrt{n}}}{1+
      \Big(\frac{r}{\sqrt{n}}\Big)^{p+1}}}\label{eq:43}
  \\
  & & \times 
  \prod_{e\in E_F(\sg)}\frac{1}{\sqrt{r_e}}\cdot\frac{1}{1+\Big(\frac{\ell_{e_+}-\ell_{e_-}}{\sqrt{r_e}}\Big)^2}
  \prod_{e\in
    E^{(1)}_N(\sg)}  \frac{1}{r_e}\cdot\frac{\frac{\ell_{e_-}+\ind_{\{e\in
        E_T(\sg)\}}}{\sqrt{r_e+\ind_{\{e\in E_T(\sg)\}}
      }}}{1+\Big|\frac{\ell_{e_-}+\ind_{\{e\in
        E_T(\sg)\}}}{\sqrt{r_e+\ind_{\{e\in E_T(\sg)\}}}}\Big|^3}
  \prod_{e\in E^{(0)}_N(\sg)}\frac{1}{r_e^{3/2}}\, ,\nonumber
\end{eqnarray}
where we let $E_N^{(0)}(\sg)$ be the set of edges in $E_N(\sg)$ with
both extremities in $V_N(\sg)$, $E_N^{(1)}(\sg)=E_N(\sg)\setminus
E_N^{(0)}(\sg)$, and finally 
$$E_F(\sg)=E(\sg)\setminus
E_N(\sg)=E_I(\sg)\cup E_O(\sg)$$ (here the subscript $F$ stands for
{\em free}). Of course, the constant $C$
above depends only on $p$, but not on $n$. As before, the sums ae over
all admissible labelings and edge-lengths. At this point, since we are
only interested in taking upper-bounds, we may (and will) in fact sum
for $(\ell_v)$ belonging to the set $\Z^{V(\sg)\setminus
  V_N(\sg)}\times\{0\}^{V_N(\sg)}$, i.e.\ we lift the positivity
constraint on vertices in $V_I(\sg)$.

Now, let $r'=\sum_{e\in
  E(\sg)\setminus E_N^{(0)}(\sg)}r_e$, so that $r'\leq r$, and note that 
$$\frac{\frac{r}{\sqrt{n}}}{1+\Big(\frac{r}{\sqrt{n}}\Big)^{p+1}}\leq
\frac{C}{1+\Big(\frac{r}{\sqrt{n}}\Big)^p}\leq
\frac{C}{1+\Big(\frac{r'}{\sqrt{n}}\Big)^p}\, .$$ Therefore, we can
sum out the edge-lengths in the last product of \eqref{eq:43}, and use
elementary inequalities in the second product (as simple as $l+1\leq
2l$ for every integer $l\geq 1$) to get the upper-bound
$$C\,  12^n\sum_{(\ell_v)}\sum_{(r_e)}'\frac{1}{1+
  \Big(\frac{r'}{\sqrt{n}}\Big)^p} \prod_{e\in E_F(\sg)}\frac{1}{\sqrt{r_e}}\cdot\frac{1}{1+\Big(\frac{\ell_{e_+}-\ell_{e_-}}{\sqrt{r_e}}\Big)^2}
\prod_{e\in E^{(1)}_N(\sg)}
\frac{1}{r_e}\cdot\frac{1}{1+\Big(\frac{\ell_{e_-}}{\sqrt{r_e}}\Big)^2}\,
,$$ the symbol $\sum_{(r_e)}'$ meaning that we sum only over all
positive integers $r_e$ indexed by edges $e\in E(\sg)\setminus
E_N^{(0)}(\sg)$. At this point, we write this as an integral
\begin{eqnarray*}
  \lefteqn{C\, 12^n\int\d(\ell_v)\int'\d(r_e)\frac{1}{1+
      \Big(\frac{\sum'_{e}\lfloor r_e\rfloor }{\sqrt{n}}\Big)^p}}\\
  &&\times
  \prod_{e\in E_F(\sg)}\frac{\ind_{\{r_{e}\geq 1\}}}{\sqrt{\lfloor
      r_e\rfloor}}\cdot\frac{1}{1+\Big(\frac{\lfloor\ell_{e_+}\rfloor-\lfloor\ell_{e_-}\rfloor}{\sqrt{\lfloor
        r_e\rfloor}}\Big)^2}
  \prod_{e\in E^{(1)}_N(\sg)}
  \frac{\ind_{\{r_e\geq 1\}}}{\lfloor r_e\rfloor }\cdot\frac{1}{1+\Big(\frac{\lfloor
      \ell_{e_-}\rfloor }{\sqrt{\lfloor r_e\rfloor }}\Big)^2}\,
  ,
\end{eqnarray*}
the first integral being with respect to the measure $\prod_{v\in
  V(\sg)\setminus V_N(\sg)}\d \ell_v\prod_{v\in V_N(\sg)}\delta_0(\d
\ell_v)$, and the second over $\prod_{e\in E(\sg)\setminus
  E_N^{(0)}(\sg)}\d r_e\ind_{\{r_e\geq 0\}}$.  Using the fact that all
the quantities $r_e$ in the integrand are greater than or equal to
$1$, the integral is bounded by
$$C\, 12^n\int\!\d(\ell_v)\!\int'\!\!\d(r_e)\!\frac{1}{1+
      \Big(\frac{r' }{\sqrt{n}}\Big)^p}
  \prod_{e\in E_F(\sg)}\frac{1}{\sqrt{
      r_e}}\cdot\frac{1}{1+\Big(\frac{\ell_{e_+}-\ell_{e_-}}{\sqrt{
        r_e}}\Big)^2}
  \prod_{e\in E^{(1)}_N(\sg)}
  \frac{1}{r_e }\cdot\frac{1}{1+\Big(\frac{
      \ell_{e_-}}{\sqrt{ r_e }}\Big)^2}\, .$$
We perform a linear change of variables, dividing $\ell_v$ by
$n^{1/4}$ and $r_e$ by $n^{1/2}$, which gives after simplification
$$C12^nn^{(\#V_F(\sg)+\#E_F(\sg))/4}\mathscr{I}(\sg)\,
,$$ where 
$$V_F(\sg)=V(\sg)\setminus V_N(\sg)=V_I(\sg)\cup V_O(\sg)\, ,$$
is the set of {\em free vertices}, and
$$\mathscr{I}(\sg)=\int\!\d(\ell_v)\!\int'\!\!\d(r_e)\frac{1}{1+
  (r')^p}\!\!  \prod_{e\in E_F(\sg)}\frac{1}{\sqrt{
    r_e}}\cdot\frac{1}{1+\Big(\frac{\ell_{e_+}-\ell_{e_-}}{\sqrt{
      r_e}}\Big)^2} \prod_{e\in E^{(1)}_N(\sg)} \frac{1}{r_e
}\cdot\frac{1}{1+\Big(\frac{ \ell_{e_-}}{\sqrt{ r_e }}\Big)^2}\, .$$
Provided $\mathscr{I}(\sg)$ is finite for every scheme $\sg$, which we
will check in a moment, we can conclude that \eqref{eq:36} has a
dominant behavior as $n\to\infty$ whenever $\sg$ is such that
$\#V_F(\sg)+\#E_F(\sg)$ is the largest possible. Assume that at least
one of the vertices of $V_N(\sg)$ has degree at least $3$. If we
``free'' this vertex by declaring it in $V_I$ instead, then we have
increased $\#V_F(\sg)$ by $1$, but to make sure that the resulting map
is also a scheme, we should also add a degree-$2$ null vertex at the
center of each edge of $E_N^{(1)}(\sg)$ incident to $v$. Each of these
operations adds an edge in $E_N^{(1)}$, but does not decrease the
cardinality of $E_F(\sg)$. Therefore, the maximal value of
$\#V_F(\sg)+\#E_F(\sg)$ is attained for schemes in which all
null-vertices have degree $2$ (and in particular, $E_N^{(0)}$ is
empty). Furthermore, if there are more than $k+1$ null-vertices with
degree $2$, then at least one of them can be removed without breaking
the condition that $\sg$ is a scheme, because at least two of them are
incident to the same face $f_i$ for some $i\in
\{1,2,\ldots,k\}$. Removing a degree-$2$ null vertex increases $\#E_F$
by $1$ and leaves $\#V_F$ unchanged. Therefore, the optimal schemes
are those having $k$ null vertices, which are all of degree $2$. But
now, the optimal schemes will be obviously those with the largest
number of vertices (or edges), and by definition, these are the
dominant ones. Since a dominant scheme has $2k-2$ free vertices and
$2k-3$ free edges, we obtain that $\#V_F(\sg)+\#E_F(\sg)\leq 4k-5$,
with equality if and only if $\sg$ is dominant. Dominant schemes thus
have a contribution $O(12^nn^{k-5/4})$ to \eqref{eq:36}, while
non-dominant ones have a contribution $O(12^nn^{k-5/4-1/4})$, and by
summing over all such schemes we get the result.

It remains to justify that $\mathscr{I}(\sg)$ is finite for every
scheme $\sg$. To see this, we first integrate with respect to the
variables $(\ell_v)$. We view $E_N^{(1)}$ and the corresponding
incident vertices as a subgraph of $\sg$. Let $\mathfrak{a}$ be a
spanning tree of $\sg/E_N^{(1)}$, that is, every vertex in $V_F(\sg)$
is linked to a vertex of the subgraph $E_N^{(1)}$ by a unique
injective chain, canonically oriented toward $E_N^{(1)}$.

We then bound
$$ \prod_{e\in E_F(\sg)}\frac{1}{\sqrt{
    r_e}}\cdot\frac{1}{1+\Big(\frac{\ell_{e_+}-\ell_{e_-}}{\sqrt{
      r_e}}\Big)^2} \leq \prod_{e\notin E(\mathfrak{a})\cup
  E_N(\sg)}\frac{1}{\sqrt{r_e}} \prod_{e\in
  E(\mathfrak{a})}\frac{1}{\sqrt{r_e}}\cdot\frac{1}{1+\Big(\frac{\ell_{e_+}-\ell_{e_-}}{\sqrt{
      r_e}}\Big)^2}\, .$$ In this form, we can then perform the change
of variables $x_{e_-}=\ell_{e_+}-\ell_{e_-}$ for every $e\in
E(\mathfrak{a})$, and $x_v=\ell_v$ whenever $v\in V_F(\sg)$ is
incident to an edge in $E_N^{(1)}$. This change of variable is
triangular with Jacobian $1$. This allows to
bound $\mathscr{I}(\sg)$ by
\begin{eqnarray*}
  \int'\d (r_e) \frac{1}{1+(r')^p}
  \prod_{e\notin E(\mathfrak{a})\cup E_N(\sg)}\frac{1}{\sqrt{r_e}}
  \prod_{e\in E(\mathfrak{a})} \int \frac{\d
    x_{e_-}/\sqrt{r_e}}{1+\Big(\frac{x_{e_-}}{\sqrt{
        r_e}}\Big)^2}\prod_{v\in V_F(\sg)}\int \prod_{\substack{e\in
      E_N^{(1)}(\sg)\\e_-=v}}\frac{1}{r_e}\cdot\frac{\d
    x_v}{1+(x_v/\sqrt{r_e})^2}\\
  = C\int'\d (r_e) \frac{1}{1+(r')^p}
  \prod_{e\notin E(\mathfrak{a})\cup E_N(\sg)}\frac{1}{\sqrt{r_e}}
  \prod_{v\in V_F(\sg)}\int \prod_{\substack{e\in
      E_N^{(1)}(\sg)\\e_-=v}}\frac{1}{r_e}\cdot\frac{\d
    x_v}{1+(x_v/\sqrt{r_e})^2}
\end{eqnarray*}
The last integral is of the form (for positive $r_1,\ldots,r_l$)
\begin{eqnarray*}
\int \d x\prod_{i=1}^l
\frac{1}{r_i}\cdot\frac{1}{1+(x/\sqrt{r_i})^2}&\leq& C\int \d
x\prod_{i=1}^l \frac{1}{r_i}\wedge \frac{1}{x^2}\\
&\leq & C\frac{\sqrt{\min(r_1,\ldots,r_l)}}{\prod_{i=1}^lr_i}\, ,
\end{eqnarray*}
for some constant $C$ depending on $l$, the second step being easy to
prove by first assuming that $r_1<r_2<\ldots<r_l$ and decomposing the
integral along the intervals
$[0,\sqrt{r_1}],[\sqrt{r_1},\sqrt{r_2}],\ldots,[\sqrt{r_l},\infty[$. 
Note that the last upper-bound is integrable on $[0,1]^l$ with respect to $\d
r_1\ldots \d r_l$, since 
$$\int_{[0,1]^l}\prod\frac{\d
  r_i}{r_i}\sqrt{\min(r_1,\ldots,r_i)}=l\int_0^1\frac{\d
  r_1}{\sqrt{r_1}}\int_{[r_1,1]^{l-1}}\prod_{i=2}^l\frac{\d
  r_i}{r_i}=l\int_0^1\frac{\d r_1}{\sqrt{r_1}}\log(1/r_1)^{l-1}\, ,$$
which is finite. By choosing $p>\#E(\sg)$ (since there are only a
finite number of elements of $\bS^{(k+1)}$, we can choose a single $p$
valid for every scheme with $k+1$ faces), we finally obtain that
$\mathscr{I}(\sg)$ is bounded from above by an integral in $(r_e)$,
whose integrand is both integrable in a neighborhood of $0$ and of
infinity. Therefore $\mathscr{I}(\sg)$ is finite, as wanted.  \cq

\medskip

We now prove Proposition \ref{sec:limit-theorems-1}. By Lemma
\ref{sec:limit-theorems-2}, it suffices to consider the restriction of
$\LM_n^{(k+1)}$ to the labeled maps whose associated scheme is in
$\bS^{(k+1)}_d$. We start by considering the image of the measure
$\LM_n^{(k+1)}$ under the mapping $(\sg,(W_e)_{e\in \vec{E}\sg})\mapsto
(\sg,(\ell_v)_{v\in V(\sg)},(r_e)_{e\in E(\sg)})$, which we still
write $\LM_n^{(k+1)}$ by abuse of notation. Let $f$ be a continuous
function on $\bS_d^{(k+1)}\times \R^{V(\sg)}\times \R^{E(\sg)}$, we
can assume that $f(\sg',(\ell_v),(r_e))$ is non-zero only when $\sg'$
equals some particular dominant scheme $\sg\in \bS_d^{(k+1)}$, and we
drop the mention of the first component of $f$ in the sequel.

Then we have, by similar arguments as in the derivation of
\eqref{eq:45},
$$\psi^{(k+1)}_n{}_*\LM^{(k+1)}_n(f)=2\cdot 3^{n+k}\sum_{(\ell_v)}\sum_{(r_e)}r\binom{2n}{n-r}\prod_{e\in
  \check{E}(\sg)}q^{(e)}_{r_e}(\ell_{e_-},\ell_{e_+})
f\bigg(\Big(\Big(\frac{9}{8n}\Big)^{1/4}\ell_v\Big),\Big(\frac{r_e}{\sqrt{2n}}\Big)\bigg)\,
,$$ where $q^{(e)}_r(a,b)=3^{-r}\mathcal{W}^{(e)}(a,b;r)$ and
$r=\sum_e r_e$. Here, the factor $3^k$ comes from the fact that there
are $k$ edges in $E_T(\sg)$, because $\sg$ is dominant, and each such
edge corresponds to a Motzkin walk with a final ``virtual'' step in
the end, participating an extra factor $3$, as explained in the proof
of Proposition \ref{sec:limit-theorems-2}. We write this as an
integral
$$2\cdot 12^n3^k\int \lambda_\sg(\d (\ell_v))\int\d
  (r_e)\frac{[r]}{4^n}\binom{2n}{n-[r]}
\prod_{e\in \check{E}(\sg)}q^{(e)}_{\lfloor r_e\rfloor}(\lfloor
\ell_{e_-}\rfloor,\lfloor \ell_{e_+}\rfloor)
f\bigg(\Big(\Big(\frac{9}{8n}\Big)^{1/4}\lfloor \ell_v\rfloor\Big),
\Big(\frac{\lfloor r_e\rfloor}{\sqrt{2n}}\Big)\bigg)\, ,$$
where
$[r]=\sum_{e\in E(\sg)}\lfloor r_e\rfloor$, and where we omitted to
write indicators $\ind_{\{r_e\geq 1\}}$ and $\ind_{\{\ell_v\geq 1\}}$
whenever $v\in V_I(\sg)$ to lighten the expression. We perform a
linear change of variables, dividing $r_e$ by $\sqrt{2n}$ and $\ell_v$
by $(8n/9)^{1/4}$, yielding
\begin{eqnarray*}
\lefteqn{2\cdot 12^n3^k(2n)^{\#E(\sg)/2}(8n/9)^{\#V_F(\sg)/4}\int \lambda_\sg(\d (\ell_v))\int\d
  (r_e)\frac{[r\sqrt{2n}]}{4^n}\binom{2n}{n-[r\sqrt{2n}]}}\\
&&
\prod_{e\in \check{E}(\sg)}q^{(e)}_{\lfloor r_e\sqrt{2n}\rfloor}\Big(\Big\lfloor
\ell_{e_-}\Big(\frac{8n}{9}\Big)^{1/4}\Big\rfloor,\Big\lfloor
\ell_{e_+}\Big(\frac{8n}{9}\Big)^{1/4}\Big\rfloor\Big)
f\bigg(\Big(\Big(\frac{9}{8n}\Big)^{1/4}\Big\lfloor \ell_v\Big(\frac{8n}{9}\Big)^{1/4}\Big\rfloor\Big),
\Big(\frac{\lfloor r_e\sqrt{2n}\rfloor}{\sqrt{2n}}\Big)\bigg)\, . 
\end{eqnarray*}
Now, by the local limit theorem, it holds that for every fixed $r>0$
and $\ell,\ell'\in \R$, we have 
$$\Big(\frac{8n}{9}\Big)^{1/4}q_{\lfloor r_e\sqrt{2n}\rfloor
}\Big(\Big\lfloor
\Big(\frac{8n}{9}\Big)^{1/4}\ell\Big\rfloor,\Big\lfloor
\Big(\frac{8n}{9}\Big)^{1/4}\ell'\Big\rfloor\Big)\build\longrightarrow_{n\to\infty}^{}p_{r_e}(\ell,\ell')\,
, $$ and similarly with $q^+$ and $p^+$ instead of $q$ and $p$ (recall
\eqref{eq:40}), whenever $\ell,\ell'>0$, while for $\ell>0$, using
\eqref{eq:41} and the local limit theorem, 
$$\sqrt{2n}q^+_{\lfloor r_e\sqrt{2n}\rfloor
}\Big(\Big\lfloor \Big(\frac{8n}{9}\Big)^{1/4}\ell\Big\rfloor,0\Big)
\build\longrightarrow_{n\to\infty}^{}\ov{p}_{r_e}(\ell,0)\, , .$$ A
final use of the local limit theorem (or Stirling's formula), as in
\eqref{eq:42} shows that for positive $r_e,e\in E(\sg)$, with same
notation as above,
$$\frac{[r\sqrt{2n}]}{4^n}\binom{2n}{n-[r\sqrt{2n}]}\build\longrightarrow_{n\to\infty}^{}
\ov{p}_1(2r,0)\, .$$ Finally, since dominant schemes have $4k-3$ edges
(of which $2k-3$ are free) and $3k-2$ vertices ($2k-2$ free) we see
after simplifications that the integral is equivalent to
$$3^k2^{k}(2/9)^{1/4}12^n
n^{k-5/4}\int\lambda_\sg(\d(\ell_v))\int\d(r_e)\ov{p}_1(2r,0)
\Big(\prod_{e\in
  \check{E}(\sg)}p^{(e)}_{r_e}(\ell_{e_-},\ell_{e_+})\Big)f((\ell_v),(r_e))\,
,$$ at least if we can justify to take the limit inside the
integral. This is done by dominated convergence: We bound $f$ by its
supremum norm, and apply the very same argument as in the proof of
Lemma \ref{sec:limit-theorems-2}, using \eqref{eq:39} and
\eqref{eq:42} to bound the integrand by the integrand of the quantity
$\mathscr{I}(\sg)$. Details are left to the reader.

In summary, by definition of $\CLM_1^{(k+1)}$, we have proven that
$$\Big(\frac{9}{2}\Big)^{1/4}\frac{
  \psi^{(k+1)}_n{}_*\LM_n^{(k+1)}}{6^k\cdot12^nn^{k-5/4}
}(F)\build\longrightarrow_{n\to\infty}^{}\Upsilon_{k+1}\CLM_1^{(k+1)}(F)\,
,$$ when $F$ is a continuous and bounded function of the form
$F(\sg,(W_e))=f(\sg,(\ell_v),(r_e))$. In order to prove the full
result, we note that Proposition \ref{sec:decomp-label-maps-1} entails
that under the measure $\LM^{(k+1)}_n$, conditionally given
$(\sg,(\ell_v),(r_e))$, the paths $(M_e,e\in \check{E}(\sg))$ are
independent discrete walks, respectively from $\ell_{e_-}$ to
$\ell_{e_+}$ and with duration $r_e$, which are conditioned to be
positive (except possibly at their final point) when $e\in
E_N(\sg)\cup E_I(\sg)$. This is in fact not perfectly right since
paths $M_e$ with $e\in E_T(\sg)$ come with one extra negative step,
but this is of no incidence to what follows. Since we have shown that
$r_e$ scales like $\sqrt{2n}$ and $\ell_v$ like $(8n/9)^{1/4}$, it is
quite standard (see \cite[Lemmas 10 and 14]{bettinelli10} for a recent
and thorough exposition) that the paths $M_e$, after applying the
operation $\psi^{(k+1)}_n$, converge in distribution to independent
Brownian bridges (conditioned to be positive if $e\in E_I(\sg)$, or
first-passage bridges when $e\in E_N(\sg)$), in accordance with the
definition of $\CLM^{(k+1)}_1$.

Lastly, conditionally on the paths $(M_e,e\in E(\sg))$ under the
discrete measure $\LM^{(k+1)}_n$, the snakes $(W_e,e\in \vec{E}(\sg))$
are associated with independent labeled forests $(F_e,L_e)$ with
respectively $r_e$ trees, conditioned on having a total number of
oriented edges equal to $2n$. Moreover, the labels are uniform among
all admissible labelings, with the constraint that the labels of the
root vertices in the forest associated with the edge $e\in
\vec{E}(\sg)$ are given by the path $M_e$. By subtracting these labels
we obtain forests with root labels $0$, and uniform labelings among
the admissible ones, so that we are only concerned in the convergence
of the discrete snakes associated with these forests to independent
Brownian snakes started respectively from the constant trajectory
equal to $0$ and duration $r_e$, and conditioned on $\sigma=\sum_{e\in
  \vec{E}(\sg)}\sigma_e=1$, under the law $\CLM^{(k+1)}_1$
conditionally on $(\sg,(\ell_v),(r_e))$.

But the forests $F_e$ can be obtained from one single forest with $2r$
trees and total number of oriented edges (comprising floor edges)
equal to $2n$, by cutting the floor into segments with respectively
$r_e$ trees, for $e\in \vec{E}(\sg)$. The random snake associated with
this forest converges, once rescaled according to the operation
$\psi^{(k+1)}_n$, to a Brownian snake with total duration $1$, starting
from the constant trajectory equal to $0$ and with duration $2r$, by
\cite[Proposition 15]{bettinelli10}. It is then easy to see that the
snakes obtained by ``cutting this snake into bits of lengths $r_e$'',
as explained in the definition of $\CLM_1^{(k+1)}$ in Section
\ref{sec:measure-clmk+1}, is indeed the limiting analog of the
discrete cutting just mentioned. A completely formal proof
requires to write the discrete analogs of \eqref{eq:46} and
\eqref{eq:47} and verify that they pass to the limit, which is a
little cumbersome and omitted.  \cq

\subsection{The case of planted
  schemes}\label{sec:case-planted-schemes}

Recall that in a planted scheme $\sg\in \dot{\bS}^{(k+1)}$ with $k+1$
faces, there is a unique vertex of degree $1$, the others being of
degree $3$. This vertex is an element of $V_O(\sg)$, and the edge
incident to this vertex is an element of $E_O(\sg)$. The elements of
$\dot{\bS}^{(k+1)}_d$ (the dominant planted schemes) have $3k$
vertices, of which $k$ are null vertices (all of degree $2$), and
$4k-1$ edges, of which $2k$ are elements of $E_N(\sg)$.

On the space 
$$\dot{\bCLM}^{(k+1)}=\Big\{(\sg,(W_e)_{e\in \vec{E}(\sg)}):\sg\in
\dot{\bS}^{(k+1)},(W_e)_{e\in \vec{E}(\sg)}\in
\CC(\CC(\R))^{\vec{E}(\sg)}\Big\}\, ,$$ the ``planted continuum
measure'' $\dot{\CLM}^{(k+1)}$ and its conditioned counterpart
$\dot{\CLM}^{(k+1)}_1$ are defined as in formulas \eqref{eq:19} and
\eqref{eq:38} for $\CLM^{(k+1)}$ and $\CLM_1^{(k+1)}$, only replacing
the counting measure $\S_d^{(k+1)}$ by the counting measure
$\dot{\S}_d^{(k+1)}$ on planted scheme, and changing
$\Upsilon^{(k+1)}$ by the proper normalization constant
$\dot{\Upsilon}^{(k+1)}$. These are measures on the space
$\bCLM^{(k+1)}$ that constitutes in the pairs $(\sg,(W_e)_{e\in
  \vec{E}(\sg)})$, where $\sg$ is a planted scheme, and $(W_e)_{e\in
  \vec{E}(\sg)}$ is an admissible family of discrete snakes.

  The analogous statements to Propositions \ref{sec:cont-meas-label-1}
  and \ref{sec:limit-theorems-2} goes as follows. We define the
  scaling operations $\dot{\Psi}_c$ and $\dot{\psi}^{(k+1)}_n$ on
  $\dot{\bCLM}^{(k+1)}$ by the same formulas as $\Psi^{(k+1)}_c$ and $\psi^{(k+1)}_n$
  in sections \ref{sec:measure-clmk+1} and
  \ref{sec:limit-theorems}. In the following statement, we use the
  decomposition of Section \ref{sec:keeping-track-root} rather than
  Proposition  \ref{sec:decomp-label-maps-1}, and we view $\LM_n$ as a
  measure on $\dot{\bCLM}^{(k+1)}$. 

\begin{prp}\label{sec:case-planted-schemes-1}
It holds that 
$$\dot{\CLM}^{(k+1)}=\int_0^\infty \d \sigma\, \sigma^{k-5/4}
\dot{\CLM}_\sigma^{(k+1)}\, .$$ Moreover, we have the following weak convergence
of finite measures on $\dot{\bCLM}^{(k+1)}$: 
$$\Big(\frac{9}{2}\Big)^{1/4}\frac{
\dot{\psi}^{(k+1)}_n{}_*\LM_n^{(k+1)}}{6^k\cdot12^nn^{k-5/4}
}\build\longrightarrow_{n\to\infty}^{}\dot{\Upsilon}^{(k+1)}\dot{\CLM}_1^{(k+1)}\, .$$
\end{prp}

The proof follows exactly the same lines as Proposition
\ref{sec:limit-theorems-1}, so we leave it to the reader. One has to
be a little careful about the small variations in the construction of
Section \ref{sec:keeping-track-root}, compared to Proposition
\ref{sec:decomp-label-maps-1}, when we are dealing with the
distinguished edge and its adjacent edges, but these variations
disappear in the scaling limit. Also, note that the reason why the
measure $\d \sigma\,  \sigma^{k-5/4}$ appears in the disintegration of
$\dot{\CLM}^{(k+1)}$, instead of $\d \sigma\, \sigma^{k-9/4}$, is that
$\dot{\CLM}^{(k+1)}$ carries intrinsically the location of a
distinguished point (corresponding to the only vertex of the scheme
that has degree $1$). This marked  point should be seen as the continuum
counterpart of the root edge in discrete maps, so it is natural that
if the continuum object has a total ``mass'' $\sigma$, then there
marking introduces a further factor $\sigma$. 

In fact, we could recover Propositions \ref{sec:cont-meas-label-1} and
\ref{sec:limit-theorems-1} from Proposition
\ref{sec:case-planted-schemes-1}, by considering the natural operation
from $\dot{\bS}^{(k+1)}$ to $\bS^{(k+1)}$ that erases the
distinguished edge incident to the degree-$1$ vertex. We leave to the
reader to properly formulate and prove such a statement, which we are
not going to need in the sequel. The reason why we are dealing with
planted schemes (which are ``richer'' objects) only now is that we
find them a little harder to understand and manipulate than
schemes. Using non-planted schemes will also simplify the proofs in
Section \ref{sec:eps-geodesic-stars}.

\section{Proof of the key lemmas}\label{sec:proof-key-lemmas}

We finally use the results of Sections \ref{sec:coding-labeled-maps}
and \ref{sec:scal-limits-label} to prove Lemmas
\ref{sec:covering-3-star} and \ref{sec:covering-3-star-1}.  Recall
from Proposition \ref{sec:geod-stars-discr} how the probabilities of
events $\mathcal{A}_1(\eps,\beta)$ and $\mathcal{A}_2(\eps)$ are
dominated by the limsup of probabilities of related events
$\mathcal{A}^{(n)}_1(\eps,\beta)$ and $\mathcal{A}^{(n)}_2(\eps)$ for
quadrangulations. The latter events were defined on a probability
space supporting the random variables $Q_n$ and marked vertices
$v_0,v_1,\ldots,v_k$, in the sequel we will see them as sets of marked
quadrangulations $(\bq,\bv)$ with $\bq\in \bQ_n$ and $\bv\in
V(\bq)^{k+1}$, by abuse of notation.

It turns out that the latter events have a tractable
translation in terms of the labeled map $(\bm,\bl)$ associated with
the random quadrangulation $Q_n$ by the multi-pointed bijection
$\Phi^{(k+1)}$ of Section \ref{sec:multi-point-scha}.

\subsection{Relation to labeled maps}\label{sec:relat-label-maps}

Let us first consider $\mathcal{A}^{(n)}_1(\eps,\beta)$. 
In the sequel, if $X\in \CC(\R)$ we let $\un{X}=\inf_{0\leq t\leq
  \zeta(X)}X(t)$, and if $W\in \CC(\CC(\R))$ we let
$\un{W}=\inf_{0\leq s\leq \sigma(W)}\un{W(s)}$.  We let
$\mathcal{B}_1(\eps,\beta)$ be the event on $\dot{\bCLM}^{(3)}$ that
\begin{itemize}
\item for every $e\in \vec{E}(\sg)$ incident to $f_0$, it holds that
  $\un{W}_e\geq -2\eps$,
\item if $e\in E(\sg)$ is incident to $f_1$ and $f_2$, then
  $\un{M}_e\geq 0$,
\item there exists $e\in E_I(\sg)$ such that, if we orient $e$ in such
  a way that it is incident to $f_1$ or $f_2$, then $\un{W}_e\leq
  -\eps^{1-\beta}$.
\end{itemize}

Let $(\bq,\bv)\in \mathcal{A}^{(n)}_1(\eps,\beta)$, $r=\lfloor
\eps(8n/9)^{1/4}\rfloor -1$, and $r'\in \{r+1,r+2,\ldots,2r\}$. We
note that $\mathcal{A}^{(n)}_1(\eps,\beta)$ contains the event
$G(r,2)$, so by Lemma \ref{sec:geodesic-r-stars-1}, the labeled map
$(\bm,\bl)=\Phi^{(3)}(\bq,\bv,\btau^{(r')})$ is an element of
$\bLM^{(3)}$. We let $(\sg,(W_e,e\in \vec{E}(\sg)))$ be the element of
$\dot{\bCLM}^{(k+1)}$ associated with $(\bm,\bl)$ as in Section
\ref{sec:keeping-track-root} (note that we choose to take planted
schemes here).  Recall the definition of the scaling functions
$\psi^{(k+1)}_n$ and $\dot{\psi}^{(k+1)}_n$ from Sections
\ref{sec:limit-theorems} and \ref{sec:case-planted-schemes}.

\begin{lmm}\label{sec:proof-key-lemmas-1} 
  With the above notation, for every $\eps>0$ small enough,
  $(\sg,(W_e,e\in \vec{E}(\sg)))$ belongs to the event
  $\mathcal{B}^{(n)}_1(\eps,\beta)=(\dot{\psi}^{(3)}_n)^{-1}(\mathcal{B}_1(\eps,\beta))$
\end{lmm}

\proof The first point is a consequence of the fact \eqref{eq:13} that the minimal
label of a corner incident to $f_0$ is equal to 
$\tau^{(r')}_0+1$, which is $-r'+1\geq -2\eps(8n/9)^{1/4}$ by our choice.

For the second point, since
$\tau^{(r')}_i=-d_\bq(v_i,v_0)-\tau^{(r')}_0$, we deduce that for every $v\in
B_{d_\bq}(v_0,\eps(8n/9)^{1/4})$, 
$$d_\bq(v,v_0)+\tau^{(r')}_0 < 0 <
d_\bq(v,v_i)+\tau^{(r')}_i\, , $$ since
$d_\bq(v_0,v_i)-d_\bq(v,v_i)\leq d_\bq(v,v_0)<-\tau^{(r')}_0$ by
the triangle inequality. Consequently, we obtain that
$B_{d_\bq}(v_0,\eps(8n/9)^{1/4})\setminus\{v_0\}$ is included in the
set of vertices of $\bm$ that are incident to the face $f_0$, and to
no other face. Under the event $\mathcal{A}^{(n)}_1(\eps,\beta)$, we
thus obtain that any geodesic from $v_1$ to $v_2$ has to visit
vertices incident to $f_0$ exclusively. 

Therefore, such a geodesic has to visit a vertex incident to $f_0$ and
to either $f_1$ or $f_2$. By definition of $\bLM^{(3)}$, such a vertex
$v$ has non-negative label $l\geq 0$, and the length of the geodesic
chain has to be at least (we let $\min_{f_1}\bl=\min_{w\in
  V(f_1)}\bl(w)$)
\begin{equation}
  \label{eq:22} d_\bq(v_1,v_2)\geq
  \Big(\bl(v)-\min_{f_1}\bl+1\Big)+\Big(\bl(v)-\min_{f_2}\bl+1\Big)\geq
  2-\min_{f_1}\bl-\min_{f_2}\bl\, .
\end{equation}
Now assume that $v$ a vertex incident to both $f_1$ and $f_2$, and let
$l$ be its label. Choose two corners $e,e'$ incident to $v$ that
belong respectively to $f_1$ and $f_2$. By drawing the leftmost
geodesic chains from $e,e'$ to $v_1,v_2$ visiting the consecutive
successors of $e,e'$ as in Section \ref{sec:reverse-construction}, and
concatenating these chains, we obtain a chain visiting only vertices
that are incident to $f_1$ or $f_2$. This chain cannot be geodesic
between $v_1$ and $v_2$ because it does not visit vertices incident to
$f_0$ exclusively. But its length is given by
$$\Big(l-\min_{f_1}\bl+1\Big) + \Big(l-\min_{f_2}\bl+1\Big)\geq d_\bq(v_1,v_2)\, .$$
Comparing with \eqref{eq:22}, we see that necessarily, $l\geq 0$.
Hence all labels of vertices incident to $f_1$ and $f_2$ are
non-negative, implying the second point in the
definition of $\mathcal{B}^{(n)}_1(\eps,\beta)$.

For the third point, we argue by contradiction, and assume that for
every edge $e\in \vec{E}(f_1)$ of $\sg$ with $\ov{e}\in \vec{E}(f_0)$,
(resp.\ $e\in \vec{E}(f_2)$ with $\ov{e}\in \vec{E}(f_0)$), the labels
in the forest $F_e$ that is branched on $e$ are all greater than or
equal to $-\eps^{1-\beta}(8n/9)^{1/4}$.  Note that every dominant
pre-scheme with $3$ faces is such that there is a unique edge incident
both to $f_0$ and $f_1$ (resp.\ $f_0$ and $f_2$). Consider a geodesic
chain $\gamma$ from $v_2$ to $v_1$. Since this chain has to visit
vertices exclusively contained in $f_0$, we can consider the last such
vertex in the chain. Then the following vertex $v$ is necessarily
incident to $f_1$: Otherwise, it would be incident to $f_2$, and the
concatenation of a leftmost geodesic from $v$ to $v_2$ with the
remaining part of $\gamma$ between $v$ and $v_1$ would be a shorter
chain from $v_2$ to $v_1$ not visiting strictly $f_0$.

Necessarily, there is a vertex $v'$ with label $0$ incident to $f_1$
and $f_0$. Let $e,e'$ be corners of $\bm$ incident to $f_1$ and
respectively incident to $v$ and $v'$. By hypothesis, when visiting
one of the intervals $[e,e']$ or $[e',e]$, we encounter only corners
with labels greater than or equal to
$-\eps^{1-\beta}(8n/9)^{1/4}$. Consider the left-most geodesics
$\tilde{\gamma},\tilde{\gamma}'$ from $e$ and $e'$ to $v_1$, visiting
their consecutive successors. By replacing the portion of $\gamma$
from $v$ to $v_1$ by $\tilde{\gamma}$, we still get a geodesic from
$v_2$ to $v_1$. On the other hand, $\tilde{\gamma}'$ is a portion of a
geodesic from $v_0$ to $v_1$, in which an initial segment of length at
most $2r\leq 2\eps(8n/9)^{1/4}$ has been removed. So we have found a
geodesic from $v_0$ to $v_1$ such that every vertex $v''$ on this
geodesic that lies outside of
$B_{d_\bq}(v_0,(2\eps+\eps^{1-\beta})(8n/9)^{1/4})$ is such that
$(v_2,v'',v_1)$ are aligned.  The same holds with the roles of $v_1$
and $v_2$ interchanged, and for $\eps$ small enough we have $2\eps\leq
\eps^{1-\beta}$ therefore, $\mathcal{A}^{(n)}_1(\eps,\beta)$ does not
hold.  \cq

\medskip

As a corollary, we obtain that for every $\beta\in (0,1)$ and $\eps>0$
small enough,
\begin{equation}\label{eq:48}
  \limsup_{n\to\infty}P(\mathcal{A}^{(n)}_1(\eps,\beta))\leq
  \frac{C}{\eps}\dot{\CLM}_1^{(3)}(\mathcal{B}_1(\eps,\beta))\, 
.
\end{equation}
To justify this, we write, since $Q_n$ is uniform in $\bQ_n$ and
$(v_0,v_1,\ldots,v_k)$ are independent uniformly chosen points in
$V(Q_n)$, if we let $r=\lfloor \eps(8n/9)^{1/4}\rfloor-1$,
\begin{eqnarray*}
  P(\mathcal{A}^{(n)}_1(\eps,\beta))&=&
  \frac{1}{n^{k+1}\#\bQ_n}\sum_{\bq\in \bQ_n, 
    \bv\in V(\bq)^{k+1}}\ind_{\{(\bq,\bv)\in
    \mathcal{A}^{(n)}_1(\eps,\beta)\}}\\
  &=& \frac{1}{rn^{k+1}\#\bQ_n}\sum_{\bq\in \bQ_n, 
    \bv\in
    V(\bq)^{k+1}}\sum_{r'=r+1}^{2r}\ind_{\{(\bq,\bv)\in\mathcal{A}^{(n)}_1(\eps,\beta)\}}\, .
\end{eqnarray*}
Since $\Phi^{(k+1)}$ is two-to-one,
and by Lemma \ref{sec:geodesic-r-stars-2}, this is bounded from above
by 
$$
\frac{2}{rn^{k+1}\#\bQ_n}\sum_{(\bm,\bl)\in
  \bLM^{(3)}_n}\ind_{\{(\Phi^{(k+1)})^{-1}(\bm,\bl)
\in\mathcal{A}^{(n)}_1(\eps,\beta)\}}\\
\leq 
  \frac{C}{\eps 12^n n^{k-5/4}}\dot{\psi}^{(3)}_n{}_*\LM_n(\mathcal{B}_1(\eps,\beta))\, ,
$$
where we finally used Lemma \ref{sec:proof-key-lemmas-1} and the
well-known fact that
$$\#\bQ_n=\frac{2}{n+2}\cdot\frac{3^n}{n+1}\binom{2n}{n}\, .$$
Finally, Proposition \ref{sec:case-planted-schemes-1} entails
\eqref{eq:48}, since $\mathcal{B}_1(\eps,\beta)$ is a closed set.  We
will estimate the upper-bound in \eqref{eq:48} in Section
\ref{sec:fast-confl-geod}.

Let us now bound the probability of $\mathcal{A}^{(n)}_2(\eps)$. We
say that a dominant scheme $\sg\in \bS_d^{(k+1)}$ is {\em predominant}
if $f_0$ has minimal degree in $\sg$. For $k=3$, we see that the
top-left pre-scheme of Figure \ref{fig:shapes} is the only one that
gives rise to predominant schemes with $4$ faces by adding a
null-vertex in the middle of the three edges incident to $f_0$: All
the dominant schemes constructed from the four others pre-schemes will
have at least $7$ oriented edges incident to $f_0$. We let
$\mathfrak{P}$ be the set of predominant schemes.

%
We consider the event $\mathcal{A}_3\subset\bLM^{(4)}$
that
\begin{itemize}
\item the scheme $\sg$ associated with the map belongs to
  $\mathfrak{P}$,
\item for every $v\in V(f_1\cap f_2)$, it holds that $\bl(v)\geq 0$,
\item for every $v\in V(f_3\cap f_1)$ and $v'\in V(f_3\cap f_2)$, and
  $e,e'$ two corners of $f_3$ incident to $v,v'$, if we take the
  convention that $[e,e']$ is the set of corners between $e$ and $e'$
  in facial order around $f_3$, that passes through the corners of the
  unique edge of $\sg$ incident to $f_3$ and $f_0$, then it holds that
\begin{equation}
\label{eq:50}
\bl(v)+\bl(v')+2\geq \min_{e''\in [e,e']}\bl(e'')\, .
\end{equation}
\end{itemize}
For instance, in Figure \ref{fig:scheme0} below, the interval $[e,e']$
is represented as the dotted contour inside the face $f_3$. 

Let $(\bq,\bv)\in \mathcal{A}^{(n)}_2(\eps)$, $r=\lfloor \eps
(8n/9)^{1/4}\rfloor-1$ and $r'\in \{r+1,r+2,\ldots,2r\}$. By Lemma
\ref{sec:geodesic-r-stars-1}, the labeled map
$(\bm,\bl)=\Phi^{(4)}(\bq,\bv,\btau^{(r')})$ is an element of
$\bLM^{(4)}$, and we let $(\sg,(W_e,e\in \vec{E}(\sg)))$ be the
element of $\bCLM^{(4)}$ associated with it via Proposition
\ref{sec:decomp-label-maps-1}.
\begin{lmm}\label{sec:relat-label-maps-1}
  Under these hypotheses, it always hold that $\min_{e\in
    \vec{E}(f_0)}\un{W}_e\geq -2\eps(8n/9)^{1/4}$, and moreover, if $\sg\in
  \mathfrak{P}$ then $(\bm,\bl)\in \mathcal{A}_3$.
\end{lmm}

\proof The lower-bound on $\min_{e\in \vec{E}(f_0)}\un{W}_e$ is just
\eqref{eq:13}. Next, let us assume that $\sg\in \mathfrak{P}$, and let
us check that $\mathcal{A}_3$ is satisfied. The second point is
derived in exactly the same way as we checked the second point of
$\mathcal{B}_1(\eps,\beta)$ in the derivation of Lemma
\ref{sec:proof-key-lemmas-1}.  Indeed, on $\mathcal{A}_2(\eps)$, for
$n$ large enough, all geodesic paths from $v_1$ to $v_2$ have to pass
through $f_0$, and so they have length at least $-\min_{v\in
  f_1}\bl(v)-\min_{v\in V(f_2)}\bl(v)+2$.
So if $\bl(v)<0$ for some $v\in V(f_1\cap f_2)$, then by drawing the
successive arcs starting from two corners of $f_1$ and $f_2$ incident
to $v$, until $v_1$ and $v_2$ are reached, we would construct a path
with length at most $-\min_{v\in V(f_1)}\bl(v)-\min_{v\in
  V(f_2)}\bl(v)$ between $v_1$ and $v_2$, a contradiction. 

For the third point, note that if we draw the two leftmost geodesic
chains from $e,e'$ to $v_3$ inside the face $f_3$, then these two
chains coalesce at a distance from $v_3$ which is precisely
$\min_{e''\in [e,e']}\bl(e'')+2$. Therefore, it is possible to build a
geodesic chain from $v_1$ to $v$, with length $\bl(v)-\min_{u\in
  V(f_1)}\bl(u)+1$, and to concatenate it with a chain of length
$\bl(v)+\bl(v')-2 \min_{e''\in [e,e']}\bl(e'')+2$ from $v$ to $v'$,
and then with a geodesic chain from $v'$ to $v_2$, with length
$\bl(v')-\min_{u\in V(f_2)}\bl(u)+1$. Since the resulting path cannot
be shorter than a geodesic from $v_1$ to $v_2$, we obtain the third
required condition.  \cq

\medskip

We now translate the event $\mathcal{A}_3$ in terms of the
encoding processes $(\sg,(W_e)_{e\in \vec{E}(\sg)})$.  If $\sg$ is a
predominant scheme with $4$ faces, then the associated pre-scheme is
the first one of Figure \ref{fig:shapes}. Therefore, there is a single
edge $e_{ij}$ incident to $f_i$ and $f_j$ for every $i<j$ in
$\{1,2,3\}$.  We let $\mathcal{B}_2(\eps)\subset \bCLM^{(4)}$ be the
event that
\begin{enumerate}
\item the scheme $\sg$ belongs to $\mathfrak{P}$, 
\item $\min_{e\in \vec{E}(f_0)}\un{W}_e\geq -2\eps$, 
\item $\un{M}_{e_{12}}\geq 0$, 
\item for every $t\in [0,r_{e_{13}}]$ and $t'\in [0,r_{e_{23}}]$, it
  holds that
\begin{equation}\label{eq:44}
  \min_{e\in \vec{E}_3(\sg):\ov{e}\in \vec{E}_0(\sg)}\un{W}_e\wedge \inf\{\un{W}_{e_{13}}^{(s)},0\leq
  s\leq t\}\wedge \inf\{\un{W}_{e_{23}}^{(s)},0\leq s\leq t'\}\leq 
  M_{e_{13}}(t)+M_{e_{23}}(t')
\end{equation}
\end{enumerate}

\begin{figure}[htb!]
\begin{center}
\includegraphics[scale=1.2]{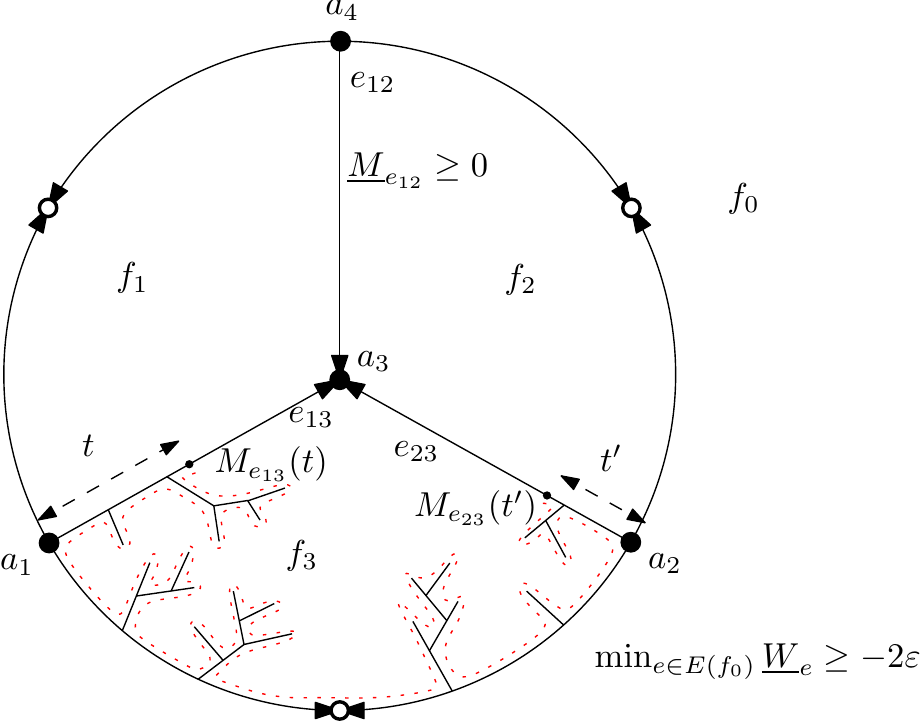}
\end{center}
\caption{Illustration for the event $\mathcal{B}_2(\eps)$. Only a
  small portion of the snakes branching on the scheme have been
  represented here. The scheme underlying this element of
  $\bCLM^{(4)}$ is a predominant scheme with $4$ faces, and all the
  others are obtained by obvious symmetries. All labels along the edge
  $e_{12}$ are non-negative, all labels in $f_0$ are greater than or
  equal to $-2\eps$, and for every $t\in [0,r_{e_{13}}],t'\in
  [0,r_{e_{23}}]$, the minimal label along the dotted contour is at
  most $M_{e_{13}}(t)+M_{e_{23}}(t')$. The blank vertices indicate
  elements of $V_N(\sg)$, and the labels $a_1,a_2,a_3,a_4$ are the
  integration variables appearing in the proof of Lemma
  \ref{sec:eps-geodesic-stars-1} below. }
\label{fig:scheme0}
\end{figure}

Here and in the remainder of the paper, we use a slightly unusual
convention, that $W^{(s)}_{e_{13}},W^{(s)}_{e_{23}}$ are the snake
excursions (the tree components in the continuum forest) branching on
$e_{13}$ and $e_{23}$ that lie inside $f_3$, where the orientation of
$e_{13},e_{23}$ points away from $f_0$. For instance, in Figure
\ref{fig:scheme0}, these are the tree components that lie to the {\em
  right} of $e_{13}$ instead of the left, as our usual conventions
would require: In this case, we should really read
$W^{(r_{e_{13}}-s)}_{\ov{e}_{13}}$ rather than $W^{(s)}_{e_{13}}$.

Although it is a little tedious, it is really a matter of definitions
to check that if $(\bm,\bl)\in \mathcal{A}_3$ and $\min_{e\in
  \vec{E}(f_0)}\un{W}_e\geq -2\eps(8n/9)^{1/4}$ then $(\sg,(W_e)_{e\in
  \vec{E}(\sg)})$ belongs to
$(\psi^{(4)}_n)^{-1}(\mathcal{B}_2(\eps))$. To be perfectly accurate,
there is a small difference coming from the $+2$ term in
\eqref{eq:50}, which does not appear anymore in \eqref{eq:44}. A way
to circumvent this would be to re-define the discrete snake processes
associated with a labeled map $(\bm,\bl)$ by shifting them by
$1$. Such a modification obviously does not change the limit theorems
of Section \ref{sec:scal-limits-label}.

By similar reasoning as in the derivation of \eqref{eq:48}, using
Proposition \ref{sec:limit-theorems-1} rather than Proposition
\ref{sec:case-planted-schemes-1}, we deduce from the above discussion
that
\begin{equation}\label{eq:49}
\limsup_{n\to\infty}P(\mathcal{A}^{(n)}_2(\eps))\leq
\frac{C}{\eps}\Big(\CLM_1^{(4)}(\sg\notin\mathfrak{P},\min_{e\in
  \vec{E}(f_0)}\un{W}_e\geq -2\eps)+\CLM_1^{(4)}(\mathcal{B}_2(\eps))\Big)\, .
\end{equation}

\subsection{Some estimates for bridges and
  snakes}\label{sec:some-estimates}

Here we gather the technical estimates that will be needed to estimate
the upper-bounds in \eqref{eq:48} and \eqref{eq:49}. Recall the
notation from Section \ref{sec:bridge-measures} and \ref{sec:snakes}
for bridge and snake measures.

\begin{lmm}\label{sec:some-estimates-2}
Fix $\lambda>0$. Then 
\begin{equation}
  \label{eq:29}
\sup_{x,y>0}\B_{x\to y}[e^{-\lambda \zeta(X)}]<\infty \qquad
\mbox{and} \qquad 
\sup_{x\in \R}\int_\R \d y\, \B_{x\to y}[e^{-\lambda \zeta(X)}]<\infty\,
. 
\end{equation}
Moreover, for every $x,y\geq
0$, it holds that 
\begin{equation}
  \label{eq:30}
  \B_{x\to y}[e^{-\lambda \zeta(X)}\ind_{\{\un{X}\geq 0\}}]\leq
  2(x\wedge y)\, ,
\end{equation}
and
\begin{equation}
  \label{eq:31}
   \B_{x\to y}[\zeta(X)e^{-\lambda \zeta(X)}\ind_{\{\un{X}\geq 0\}}]\leq
  \sqrt{\frac{2}{\lambda}}\, xy\, .
\end{equation}
\end{lmm}

\proof 
On the one hand, we have
$$\B_{x\to y}[e^{-\lambda \zeta(X)}]=\int_0^\infty \d r\,
p_r(x,y)e^{-\lambda r}\leq \int_0^\infty\frac{\d r\, e^{-\lambda
    r}}{\sqrt{2\pi r}}\, ,$$ which is finite and independent of
$x,y$. Moreover, by Fubini's theorem and the fact that $p_r(x,y)$ is a
probability density, we have
$$\int_{\R}\d y\, \B_{x\to y}[e^{-\lambda\zeta(X)}]=\int_0^\infty\d
r\, e^{-\lambda r}\, ,$$ which is again independent of $x$. This gives
\eqref{eq:29}.

Next, as in Section \ref{sec:bridge-measures}, we use the following
consequence of the reflection principle: 
$$\P_{x\to y}^r(\un{X}\geq 0)=1-e^{-2xy/r}\, ,\qquad \mbox{ for every
}\quad x,y\geq 0\mbox{ and} r>0\, .$$
Assuming first that $x,y>0$ and $x\neq y$, this gives
\begin{eqnarray*}
  \mathbb{B}_{x\to y}(e^{-\lambda\zeta(X)}\ind_{\{\un{X}\geq
    0\}})&=&\int_0^\infty \d r\, p_r(x,y)e^{-\lambda r}-\int_0^{\infty}\d r\,
  p_r(x,y)e^{-\lambda r}e^{-2xy/r}\\
  &=&
  \E_0\Big[\frac{T_{|x-y|}}{|x-y|}e^{-\lambda
    T_{|x-y|}}\Big]-\E_0\Big[\frac{T_{x+y}}{x+y}e^{-\lambda T_{x+y}}\Big]\, ,
\end{eqnarray*}
where we used the well-known fact that $(a/r)p_r(0,a)\d r=\P_0(T_a\in \d
r)$. From the Laplace transform of $T_a$, given by $\E_0[e^{-u
  T_a}]=e^{-a\sqrt{2u}}$, we immediately get that for $a>0$, 
$$\E_0\Big[\frac{T_a}{a}e^{-u T_a}\Big]=\frac{1}{\sqrt{2u}}e^{-a\sqrt{2u}}\,
,$$
from which we obtain
$$\mathbb{B}_{x\to y}(e^{-\lambda\zeta(X)}\ind_{\{\un{X}\geq
  0\}})=\frac{e^{-|x-y|\sqrt{2\lambda}}}{\sqrt{2\lambda}}(1-e^{-2\sqrt{2\lambda}(x\wedge
  y)})\leq 2(x\wedge y)\, ,$$ as wanted. This remains true for $x=y$
or for $xy=0$ by a continuity argument, yielding \eqref{eq:30}.  The
proof of \eqref{eq:31} is similar, writing
\begin{eqnarray*}
  \mathbb{B}_{x\to y}(\zeta(X)e^{-\lambda \zeta(X)}\ind_{\{\un{X}\geq
    0\}})&=&\int_0^\infty \d r\, r\, p_r(x,y)e^{-r}-\int_0^{\infty}\d r\,
  r\, p_r(x,y)e^{-r}e^{-2xy/r}\\
  &=&
  \E_0\Big[\frac{T_{|x-y|}^2}{|x-y|}e^{-\lambda
    T_{|x-y|}}\Big]-\E_0\Big[\frac{T_{x+y}^2}{x+y}e^{-\lambda T_{x+y}}\Big]\, ,
\end{eqnarray*}
and using again Laplace transforms to get, for $a,u>0$, 
$$\E_0\Big[\frac{T_a^2}{a}e^{-uT_a}\Big]=\frac{1}{2\sqrt{2u^3}}(a\sqrt{2u}+1)e^{-a\sqrt{2ua}}\,
.$$ Now assume without loss of generality that $0<y\leq x$ and write
\begin{eqnarray*}
  \lefteqn{\mathbb{B}_{x\to y}(\zeta(X)e^{-\lambda\zeta(X)}\ind_{\{\un{X}\geq
      0\}})}\\
  &=&\frac{1}{2\sqrt{2\lambda^3}}\bigg(((x-y) \sqrt{2\lambda}+1)e^{(-x+y)\sqrt{2\lambda}}-(  (x+y) \sqrt{2\lambda}+1)e^{-(x+y)\sqrt{2\lambda}}\bigg)\\
  &=&
  \frac{e^{-x\sqrt{2\lambda}}\cosh(y\sqrt{2\lambda})}{\sqrt{2\lambda^3}}\bigg((
  x \sqrt{2\lambda}+1)\tanh(y\sqrt{2\lambda})- y \sqrt{2\lambda}\bigg)\\
  &\leq & \sqrt{\frac{2}{\lambda}}xy\, ,
\end{eqnarray*}
as wanted. \cq

\begin{lmm}\label{sec:some-estimates-1}
Let $x,y,z$ be positive real numbers. Then 
\begin{equation}
  \label{eq:32}
  \E_x^{(y,\infty)}[\Q_X(\un{W}\geq 0)]\leq\Big(
  \frac{y}{x}\Big)^2\qquad \mbox{ if }\quad 0<y<x\, ,
\end{equation}
and there exists a finite $C>0$ such that
\begin{equation}
  \label{eq:33}
  \E_x^{(y,\infty)}[\Q_X(\un{W}\geq 0)\Q_X(\un{W} < -z)]\leq
  Cy^2\Big(\frac{1}{z} \wedge \frac{1}{x}\Big)^2\qquad \mbox{ if }\quad 0<y<x\, .
\end{equation}
Finally, for every $\beta\in [-2,3]$, it holds that 
\begin{equation}
  \label{eq:34}
  \B_{x\to y}[\Q_X(\un{W}\geq 0)]\leq \frac{1}{5}x^\beta y^{1-\beta}\, .
\end{equation}
\end{lmm}

\proof We use the Poisson point process description of the Brownian
snake. Namely, recall the notation from Section \ref{sec:snakes} and
the fact that the snake $W$ under $\Q_X$ can be decomposed in
excursions $W^{(r)},0\leq r\leq \zeta(X)$, in such a way that
$$\sum_{0\leq r\leq \zeta(X)
}\delta_{(r,W^{(r)})}\ind_{\{T_r>T_{r-}\}}$$ is a Poisson random
measure on $[0,\zeta(X)]\times \CC(\R)$, with intensity measure given
by
$$2\, \d r\, \ind_{[0,\zeta(X)]}(r)\, \N_{X(r)}(\d W)\, .$$
From this and the known formula \cite{legweill}
$$\N_0(\un{W}<-y)=\frac{3}{2y^2}\, ,\qquad y>0\, ,$$
we obtain, using standard properties of Poisson measures,
\begin{eqnarray}
  \Q_X(\un{W}\geq
  0)&=&\exp\bigg(-2\int_0^{\zeta(X)}\!\!\!\d
  r\, \N_{X(r)}(\un{W}<0)\bigg)\nonumber\\
  &=&\exp\bigg(-2\int_0^{\zeta(X)}\!\!\!\d
  r  \, \N_{0}(\un{W}<-X(r))\bigg)\nonumber\\
  &=&\exp\bigg(-\int_0^{\zeta(X)}\frac{3\, \d
    r}{X(r)^2}\bigg)\, .\label{eq:24}
\end{eqnarray}
We deduce that for every $x>y>0$,
$$\E^{(y,\infty)}_x[\Q_X(\un{W}\geq 0)]=
 \E_x^{(y,\infty)}\bigg[\exp\bigg(-\int_0^{T_y}\frac{3\, \d 
r}{X(r)^2}\bigg)\bigg]\, .
$$
Recalling that reflected Brownian motion is a $1$-dimensional Bessel
process, we now use the absolute continuity relations between Bessel
processes with different indices, due to Yor \cite[Exercise
XI.1.22]{revyor} (see also \cite{legweill} for a similar use of these
absolute continuity relations).  The last expectation then equals
$$\Big(\frac{x}{y}\Big)^3\P_x^{\langle 7\rangle}(T_y<\infty)\, ,$$
where, for $\delta\geq 0$, $\P_x^{\langle\delta\rangle}$ is the law of the
$\delta$-dimensional Bessel process started from $x>0$. Recall that
for $\delta\geq 2$, the $\delta$-dimensional Bessel process started
from $x>0$ is the strong solution (starting from $x$) of the
stochastic differential equation driven by the standard Brownian
motion $(B_t,t\geq 0)$:
$$\d Y_t=\d B_t+\frac{\delta-1}{2Y_t}\d t\, .$$
One can show that $Y_t>0$ for all $t$, hence that the drift term is
well-defined, whenever $\delta\geq 2$. Showing that
$\P^{\langle\delta\rangle}_x(T_y<\infty)=(y/x)^{\delta-2}$ for every
$x>y>0$ is now a simple exercise, using the fact that
$(Y_t^{2-\delta},t\geq 0)$ is a local martingale by Itô's
formula. Putting things together, we get
\eqref{eq:32}. 

Let us now turn to \eqref{eq:33}. By \eqref{eq:24} and an easy
translation invariance argument, we have
\begin{eqnarray*} 
\lefteqn{\E_x^{(y,\infty)}[\Q_X(\un{W}\geq 0)\Q_X(\un{W}< -z)]}\\
&=&
 \E_x^{(y,\infty)}\bigg[\exp\Big(-\int_0^{T_y}\frac{3\, \d 
r}{X(r)^2}\Big)\bigg(1-\exp\Big(-\int_0^{T_y}\frac{3\, \d 
r}{(X(r)+z)^2}\Big)\bigg)\bigg]\\
&\leq & 3\,  \E_x^{(y,\infty)}\bigg[\exp\Big(-\int_0^{T_y}\frac{3\, \d 
r}{X(r)^2}\Big)\Big(1\wedge \frac{T_y}{z^2}\Big)\bigg]\\
&=&3\Big(\frac{x}{y}\Big)^3  
\E_x^{\langle 7\rangle}\bigg[\ind_{\{T_y<\infty\}}\Big(1\wedge
\frac{T_y}{z^2}\Big)\bigg]\\
&\leq & 3\Big(\frac{x}{y}\Big)^3  
\Big(\P_x^{\langle 7\rangle}(T_y<\infty)\wedge
\frac{1}{z^2}\E_x^{\langle 7\rangle}[T_y\ind_{\{T_y<\infty\}}]\Big) \, ,
\end{eqnarray*}
where we have used again the absolute continuity relations for Bessel
processes at the third step. We already showed that $\P_x^{\langle
  7\rangle}(T_y<\infty)=(y/x)^5$, so to conclude it suffices to show
that $\E_x^{\langle 7\rangle}[T_y\ind_{\{T_y<\infty\}}]\leq Cy^5/x^3$
for some finite constant $C$. By the Markov property, and
using again the formula for the probability that $T_y<\infty$ for the
$7$-dimensional Bessel process, we have
\begin{eqnarray*}
  \E_{x}^{\langle 7\rangle}[T_y\ind_{\{T_y<\infty\}}] & =&
  \int_0^\infty\d s\, \P_{x}^{\langle 7\rangle}(s<T_y<\infty)\\
  &=& \int_0^\infty\d s\, \E_{x}^{\langle
    7\rangle}\big[\ind_{\{s<T_y\}}\P_{X(s)}^{\langle
    7\rangle}(T_y<\infty)\big]\\
  &\leq & \int_0^\infty\d s\, \E_x^{\langle
    7\rangle}\Big[\Big(\frac{y}{X(s)}\Big)^5\Big]\, .
\end{eqnarray*}
Using the fact that the $7$-dimensional Bessel process has same
distribution as the Euclidean norm of the $7$-dimensional Brownian
motion, and the known form of the latter's Green function, we obtain
that if $u=(1,0,0,0,0,0,0)\in \R^7$, it holds that 
$$\E_x^{\langle 7\rangle}[T_y\ind_{\{T_y<\infty\}}] \leq
\int_{\R^7}\frac{\d
  z}{|z-xu|^5}\cdot\Big(\frac{y}{|z|}\Big)^5=\frac{y^5}{x^3}\int_{\R^7}\frac{\d
z}{|z-u|^5|z|^5}\, ,$$
and the integral is finite, as wanted. 

We now prove \eqref{eq:34}, by using the agreement formula
\eqref{eq:27}, entailing that 
\begin{eqnarray*}
  \B_{x\to y}[\Q_X(\un{W}\geq 0)]&=&\int_{-\infty}^{x\wedge y}\d z\,
  (\E_x^{(z,\infty)}\bowtie\widehat{\E}_y^{(z,\infty)})[\Q_X(\un{W}\geq
  0)]\\
  &=&\int_0^{x\wedge y}\d z\, \E_x^{(z,\infty)}[\Q_X(\un{W}\geq 0)] 
\E_y^{(z,\infty)}[\Q_X(\un{W}\geq 0)]\\
&\leq &\int_0^{x\wedge y}\d z\,
\Big(\frac{z}{x}\Big)^2\cdot\Big(\frac{z}{y}\Big)^2=\frac{(x\wedge
  y)^5}{5x^2y^2}\, ,
\end{eqnarray*}
where we used \eqref{eq:32} in the penultimate step. The conclusion follows
easily.  \cq

\medskip

Let us consider once again the Poisson point measure representation
$(W^{(t)},0\leq t\leq \zeta(X))$ of the Brownian snake under the law
$\Q_X$, as explained around \eqref{eq:17}. Fix $y>0$. We are
interested in the distribution of the random variable $\inf_{0\leq
  t\leq T_{-z}}\un{W}^{(t)}$ for $0\leq z\leq y$, as well as in
bounding expectations of the form
$$\E_0^{(-y,\infty)}\Big[\Q_X\Big(-x\wedge \inf_{0\leq r\leq T_{-z}}\un{W}^{(r)}\leq
-2z\mbox{ for every }z\in [0,y]\Big)\Big]\, .$$ To this end, we
perform yet another Poisson measure representation for these
quantities.  

\begin{lmm}\label{sec:some-estimates-3}
  Let $y>0$ be fixed. Let $X$ be the canonical process on $\CC(\R)$,
  and $W$ be the canonical process on $\CC(\CC(\R))$ started from
  $W(0)=X$. Recall that $T_x=\inf\{t\geq 0:X(t)=x\}$ for every $x\in
  \R$. For every $z\geq 0$, let 
$$I_z= -z-\inf\{\un{W}^{(t)} : T_{(-z)-}\leq t\leq T_{-z}\}\, ,$$
which is taken to be $0$ by convention if $T_{-z}=\infty$.  Then under
$\P^{(-y,\infty)}_0(\d X)\Q_X(\d W)$, the point measure
$$\sum_{0\leq z\leq y}\delta_{(z,I_z)}\ind_{\{I_z>0\}}$$
 is a Poisson random measure on
$[0,y]\times \R_+$ with intensity $\d z\ind_{[0,y]}(z)\otimes 2\d
a/a^2$.
\end{lmm}

\proof By Itô's excursion theory, under the distribution
$\E^{(-y,\infty)}_0(\d X)$, if we let
$$X^{(z)}=z+X(T_{(-z)-}+t)\, ,\qquad 0\leq t\leq T_{-z}-T_{(-z)-}$$
for every $z\in (0,y)$ such that $T_{-z}>T_{(-z)-}$, then the measure
$$\sum_{0\leq z\leq y}\delta_{(z,X^{(z)})}\ind_{\{T_{(-z)}>T_{(-z)-}\}}$$
is a Poisson point measure on $\R_+\times \CC(\R)$ with intensity $\d
z\ind_{(0,y)}(z)\otimes 2n(\d X)$. For every $z\in [0,y]$ such that
$T_{(-z)}>T_{(-z)-}$, we can interpret $(z+W^{(t+T_{(-z)-})},0\leq t\leq
T_{-z}-T_{(-z)-})$ as an independent mark on the excursion $X^{(z)}$, given by a
snake with distribution $\Q_{X^{(z)}}(\d W)$. By the marking
properties for Poisson measures and symmetry, we obtain that
$\sum_{0\leq z\leq y}\delta_{(z,I_z)}$ is itself a Poisson measure
with intensity
$$ \d z\ind_{\{0\leq z\leq y\}}\otimes \int_{\CC(\R)}2n(\d
X)\Q_X(-\un{W}\in \d a)\, .$$ Now, we use \eqref{eq:24} and the fact
(\cite[Exercise XII.2.13]{revyor}) that the image measure of $n(\d X)$
under the scaling operation $X\mapsto a^{-1}X(a^2\cdot)$ is
$a^{-1}n(\d X)$, to get
\begin{eqnarray*}
\int_{\CC(\R)}2n(\d
X)\Q_X(-\un{W}>a)&=&\int_{\CC(\R)}2n(\d
X)\bigg(1-\exp\Big(-\int_0^{\zeta(X)}\frac{3\d
  r}{(X(r)+a)^2}\Big)\bigg)\\
& = &\int_{\CC(\R)}2n(\d X) \bigg(1-\exp\Big(-\int_0^{\zeta(X)/a^2}\frac{3\d
  r}{(a^{-1}X(a^2r)+1)^2}\Big)\bigg)\\
&=& \frac{K}{a}\, ,
\end{eqnarray*}
where 
$$K=\int_{\CC(\R)}2n(\d X) \bigg(1-\exp\Big(-\int_0^{\zeta(X)}\frac{3\d
  r}{(X(r)+1)^2}\Big)\bigg)\, .$$ To compute it explicitly, write 
$$K=\int_{\CC(\R)}2n(\d X)\int_0^{\zeta(X)} \frac{3\d t}{(X(t)+1)^2}\exp\Big(-\int_t^{\zeta(X)}\frac{3\d
  r}{(X(r)+1)^2}\Big)$$ and use the Bismut decomposition \cite[Theorem
XII.4.7]{revyor} to obtain
$$K=6\int_0^\infty \frac{\d
  a}{(a+1)^2}\E_a^{(0,\infty)}\Big[\exp\Big(-\int_0^{T_0}\frac{3\d
  r}{(X(r)+1)^2}\Big)\Big]\, .$$ By translating the process $X$ by $1$
and arguing as in the proof of Lemma \ref{sec:some-estimates-1}, we
obtain that the expectation inside the integral equals
$(a+1)^{-2}$. This ends the proof. \cq

\medskip

By definition of $I_z,z\geq 0$, we have that the process $(\inf_{0\leq
  r\leq T_{-z}}\un{W}^{(r)},0\leq z\leq y)$ under $\P^{(-y,\infty)}_z$
is equal to $(-\sup_{0\leq r\leq z}I_r,0\leq z\leq y)$.  Dealing with
such random variables and processes is classical in extreme values
theory \cite{resnick87}. In the sequel, we will let
$\sum_{z}\delta_{(z,\Delta_z)}\ind_{\{\Delta_z>0\}}$ be a Poisson
random measure on $\R_+\times (0,\infty)$ with intensity $\d z\otimes
K\d a/a^2$, for some $K>0$ (not necessarily equal to $2$), and defined
on some probability space $(\Omega,\FF,P)$.  We also let
$\ov{\Delta}_z=\sup_{0\leq r\leq z}\Delta_r$ (the process
$\ov{\Delta}$ is called a {\em record process}). Note that for every
$t>0$, the process $\ov{\Delta}$ remains constant equal to
$\ov{\Delta}_t$ on a small neighborhood to the right of $t$, and that
infinitely many jumps times accumulate near $t=0$. Hence, the process
$(\ov{\Delta}_t,t>0)$ is a jump-hold process.

By standard properties of Poisson random measures, the one-dimensional
marginal laws of this process are so-called {\em Fr\'echet laws},
given by
$$P(\ov{\Delta}_t\leq x)=\exp\Big(-\frac{Kt}{x}\Big)\, ,\qquad
t,x\geq 0\, .$$
Moreover, the process $\ov{\Delta}$ satisfies the homogeneous scaling
relation 
$$\Big(\frac{\ov{\Delta}_{at}}{a},t\geq
0\Big)\build=_{}^{(d)}(\ov{\Delta}_t,t\geq 0)\, .$$

For $x,t\geq 0$, consider the event
$$H(x,t)=\{x\vee \ov{\Delta}_s\geq s\, ,\, 0\leq s\leq
t\}=\{\ov{\Delta}_s\geq s\, ,\, x\leq s\leq t\}\, .$$ By scaling, we
have $P(H(x,t))=P(H(x/t,1))$.

\begin{lmm}\label{sec:fast-coal-geod}
For every
  $0\leq x\leq 1$, it holds that 
$$P(H(x,1))\leq x^{e^{-K}/K}\, .$$
\end{lmm}

\proof Let $(J_n,D_n,n\in \Z)$ be a (measurable) enumeration of the jump
times of the process $\ov{\Delta}$, and values of $\ov{\Delta}$ at these
jump times, in such a way that 
$$\ldots <J_{-2}<J_{-1}<J_0<J_1<J_2<\ldots\, ,\qquad
D_n=\ov{\Delta}_{J_n}\, ,\qquad n\in \Z\, .$$ In particular, note that
$J_n$ is always the coordinate of the first component of the Poisson
measure used to construct $\ov{\Delta}$, and $D_n$ is the
corresponding second component, since at time $J_n$, by definition,
the process $\ov{\Delta}$ achieves a new record.

We use the fact \cite{resnick87} that the measure
$\mathcal{M}=\sum_{n\in \Z}\delta_{(D_n,J_{n+1}-J_n)}$ is a Poisson
random measure on $(0,\infty)^2$, with intensity measure given by
$$\mu(\d y\d u)=\frac{e^{-Ku/y}}{y^2}\d y\d u\, .$$
Note that on the event $H(x,1)$, it must hold that for every $n$ such
that $D_n\in [x,1]$, we have $J_{n+1}-J_n\leq D_n$ (since otherwise,
we have $\ov{\Delta}_{J_{n+1}-}=D_n<J_{n+1}$ with $D_n\in [x,1]$, so
$H(x,1)$ cannot hold). This means that the measure $\mathcal{M}$ has no
atom in $\{(y,u):y\in[x,1],u>y\}$.  Therefore, we have
\begin{eqnarray*}
  P(H(x,1))&\leq &\exp(-\mu(\{(y,u):y\in [x,1],u>y\}))\\
  &=& \exp(K^{-1}e^{-K}\log(x))\, ,
\end{eqnarray*}
as wanted.  \cq

\subsection{Fast confluence of geodesics}\label{sec:fast-confl-geod}

We now turn to the proof of Lemma \ref{sec:covering-3-star}. 

\begin{figure}[htb!]
\begin{center}
\includegraphics{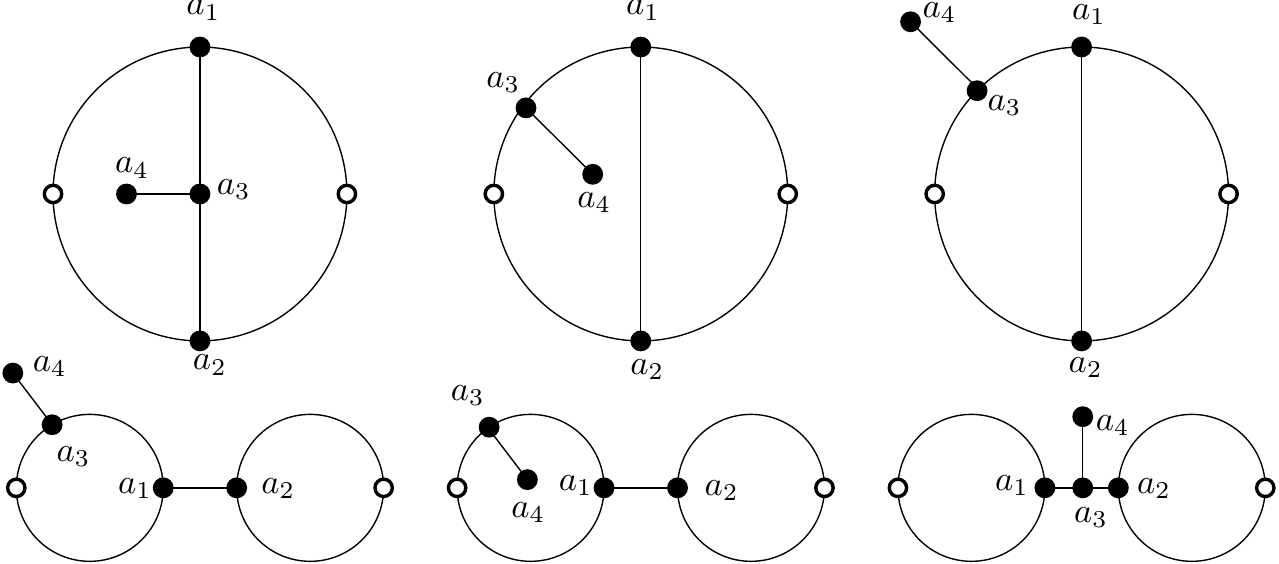}
\end{center}
\caption{The dominant planted schemes with $3$ faces, where $f_0$ is
  the outside face, considered up re-labeling of faces, and obvious
  symmetries. Elements of $V_N(\sg)$ are indicated by blank vertices.  The
  first one is a predominant scheme, the other ones are not. The
  notation $a_1,a_2,a_3,a_4$ refers to the integration variables
  appearing in the proof of Lemmas \ref{sec:fast-coal-geod-1} and
  \ref{sec:fast-coal-geod-2}. }
\label{fig:schemes2}
\end{figure}

If $\sg$ is a planted dominant scheme, we call it predominant if the
degree of $f_0$ is minimal, as for non-planted schemes. Figure
\ref{fig:schemes2} displays in first position the possible predominant
schemes (without face labels: There are two possible such labelings
depending on the location of $f_1$ and $f_2$ as inside faces).  


\begin{lmm}\label{sec:fast-coal-geod-1}There exists a constant $C\in
  (0,\infty)$ such that for every $\eps>0$,
  $$\dot{\CLM}_1^{(3)}\Big(\sg\notin\mathfrak{P},\min_{e\in \vec{E}(f_0)}\un{W}_e\geq
  -\eps, \min_{e\in E(f_1\cap f_2)}\un{M}_e\geq 0\Big)\leq C\eps^5\,
  .$$
\end{lmm}

\proof 
By an elementary scaling argument using Proposition
\ref{sec:cont-meas-label-1}, we have,
\begin{eqnarray}
  \lefteqn{
\dot{\CLM}^{(3)}_1\Big(\sg\notin \mathfrak{P},\min_{e\in \vec{E}(f_0)}\un{W}_e\geq
  -\eps,\min_{e\in E(f_1\cap f_2)}\un{M}_e\geq 0\Big)}\label{eq:51}\\
&\leq& 
  C\cdot\dot{\CLM}^{(3)}\Big(\ind_{[1/2,1]}(\sigma)\, ;\,
  \sg\notin \mathfrak{P}\, ,\, \min_{e\in \vec{E}(f_0)}\un{W}_e\geq
  -\eps,\min_{e\in E(f_1\cap f_2)}\un{M}_e\geq 0\Big)\nonumber\\
  &\leq & C\cdot \dot{\CLM}^{(3)}\Big(e^{-\sigma(W)}\, ;\, \sg\notin
  \mathfrak{P}\, ,\, 
  \min_{e\in \vec{E}(f_0)}\un{W}_e\geq
  -\eps,\min_{e\in E(f_1\cap f_2)}\un{M}_e\geq 0\Big)\, .\nonumber
\end{eqnarray}
In this form, we can take advantage of the fact that
$\dot{\CLM}^{(3)}$ is a sum over $\dot{\bS}^{(3)}_d$ of product
measures, as expressed in \eqref{eq:19}. The contributions of the second and third
schemes in Figure \ref{fig:schemes2} to the above upper-bound are then at most
\begin{eqnarray*}
  \lefteqn{C\int_{\R_+^4}\d a_1\, \d a_2\, \d a_3\, \d a_4
    \E_{a_1}^{(0,\infty)}[\Q_X(\un{W}\geq -\eps)]
    \E_{a_3}^{(0,\infty)}[\Q_X(\un{W}\geq -\eps)]
   }
\\
  & &\times \E_{a_2}^{(0,\infty)}[\Q_X(\un{W}\geq -\eps)]^2
\B_{a_1\to a_3}(\un{W}\geq -2\eps)
\mathbb{B}_{a_1\to
    a_2}(e^{-\zeta(X)}\ind_{\{\un{X}\geq 0\}})\mathbb{B}_{a_3\to
    a_4}(e^{-\zeta(X)}) 
\end{eqnarray*}
In this expression, we can integrate out the variable $a_4$ using the
second expression of \eqref{eq:29}.  We bound the terms of the form
$\E^{(0,\infty)}_a[\Q_X(\un{W}\geq -\eps)]$ by using \eqref{eq:32},
and the terms involving measures $\B_{a\to b}$ by using \eqref{eq:34}
with $\beta=1/3$, and \eqref{eq:30} together with the fact that
$a_1\wedge a_2\leq \sqrt{a_1 a_2}$. This gives an upper-bound 
\begin{eqnarray*}
C\int_{\R_+^3}\d a_1\, \d a_2\, \d a_3
  \Big(\frac{\eps}{\eps+a_1}\Big)^2
  \Big(\frac{\eps}{\eps+a_2}\Big)^4\Big(\frac{\eps}{\eps+a_3}\Big)^2
  (a_1+\eps)^{1/3}(a_3+\eps)^{2/3}\sqrt{a_1 a_2}  \\
\leq C\eps^5\int_{\R_+^3}\frac{\d a_1\, \d a_2\, \d
a_3}{(1+a_1)^{7/6}(1+a_2)^{7/2}(1+a_3)^{4/3}}\, ,
\end{eqnarray*}
and the integral is finite, as wanted.

It remains to evaluate the contribution of the bottom three schemes of
Figure \ref{fig:schemes2}. These do not have edges that are incident
both to $f_1$ and $f_2$, so that the condition $\min_{e\in E(f_1\cap
  f_2)}\un{M}_e\geq 0$ can be removed. The contribution of the
fourth and fifth schemes is then bounded by (we skip some intermediate
steps analogous to the above) 
\begin{eqnarray*}
  \lefteqn{C\int_{\R_+^4}\d a_1\, \d a_2\, \d a_3\, \d a_4
    \Big(\frac{\eps}{\eps+a_1}\Big)^2
    \Big(\frac{\eps}{\eps+a_2}\Big)^4\Big(\frac{\eps}{\eps+a_3}\Big)^2 }\\
  & &\times (a_1+\eps)^\beta(a_2+\eps)^{1-\beta}(a_1+\eps)^{\beta'}(a_3+\eps)^{1-\beta'}\mathbb{B}_{a_3\to
    a_4}(e^{-\zeta(X)}) \\
  &\leq & C\int_{\R_+^3}\d a_1\, \d a_2\, \d a_3
  \Big(\frac{\eps}{\eps+a_1}\Big)^2
  \Big(\frac{\eps}{\eps+a_2}\Big)^4\Big(\frac{\eps}{\eps+a_3}\Big)^2
  \\
  &&\times 
  (a_1+\eps)^{-1}(a_2+\eps)^2(a_1+\eps)^{1/2}(a_3+\eps)^{1/2}
\end{eqnarray*}
taking $\beta=-1$ and $\beta'=1/2$, and we conclude similarly.  The
contribution of the sixth scheme of Figure \ref{fig:schemes2} is
bounded above by
\begin{eqnarray*}
C\int_{\R_+^4}\d a_1\, \d a_2\, \d a_3\, \d a_4
    \Big(\frac{\eps}{\eps+a_1}\Big)^4
    \Big(\frac{\eps}{\eps+a_2}\Big)^4(a_1+\eps)^{\beta}(a_2+\eps)^{\beta'}
(a_3+\eps)^{2-\beta-\beta'}\mathbb{B}_{a_3\to
      a_4}(e^{-\zeta(X)}) \\
  \leq  C\int_{\R_+^3}\d a_1\, \d a_2\, \d a_3
  \Big(\frac{\eps}{\eps+a_1}\Big)^4
  \Big(\frac{\eps}{\eps+a_2}\Big)^4\Big(\frac{\eps}{\eps+a_3}\Big)^2 
  (a_1+\eps)^2(a_2+\eps)^2(a_3+\eps)^{-2}\, ,
\end{eqnarray*}
taking $\beta=\beta'=2$. Once again, we conclude
in the same way as above.
\cq

\begin{lmm}
  \label{sec:fast-coal-geod-2}
  For every $c,\beta>0$, there exists a finite $C>0$ such that for every
  $\eps>0$,
\begin{equation}
  \label{eq:23}
  \dot{\CLM}_1^{(3)}\Big(\sg\in \mathfrak{P},\min_{e\in \vec{E}(f_0)}\un{W}_e\geq -\eps,
  \min_{e\in E(f_1\cap f_2)}\un{M}_e\geq 0,\exists e\in E_I(\sg),\un{W}_e\leq
  -c\eps^{1-\beta}\Big)\leq C\eps^{4+\beta}\, , 
\end{equation}
where in the last part of the event, it is implicit that $e$ is
oriented so that it is incident to $f_1$ or $f_2$, but not $f_0$.
\end{lmm}

\proof By using the same scaling argument as in \eqref{eq:51}, up to
changing $c$ by $2^{1/4}c$ and $C$ by a larger constant, it suffices
to prove a similar bound for the measure
$\dot{\CLM}^{(3)}(e^{-\sigma}\cdot)$ instead of $\dot{\CLM}_1^{(3)}$,
and use formula \eqref{eq:19} to estimate this quantity. Furthermore,
for obvious symmetry reasons, up to increasing the constant $C$ by a
factor $8$ in the end, it suffices to estimate the contribution of the
scheme $\sg$ corresponding to the first picture of Figure
\ref{fig:schemes2}, where $f_1$ is the left internal face and $f_2$ is
the right internal face, and specifying that among the four edges
$e\in E_I(\sg)$, the top-left one is such that $\un{W}_e\leq
-\eps^{1-\beta}$, the others being unconstrained. The contribution of
this event to \eqref{eq:23} is then bounded by
\begin{eqnarray*}
  \lefteqn{C\int_{\R_+^2}\d a_1\, \d a_2
    \Big(\frac{\eps}{a_1+\eps}\Big)^2
    \Big(\frac{\eps}{a_2+\eps}\Big)^4   
    \E_{a_1}^{(0,\infty)}[\Q_X(\un{W}\geq -\eps)\Q_X(\un{W}\leq
    -\eps^{1-\beta})] }\\
  && \times\int_{\R}\d a_3\, 
  \B_{a_2\to a_3}[e^{-\zeta(X)}\ind_{\{\un{X}\geq 0\}}]\B_{a_1\to a_3}[e^{-\zeta(X)}\ind_{\{\un{X}\geq
    0\}}] 
  \int_{\R} d a_4\,  \B_{a_3\to
    a_4}[e^{-\zeta(X)}]\, .
\end{eqnarray*}
The last integral is bounded by \eqref{eq:29}, independently on
$a_1,a_2,a_3$. We then apply \eqref{eq:28} to the integral in the
variable $a_3$, which gives rise to the factor $\B_{a_1\to
  a_2}[\zeta(X)e^{-\zeta(X)\sqrt{2}}\ind_{\{\un{X}\geq 0\}}]$, and
this is bounded by $Ca_1a_2$ by \eqref{eq:31}. After an elementary
change of variables and a translation and scaling in last remaining
expectation, we obtain the bound 
$$C\eps^4\int_{\R_+^2}\d a_1\, \d
a_2\frac{a_1a_2}{(a_1+1)^2(a_2+1)^4}\E_{a_1+1}^{(1,\infty)}[\Q_X(\un{W}\geq
0)\Q_X(\un{W}\leq -\eps^{-\beta}+1)]\, ,$$ in which the contribution
of $a_2$ can be integrated out. Therefore, taking $\eps$ small enough
so that $\eps^{-\beta}-1>\eps^{-\beta}/2$ and using \eqref{eq:33}, this
bound is less than or equal to 
\begin{eqnarray*}
  \lefteqn{C\eps^4\int_0^\infty\d
    a_1\frac{a_1}{(a_1+1)^2}\Big(\frac{1}{a_1+1}\wedge
    \frac{1}{\eps^{-\beta}-1}\Big)^2}\\
  &\leq & C\eps^4\bigg(\eps^{2\beta}\int_0^{\eps^{-\beta}}\frac{\d
    a_1}{a_1+1}+\int_{\eps^{-\beta}}^\infty\frac{\d a_1}{(a_1+1)^3}\bigg)
  \\
  &\leq & C\eps^4(\eps^{2\beta}\log(1+\eps^{-\beta}) +
  (\eps^{-\beta}+1)^{-2})\\
&\leq & C\eps^{4+\beta}
\end{eqnarray*}
This bound remains true for every $\eps>0$, possibly up to changing the
constant $C$. Note that we could have in fact obtained a bound of the
form $C\eps^{4+2\beta'}$ for any $\beta'\in(0,\beta)$ with this
method. \cq

\medskip

Combining Lemmas \ref{sec:fast-coal-geod-1} and
\ref{sec:fast-coal-geod-2} with \eqref{eq:48}, this completes the
proof of Lemma \ref{sec:covering-3-star} (note that we changed the
second condition in the definition of the event
$\mathcal{B}_1(\eps,\beta)$ by $\min_{e\in \vec{E}(f_0)}\un{W}_e\geq
-\eps$, this was for lightening the notation but is of no impact, as
is easily checked).

\subsection{$\eps$-geodesic stars}\label{sec:eps-geodesic-stars}

We finally prove Lemma \ref{sec:covering-3-star-1}. In these proofs,
up to considering $\mathcal{B}_2(\eps/2)$ instead of
$\mathcal{B}_2(\eps)$, we will replace the condition that $\min_{e\in
  \vec{E}(\sg)}\un{W}_e\geq -2\eps$ by the similar bound with $-\eps$
without impacting the result.

\begin{lmm}\label{sec:proof-lemm-refs}
There exists some constant $C>0$ such that 
$$\CLM^{(4)}_1\Big(\sg\notin \mathfrak{P},\min_{e\in \vec{E}(f_0)}\un{W}_e\geq
-\eps\Big)\leq C\eps^5\, .$$
\end{lmm}

\proof As in Lemma \ref{sec:fast-coal-geod-1}, a scaling argument
shows that it suffices to prove the same bound with the measure
$\CLM^{(4)}(e^{-\sigma}\cdot)$ replacing $\CLM^{(4)}_1$. 
Note that for any $k\geq 2$, we have
\begin{eqnarray*}
  \lefteqn{\CLM^{(k+1)}\Big(e^{-\sigma}\, ;\, \sg\notin \mathfrak{P}\, ,\, \min_{e\in \vec{E}(f_0)}\un{W}_e\geq
    -\eps\Big)}\\ & = & \sum_{\sg\notin \mathfrak{P}}\int
  \lambda_\sg(\d(\ell_v)_{v\in V(\sg)})\!\! \prod_{e\in
    E_N(\sg)}\E^{(0,\infty)}_{\ell_{e_-}}\big[\Q_X[e^{-\sigma(W)}\ind_{\{\un{W}\geq
    -\eps\}}]\Q_X[e^{-\sigma(W)}]\big]\\
  & &\times \prod_{e\in E_J(\sg)}\mathbb{B}^+_{\ell_{e_-}\!\!\to
    \ell_{e_+}}\big[\Q_X[e^{-\sigma(W)}\ind_{\{\un{W}\geq
    -\eps\}}]\Q_X[e^{-\sigma(W)}]\big]
  \\
  & &\times \prod_{e\in E_I(\sg)}\mathbb{B}_{\ell_{e_-}\!\!\to
    \ell_{e_+}}\big[\Q_X[e^{-\sigma(W)}\ind_{\{\un{W}\geq
    -\eps\}}]^2\big] \prod_{e\in E_O(\sg)}\mathbb{B}_{\ell_{e_-}\!\!\to
    \ell_{e_+}}\big[\Q_X[e^{-\sigma(W)}]^2\big]\, .
\end{eqnarray*}

By Lemmas \ref{sec:some-estimates-2} and \ref{sec:some-estimates-1},
we obtain
\begin{eqnarray*}
  \lefteqn{\CLM^{(k+1)}\Big(e^{-\sigma}\, ;\, \sg\notin \mathfrak{P}\, ,\, \min_{e\in \vec{E}(f_0)}\un{W}_e\geq
    -\eps\Big)}\\ & \leq  & C\sum_{\sg\notin \mathfrak{P}}\int
  \lambda_\sg(\d(\ell_v)_{v\in V(\sg)})\!\! \prod_{e\in
    E_N(\sg)}\Big(\frac{\eps}{\ell_{e_-}+\eps}\Big)^2\prod_{e\in
    E_I(\sg)\cup E_J(\sg)}(\ell_{e_-}+\eps)^{\beta_e}(\ell_{e_+}+\eps)^{1-\beta_e}
\end{eqnarray*} 
for any choice of $\beta_e\in [-2,3]$, that can depend on $e$.  
We then divide every variable
$\ell_v$, with $v$ incident to $f_0$, by $\eps$. We
obtain
\begin{eqnarray}
  \lefteqn{\CLM^{(k+1)}\Big(e^{-\sigma}\, ;\, \sg\notin \mathfrak{P}\, ,\, \min_{e\in \vec{E}(f_0)}{W}_e\geq
    -\eps\Big) \leq  C\sum_{\sg\notin \mathfrak{P}}
    \eps^{\#V_I(\sg)+\#E_I(\sg)+\#E_J(\sg)}} \label{eq:16}\\
  &&\times   \int
  \lambda_\sg(\d(\ell_v)_{v\in V(\sg)})\!\! \prod_{e\in
    E_N(\sg)}\frac{1}{(\ell_{e_-}+1)^2}\prod_{e\in E_I(\sg)\cup
    E_J(\sg)}(\ell_{e_-}+1)^{\beta_e}(\ell_{e_+}+1)^{1-\beta_e}\, . \nonumber
\end{eqnarray} 
Now let us focus again on the case $k=3$. The predominant schemes are
the ones that are obtained from the first pre-scheme of Figure
\ref{fig:shapes}, by adding three null vertices inside each edge
incident to $f_0$, and then labeling the three ``interior'' faces by
$f_1,f_2,f_3$, and choosing a root. 

\begin{figure}
  \begin{center}
    \includegraphics{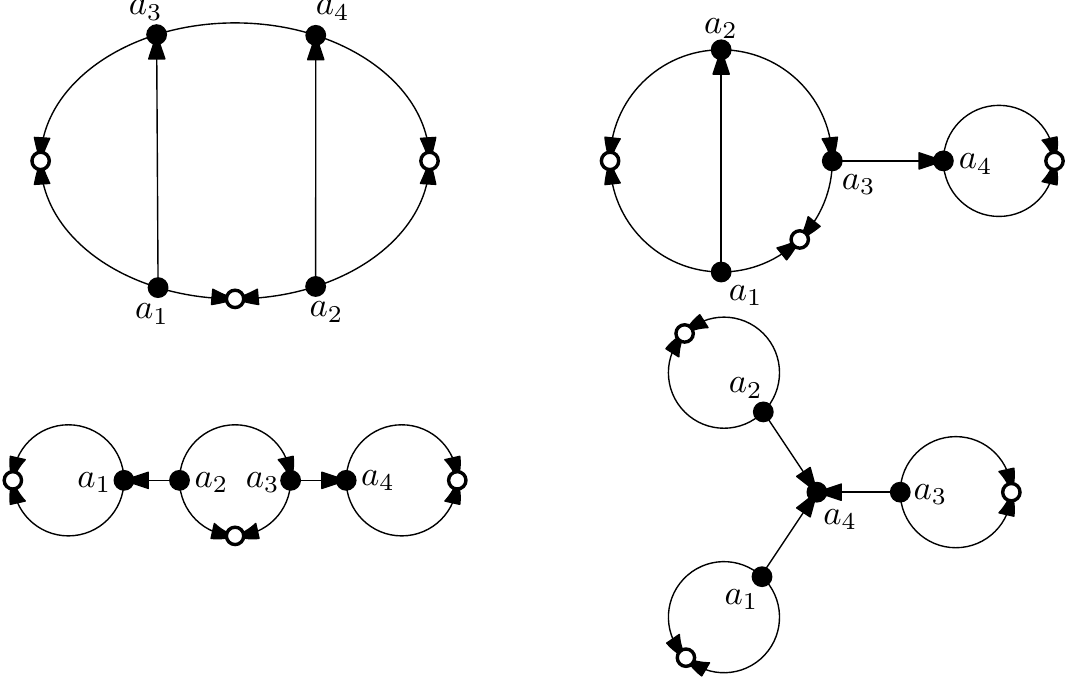}
  \end{center}
  \caption{The four dominant, non-predominant schemes with $4$ faces,
    $f_0$ being the external face, and considered up to obvious
    symmetries. Blank vertices indicate elements of $V_N(\sg)$. The
    labels $a_1,a_2,a_3,a_4$ are the integration variables
    used in the proof of Lemma \ref{sec:proof-lemm-refs}.}
\label{fig:schemes3}
\end{figure}

All other (dominant) schemes are
not predominant, and are indicated in Figure \ref{fig:schemes3}.
Let us consider the first scheme in this figure.  In
this case, we have $\#V_I(\sg)=4,\#E_J(\sg)=1,\#E_I(\sg)=0$, and
taking $\beta_e=1/2$, where $e$ is the unique edge of $\#E_J(\sg)$,
the contribution to the upper-bound \eqref{eq:16}, is bounded by
$$\eps^5\int_{\R_+^4}\frac{\d a_1\d a_2\d a_3\d a_4}{(a_1+1)^4
  (a_2+1)^4(a_3+1)^{3/2}(a_4+1)^{3/2}}\leq C\eps^5\, .$$ Let us turn
to the dominant schemes corresponding to the third pre-scheme of
Figure \ref{fig:shapes}. In this case, one has
$\#V_I(\sg)=4,\#E_J(\sg)=1,\#E_I(\sg)=1$, and, taking $\beta_e=1/3$, the
contribution to \eqref{eq:16} is bounded by
$$\eps^6\int_{\R_+^4}\frac{\d a_1\d a_2\d a_3\d a_4}{(a_1+1)^4
  (a_2+1)^{4/3}(a_3+1)^{4/3}(a_4+1)^{10/3}}\leq C\eps^6\, .$$ The
dominant schemes corresponding to the fourth pre-scheme of Figure
\ref{fig:shapes} have $\#V_I(\sg)=4,\#E_J(\sg)=1,\#E_I(\sg)=2$, and
taking alternatively $\beta_e\in \{1/3,1/2\}$ for
the three edges in $E_I(\sg)\cup E_J(\sg)$, the
contribution to \eqref{eq:16} is at most 
$$\eps^7\int_{\R_+^4}\frac{\d a_1\d a_2\d a_3\d
  a_4}{(a_1+1)^{10/3}(a_2+1)^{7/6}(a_3+1)^{7/6}(a_4+1)^{10/3}}\leq
C\eps^7\, .$$ The dominant schemes corresponding to the fifth
pre-scheme of Figure \ref{fig:shapes} have
$\#V_I(\sg)=4,\#E_J(\sg)=0,\#E_I(\sg)=3$, and taking $\beta_e=-1$ for
every edge in $E_I(\sg)$, the contribution to \eqref{eq:16} is bounded
by
$$\eps^7\int_{\R_+^4}\frac{\d a_1\d a_2\d a_3\d
  a_4}{(a_1+1)^2(a_2+1)^2(a_3+1)^2(a_4+1)^3}\leq C\eps^7\, .$$ This
entails the result.  \cq

\begin{lmm}\label{sec:eps-geodesic-stars-1}
  There exist finite $C,\chi>0$ such that for every $\eps>0$,
$$\CLM_1^{(4)}(\mathcal{B}_2(\eps))\leq C\eps^{4+\chi}\, .$$
\end{lmm}

\proof In all this proof, the mention of 1., 2., 3., 4., will refer to
the four points in the definition of $\mathcal{B}_2(\eps)$, for some
$\eps>0$. 

Using scaling again as in \eqref{eq:19}, it suffices to prove a
similar bound for
$\CLM^{(4)}(e^{-\sigma}\ind_{\mathcal{B}_2(\eps)})$. Let us introduce
some notation. We let $T^{e_{13}}_z=\inf\{s\geq 0:M_{e_{13}}(s)=z\}\in
[0,\infty]$ for every $z\in \R$. We define $T^{e_{23}}_z$ in an
analogous way.  Let 
$$\Xi=-\Big(\min_{e\in \vec{E}(f_3):\ov{e}\in \vec{E}(f_0)}\un{W}_e\wedge\inf\{\un{W}_{e_{13}}^{(s)},0\leq
s\leq T^{e_{13}}_0\}\wedge \inf\{\un{W}_{e_{23}}^{(s)},0\leq
s\leq T^{e_{23}}_0\}\Big)$$
and $\ov{\Delta}_y=\sup_{0\leq z\leq y}\Delta_z$, where 
$$\Delta_z=-z-\Big(\inf\{\un{W}_{e_{13}}^{(t)},T^{e_{13}}_{(-z)-}\leq
t\leq T^{e_{13}}_{-z}\}\wedge
\inf\{\un{W}_{e_{23}}^{(t)},T^{e_{23}}_{(-z)-}\leq s\leq
T^{e_{23}}_{-z}\}\Big)\, .$$
Then, by taking $t$ in \eqref{eq:44} to be of the form
$T^{e_{13}}_{-y}$, as long as $y>0$ is such that
$T^{e_{13}}_{-y}<\infty$, and by choosing $t'=T^{e_{23}}_{-y}$ in a
similar way, we obtain that 
$$(-\Xi)\wedge (-\ov{\Delta}_y-y) \leq (-\Xi)\wedge \inf_{0\leq z\leq
  y}(-\Delta_z-z)\leq  -2y\, .$$
Let $\alpha,\eta,\eta',\eta''$ be positive numbers, all strictly
larger than $\eps$, and such that $\eta>\eta'$. Their values will
be fixed later to be appropriate powers of $\eps$. We observe that
$\mathcal{B}_2(\eps)$ is contained in the union of the following three
events
$\mathcal{B}_2'(\eps),\mathcal{B}_2''(\eps),\mathcal{B}_2'''(\eps)$,
which are defined by points 1., 2., 3.\ in the definition of
$\mathcal{B}_2(\eps)$, together with $4'$., $4''.$\ and $4'''.$\,
respectively, where
\begin{enumerate}
\item[$4'.$] it holds that either $r_{e_{13}}\leq \eta$ or
  $r_{e_{23}}\leq\eta$ or $T^{e_{13}}_{-\alpha}>\eta'$ or
  $T^{e_{23}}_{-\alpha}>\eta'$,
\item[$4''.$] it holds that $r_{e_{13}}\wedge r_{e_{23}}>\eta$,
  $T^{e_{13}}_{-\alpha}\vee T^{e_{23}}_{-\alpha}\leq
  \eta'$, and $\Xi>\eta''$,
\item[$4'''.$] it holds that $r_{e_{13}}\wedge r_{e_{23}}>\eta$,
  $T^{e_{13}}_{-\alpha}\vee T^{e_{23}}_{-\alpha}\leq
  \eta'$, and $\eta''\vee
  \ov{\Delta}_y\geq y$ for every $y\in [0,\eta']$
\end{enumerate}
It remains to estimate separately the quantities
$\CLM^{(4)}(e^{-\sigma}\ind_{\mathcal{B}'_2(\eps)}),\CLM^{(4)}(e^{-\sigma}\ind_{\mathcal{B}''_2(\eps)})$
and $\CLM^{(4)}(e^{-\sigma}\ind_{\mathcal{B}'''_2(\eps)})$. For this,
it suffices to restrict our attention to the event that $\sg$ is the
predominant scheme of Figure \ref{fig:scheme0}, since the other
predominant schemes are the same up to symmetries. Moreover, let us
observe that the points $4'.$, $4''.$ and $4'''.$ only involve the
snakes $W_e$ where $e$ or its reversal is incident to $f_3$, and not
the others. Therefore, when writing the above three quantities
according to the defintion \eqref{eq:19} of $\CLM^{(4)}$, there are
going to be a certain number of common factors, namely, those which
correspond to the contribution of points 2., 3.\ to the $5$ edges of
$\sg$ that are not incident to $f_3$. Renaming the labels $\ell_v,v\in
V(\sg)\setminus V_N(\sg)$ as $a_1,a_2,a_3,a_4$ as indicated in Figure
\ref{fig:scheme0}, we obtain that these common factors are
\begin{itemize}
\item a factor 
$$\E^{(0,\infty)}_{a_4}\big[\Q_X[e^{-\sigma(W)}]\Q_X[e^{-\sigma(W)}\ind_{\{\un{W}\geq
  -\eps\}}]\big]^2\leq \Big(\frac{\eps}{a_4+\eps}\Big)^4$$ corresponding to
the contribution of 2.\ to the two edges incident to $f_0$ and ending
at the vertex with label $a_4$, where we used \eqref{eq:32}, 
\item a factor
\begin{eqnarray*}
\lefteqn{\!\!\!\!\!\!\!\!\!\!\!\!\!\!\!\!\!\!\!\!\!\!\!\!\!\!\!\!\!\!\!\!\!\!\!\!\!\!\!\!\!\!\!\!\!\!\!\!\!\!\!\!\!\!\!\!\!\!\!\!\!\!\!\!\!\!\!\!\!\!\!\!\!\!\E^{(0,\infty)}_{a_1}\big[\Q_X[e^{-\sigma(W)}]\Q_X[e^{-\sigma(W)}\ind_{\{\un{W}\geq
  -\eps\}}]\big]\,
\E^{(0,\infty)}_{a_2}\big[\Q_X[e^{-\sigma(W)}]\Q_X[e^{-\sigma(W)}\ind_{\{\un{W}\geq
  -\eps\}}]\big]}\\
& \leq& \Big(\frac{\eps}{a_1+\eps}\Big)^2
\Big(\frac{\eps}{a_2+\eps}\Big)^2
\end{eqnarray*}
corresponding to the contribution
to 2.\ of the two edges incident to $f_1$ and $f_2$ that end at the
vertices with labels $a_1$ and $a_2$, where we used \eqref{eq:32}
again, 
\item
a factor
$$\B_{a_4\to a_3}\big[\Q_X[e^{-\sigma(W)}]^2\ind_{\{\un{X}\geq
  0\}}\big]\leq 2a_4$$ corresponding to the contribution of 3.\ to the
edge $e_{12}$, and where we used the fact that
$\Q_X[e^{-\sigma(W)}]=e^{-\zeta(X)\sqrt{2}}$ and \eqref{eq:30}. 
\end{itemize}
Let us now bound
$\CLM^{(4)}(e^{-\sigma}\ind_{\mathcal{B}'_2(\eps)})$. The condition
$4'.$ involves only the edges $e_{13}$ and $e_{23}$, so that the
contribution of 2.\ to the two edges incident to $f_0$ and $f_3$ will
bring further factors bounded by $(\eps/(\eps+a_1))^2$ and
$(\eps/(\eps+a_2))^2$. Also, by symmetry, up to a factor $2$, we can
estimate the contribution only of the event $\{r_{13}\leq \eta\mbox{
  or }T^{e_{13}}_{-\alpha}>\eta'\}$, and ignore the rest, so that the
edge $e_{23}$ participates by a factor $\B_{a_1\to
  a_3}[\Q_X[e^{-\sigma(W)}]^2]$, which is bounded by \eqref{eq:29},
while $e_{13}$ contributes by a factor
\begin{eqnarray}
  \lefteqn{\B_{a_1\to
      a_3}\big[\Q_X[e^{-\sigma(W)}]^2\ind_{\{\zeta(X)\leq \eta\mbox{
        or }T_{-\alpha}>\eta'\}}\big]}\nonumber\\
  &\leq & 
  \B_{a_1\to
    a_3}\big[\Q_X[e^{-\sigma(W)}]^2(\ind_{\{\zeta(X)\leq \eta\}}+
  \ind_{\{\zeta(X)>\eta,T_{-\alpha}>\eta'\}})\big]\nonumber\\
  & \leq & \B_{ a_1\to
    a_3}[\zeta(X)\leq \eta]+\B_{ a_1\to
    a_3}[e^{-\zeta(X)\sqrt{2}}\ind_{\{\zeta(X)>\eta,
    T_{-\alpha}>\eta'\}}]\label{eq:25}
\end{eqnarray}
Now on the one hand, since $p_r(x,y)$ is a probability density, we
have
$$\int_{\R_+}\d a_3\, \B_{ a_1\to
  a_3}[\zeta(X)\leq \eta]\leq \int_0^{\eta}\d r\int_\R \d
a_3\, p_r(a_1,a_3)=\eta\, ,$$ and on the other hand
$$\B_{a_1\to
  a_3}[e^{-\zeta(X)\sqrt{2}}\ind_{\{\zeta(X)>\eta,T_{-\alpha}>\eta'\}}]
\leq \int_{\eta}^{\infty}\d r\, e^{-r\sqrt{2}}p_r(
a_1,a_3)\P_{ a_1}^r(T_{-\alpha}>\eta')\, .$$ From the
Markovian bridge construction of \cite{fpy92}, we have, for every
$r>\eta$,
$$\P_{
  a_1\to a_3}^r(T_{-\alpha}>\eta')=\E_{
    a_1}\Big[\ind_{\{T_{-\alpha}>\eta'\}}\frac{ p_{r-\eta'}(X(\eta'),a_3)}{p_r(
  a_1,a_3)}\Big]\, ,$$
so that 
$$\int_{\R_+}\d a_3\int_{\eta}^{\infty}\d r\, e^{-r\sqrt{2}}p_r( a_1,a_3)\P_{
  a_1\to a_3}^r(T_{-\alpha}>\eta')\leq \P_{
  a_1}(T_{-\alpha}>\eta')\int_0^\infty\d r\, e^{-r\sqrt{2}}\,
.$$ Now, we have, by symmetry and scaling of Brownian motion, and
since $T_1$ has same distribution as $X_1^{-2}$ under $\P_0$,
\begin{eqnarray*}
\P_{ a_1}(T_{-\alpha}>\eta')&=&\P_0(T_{
  a_1+\alpha}>\eta')\\
&=& \P_0\Big(T_1>\frac{\eta'}{( a_1+\alpha)^2}\Big)\\
&=&\P_0\Big(|X_1|<\frac{ a_1+\alpha}{\sqrt{\eta'}}\Big)\\
&\leq& C \, \frac{ a_1+\alpha}{\sqrt{\eta'}}\, ,
\end{eqnarray*}
for every $\eps\in (0,1)$ and $a_1>0$. Hence the integral of
\eqref{eq:25} with respect to $a_3\in \R$ is bounded by
$C(\eta+(a_1+\alpha)/\sqrt{\eta'})$. By putting together all the
factors, recalling that $\alpha>\eps$, we have obtained,
\begin{eqnarray}
  \lefteqn{  \CLM^{(4)}(e^{-\sigma}\ind_{\mathcal{B}_2'(\eps)})}\nonumber\\
  &\leq &
  C\int_{\R_+^3}\d a_1\d a_2\d a_4\Big(\frac{\eps}{a_1+\eps}\Big)^4
  \Big(\frac{\eps}{a_2+\eps}\Big)^4\Big(\frac{\eps}{a_4+\eps}\Big)^4a_4\Big(\eta+\frac{a_1+\alpha}{\sqrt{\eta'}}\Big)\nonumber\\
  &\leq & C\eps^4\int_{\R_+}\frac{\d a_1}{(a_1+1)^4}\Big(\eta+\frac{\eps
    a_1+\alpha}{\sqrt{\eta'}}\Big)\, .\nonumber\\
  &\leq & C\eps^4\Big(\eta\vee
  \frac{\alpha}{\sqrt{\eta'}}\Big)\int_{\R_+}\frac{\d
    a_1}{(a_1+1)^4}\Big(2+\frac{\eps}{\alpha}a_1\Big)\nonumber \\
&\leq & C\eps^4\Big(\eta\vee
  \frac{\alpha}{\sqrt{\eta'}}\Big)\label{eq:35}
\end{eqnarray}

Let us now turn to the estimation of
$\CLM^{(4)}(e^{-\sigma}\ind_{\mathcal{B}''_2(\eps)})$. We first
observe the following absolute continuity-type bound: For any
non-negative measurable function $F$ and $\lambda>0$
\begin{equation}
  \int_\R\d y\, \B_{x\to y}[e^{-\lambda \zeta(X)}F(X(s),0\leq s\leq
  T_{-\alpha})
  \ind_{\{\zeta(X)>\eta,T_{-\alpha}\leq \eta'\}}]
  \leq \frac{1}{\lambda}\E_x^{(-\alpha,\infty)}[F(X)]\, .\label{eq:26}
\end{equation}
Indeed, note
that by using again the Markovian bridge description of \cite{fpy92}, 
\begin{eqnarray*}
  \lefteqn{\int_\R\d y\, \B_{x\to y}[e^{-\lambda \zeta(X)}F(X(s),0\leq s\leq T_{-\alpha})
    \ind_{\{\zeta(X)>\eta,T_{-\alpha}\leq \eta'\}}]}\\
  &=&\int_\R\d y\int_{\eta}^{\infty}e^{-\lambda r}\d r\,
  \E_x[F(X(s),0\leq s\leq T_{-\alpha})\ind_{\{T_{-\alpha}\leq
    \eta'\}}p_{r-\eta'}(X(\eta'),y)]\\
&\leq &\int_0^\infty e^{-\lambda r}\d r\, \E_x[F(X(s),0\leq s\leq
T_{-\alpha})]\, ,
\end{eqnarray*}
as wanted.  Back to
$\CLM^{(4)}(e^{-\sigma}\ind_{\mathcal{B}''_2(\eps)})$, note that there
are exactly two edges incident to $f_3$ with reversal incident to
$f_0$, and we let $e_1$ be the one that is linked to the vertex with
label $a_1$, and $e_2$ the one linked to the vertex with label $a_2$.
Note that the event $\{\Xi>\eta''\}$ is then included in the union
$$\{\un{W}_{e_1}<-\eta''\}\cup
\{\un{W}_{e_2}<-\eta''\}\cup\Big\{\inf_{0\leq s\leq
  T_0^{e_{13}}}\un{W}^{(s)}_{e_{13}}<-\eta''\Big\}\cup
\Big\{\inf_{0\leq s\leq
  T_0^{e_{23}}}\un{W}^{(s)}_{e_{23}}<-\eta''\Big\}\, .$$ For
symmetry reasons, the first two have same contribution (after
integrating with respect to $a_1,a_2,a_3,a_4$), as well as the last
two. The contribution of the edges $e_1$ and $e_2$ to
$\{\un{W}_{e_1}<-\eta''\}$ is then bounded by (we use also 2., and
recall that $\eta>\eta'>\eps$)
$$\E_{a_1}^{(0,\infty)}\big[\Q_X(\un{W}\geq -\eps)\Q_X(\un{W}\leq
-\eta'')\big] \Big(\frac{\eps}{a_2+\eps}\Big)^2 \leq
C\Big(\frac{\eps}{\eta''-\eps}\Big)^2\Big(\frac{\eps}{a_2+\eps}\Big)^2$$
by \eqref{eq:33}, while the edges $e_{13}$ and $e_{23}$ contribute
factors of the form $\B_{a_1\to a_3}[\Q_X[e^{-\sigma(W)}]^2]$ and
$\B_{a_2\to a_3}[\Q_X[e^{-\sigma(W)}]^2]$, which are bounded after
integration of the variable $a_3$, by \eqref{eq:29}. Hence, 
\begin{eqnarray*}
  \lefteqn{\CLM^{(4)}(e^{-\sigma}\ind_{\mathcal{B}''_2(\eps)}\ind_{\{\un{W}_{e_1}<-\eta''\}})}
  \\
  &\leq & 
  C\Big(\frac{\eps}{(\eta''-\eps)}\Big)^2\int_{\R_+^3}\d a_1\, \d a_2\, \d a_4\,
  \Big(\frac{\eps}{a_1+\eps}\Big)^2
  \Big(\frac{\eps}{a_2+\eps}\Big)^4 
  \Big(\frac{\eps}{a_4+\eps}\Big)^4 a_4 \\
&\leq &C\eps^4\Big(\frac{\eps}{\eta''-\eps}\Big)^2\, ,
\end{eqnarray*}
On the other hand, the edges $e_1$ and $e_2$ do
not contribute to $\{\inf\{\un{W}^{(s)}_{e_{13}},0\leq s\leq
T_0^{e_{13}}\}<-\eta''\}$, and involve only, via 2., a factor
$(\eps/(a_1+\eps))^2(\eps/(a_2+\eps))^2$. The contribution of $e_{13}$
and $e_{23}$, integrated in $a_3$, is bounded by
\begin{eqnarray*}
&& \!\!\!\!\!\!\!\!\!\!\int_\R \d a_3\, \B_{a_1\to
  a_3}\Big[\Q_X[e^{-\sigma}]\Q_X\big[\inf_{0\leq s\leq
  T_0}\un{W}^{(s)}<-\eta''\big]\ind_{\{\zeta(X)>\eta,T_{-\alpha}\leq
  \eta'\}}\Big]\B_{a_2\to a_3}[e^{-\sigma}]\\
&&\leq  C \E_{a_1}^{(-\alpha,\infty)}\Big[\Q_X\big[\inf_{0\leq s\leq
  T_0}\un{W}^{(s)}<-\eta''\big]\Big]\\
&&\leq C P(\ov{\Delta}_{a_1}>\eta''+a_1)\\
&&\leq 
C\Big(1-\exp\Big(-\frac{2a_1}{a_1+\eta''}\Big)\Big)\\
&&\leq C\frac{a_1}{\eta''}\, ,
\end{eqnarray*}
where we used \eqref{eq:25} in the second step, and Lemma
\ref{sec:some-estimates-3} in the third step: Here, under $P$,
$(\Delta_t,t\geq 0)$ is a Poisson process with intensity $2\d a/a^2$,
and $\ov{\Delta}$ is its supremum process. We have obtained
\begin{eqnarray*}
  \lefteqn{\CLM^{(4)}(e^{-\sigma}\ind_{\mathcal{B}''_2(\eps)}\ind_{\{\inf\{\un{W}^{(s)},0\leq
      s\leq T^{e_{13}}_0\}<-\eta''\}})}\\
  &\leq & C\int_{\R_+^3}\d a_1\, \d a_2\, \d
  a_4\Big(\frac{\eps}{a_1+\eps}\Big)^4
  \frac{a_1}{\eta''}\Big(\frac{\eps}{a_2+\eps}\Big)^4\Big(\frac{\eps}{a_4+\eps}\Big)^4a_4\\
&\leq & C\frac{\eps^{5}}{\eta''}\int_{\R_+}\d
a_1\frac{a_1}{(a_1+1)^4}\\
&\leq & C\frac{\eps^5}{\eta''}\, .
\end{eqnarray*}
These
estimations, together with our preliminary remarks, entail that
\begin{equation}
\label{eq:37}
\CLM^{(4)}(e^{-\sigma}\ind_{\mathcal{B}_2''(\eps)})\leq
C\eps^4\bigg(\Big(\frac{\eps}{\eta''-\eps}\Big)^2+\frac{\eps}{\eta''}\bigg)\,
.
\end{equation}

Finally, let us consider
$\CLM^{(4)}(e^{-\sigma}\ind_{\mathcal{B}_2'''(\eps)})$. Points 2.\ and
3.\ induce contributions of the edges incident to $f_0$, as well as
the edge $e_{12}$, that are bounded by 
$$2\Big(\frac{\eps}{a_1+\eps}\Big)^4 \Big(\frac{\eps}{a_2+\eps}\Big)^4
\Big(\frac{\eps}{a_4+\eps}\Big)^4 a_4\, .$$
Now, point $4'''.$ involves only the edges $e_{13}$ and $e_{23}$, and
contributes by a factor bounded above by 
$$\B_{a_1\to
  a_3}\Big[e^{-\zeta(X)\sqrt{2}}\ind_{\{\zeta(X)>\eta,T_{-\alpha}\leq
  \eta'\}}\int_{\CC(\CC(\R))}\Q_X(\d W)\Q_X(\d W')\ind_{\{\eta''\vee \ov{I}_y\vee
\ov{I}'_y\geq y, 0\leq y\leq \alpha\}}\Big]\, ,$$ where $I$ was defined
in Lemma \ref{sec:some-estimates-3}, $\ov{I}_y=\sup_{0\leq z\leq
  y}I_y$, and $I',\ov{I}'$ are defined in a similar way from the
trajectory $W'$. Now we use again a bound with same spirit as
\eqref{eq:25}. Namely, for every $\lambda >0$, every $x,y\in \R$ and
for non-negative measurable $F$, we have
\begin{eqnarray*}
  \lefteqn{\int_{\R}\d z\, \int_{\CC(\R)^2}\B_{x\to z}(\d X) \B_{y\to z}(\d X')\,
  e^{-\lambda \zeta(X)-\lambda \zeta(X')}}\\
& &\times   \ind_{\{\zeta(X)\wedge \zeta(X')>\eta,T_{-\alpha}\vee T'_{-\alpha}\leq
    \eta'\}} F\big((X(s))_{0\leq s\leq
    T_{-\alpha}},(X'(s))_{0\leq s\leq
    T'_{-\alpha}}\big)\\
 & \leq  &
  \frac{1}{\lambda^2\sqrt{2\pi(\eta-\eta')}}\int_{\CC(\R)^2}
\P^{(-\alpha,\infty)}_x(\d X)\P^{(-\alpha,\infty)}_y(\d X')F(X,X')\, , 
\end{eqnarray*}
where $T'_z$ is the first hitting time of $z$ by $X'$. 
This is obtained in a similar way to \eqref{eq:25}, writing the
left-hand side as  
\begin{eqnarray*}
  \lefteqn{\int_{\R}\d z\int_{(\eta,\infty)^2}\d r\d r'e^{-\lambda
      (r+r')}\int_{\CC(\R)^2}\P_x(\d X)\P_y(\d X') }\\
  & & \times \ind_{\{T_{-\alpha}\vee T'_{-\alpha}<\eta'\}}F\big((X(s))_{0\leq s\leq
    T_{-\alpha}},(X'(s))_{0\leq s\leq
    T'_{-\alpha}}\big)p_{r-\eta'}(X(\eta'),z) p_{r-\eta'}(X'(\eta'),z)
\end{eqnarray*}
This is obtained by first checking this for $F$ of a product form and
using the Markovian description of bridges, and then applying a
monotone class argument.  We then use the bound
$p_{r-\eta'}(X'(\eta'),z)\leq (2\pi (\eta-\eta'))^{-1/2}$, valid for
$r\geq \eta>\eta'$, and use Fubini's theorem to integrate
$p_{r-\eta'}(X(\eta'),z)$ with respect to $z$, as in the derivation of
\eqref{eq:26}. Therefore, after integrating with respect to the
variables $a_3$, we obtain that the edges $e_{13}$ and $e_{23}$
together contribute by
$$C\int_{\CC(\R)^2}\!\!\!\P_{a_1}^{(-\alpha,\infty)}(\d X)
\P_{a_2}^{(-\alpha,\infty)}(\d X') \int_{\CC(\CC(\R))^2}\!\!\!\Q_X(\d
W)\Q_{X'}(\d W')\ind_{\{\eta''\vee \ov{I}_y\vee \ov{I}'_y\geq y\, ,\,
  0\leq y\leq \alpha\}}\, ,$$ where $I'$ is defined from $X'$ as $I$
was defined from $X$. This equals
$$C\int_{\CC(\R)^2}\!\!\!\P_0^{(-\alpha,\infty)}(\d X)
\P_0^{(-\alpha,\infty)}(\d X') \int_{\CC(\CC(\R))^2}\!\!\!\Q_X(\d
W)\Q_{X'}(\d W')\ind_{\{\eta''\vee \ov{I}_y\vee \ov{I}'_y\geq y\, ,\,
  0\leq y\leq \alpha\}}\, ,$$ by an application of the Markov
property, noticing that $\ov{I}_y$ and $\ov{I}'_y$ only involve the
processes $W^{(s)}$ for $T_0\leq s\leq T_{-y}$, and similarly for
$\ov{I}'_y$ (we skip the details).  By Lemma
\ref{sec:some-estimates-3} and standard properties of Poisson
measures, the process $(I_y\vee I'_y,0\leq y\leq \alpha)$ under the
law $\E^{(-\alpha,\infty)}_0\otimes \E^{(-\alpha,\infty)}_0[\Q_X(\d
W)\Q_X(\d W')]$ is a Poisson process on the time-interval $[0,\alpha]$
with intensity $4\d a/a^2$.  By Lemma \ref{sec:fast-coal-geod} we
obtain that the last displayed quantity is less than
$(\eta''/\alpha)^{e^{-4}/4}$, and we conclude that
$$\CLM^{(4)}(e^{-\sigma}\ind_{\mathcal{B}'''_2(\eps)})\leq C\eps^4 \frac{1}{\sqrt{\eta-\eta'}}\Big(\frac{\eta''}{\alpha}\Big)^{e^{-4}/4} $$

This, together with \eqref{eq:35} and \eqref{eq:37}, finally entails
that
$$\CLM^{(4)}(e^{-\sigma}\ind_{\mathcal{B}_2(\eps)})\leq C\eps^4\bigg(\eta\vee
\frac{\alpha}{\sqrt{\eta'}} +
\Big(\frac{\eps}{\eta''-\eps}\Big)^2+\frac{\eps}{\eta''} +
\frac{1}{\sqrt{\eta-\eta'}}\Big(\frac{\eta''}{\alpha}\Big)^{e^{-4}/4} \bigg)\, .$$
Let us now choose $\alpha,\eta,\eta',\eta''$ of the form 
$$\alpha=\eps^\beta\, ,\quad \eta=\eps^\nu\, ,\quad
\eta'=\eps^{\nu'}\, ,\quad \eta''=\eps^{\nu''}\, ,$$ with
$\beta,\nu,\nu',\nu''\in (0,1)$, and let us assume for the
time being that $\eps<1$. Then the condition $\eta>\eta'$
amounts to $\nu<\nu'$, and by picking $\eps$ even smaller if
necessary (depending on the choice of $\nu''$) we can assume that
$\eta-\eta'>\eta/2$ and $\eta''-\eps>\eta''/2$. The above bound then
becomes
$$C\eps^4(\eps^{\nu}+\eps^{\beta-\nu'/2}+\eps^{1-\nu''}+\eps^{(\nu''-\beta)e^{-4}/4-\nu/2})\, .$$
Therefore, it suffices to choose $\nu',\beta,\nu''$ so that
$0<\nu'<2\beta<2\nu''<2$, and then $\nu$ so that $0<\nu <
\nu'\wedge (\nu''-\beta)e^{-4}/2$, to obtain a bound of the form
$C\eps^{4+\chi}$ with positive $\chi$, as wanted. One the choice is
made, this bound remains obviously valid without restriction on
$\eps$, by taking $C$ larger if necessary.  \cq

\medskip

The combination of \eqref{eq:49} and Lemmas \ref{sec:proof-lemm-refs}
and \ref{sec:eps-geodesic-stars-1} finally entail Lemma
\ref{sec:covering-3-star-1}.

\section{Concluding remarks}\label{sec:concluding-remarks}

\paragraph{Beyond quadrangulations. }
An important problem in random map theory is the question of {\em
  universality}. It is expected that Theorem \ref{sec:introduction-4}
remains valid for much more general families of random maps than
quadrangulations, namely, the {\em regular critical} Boltzmann maps
introduced in \cite{jfmgm05,mierinv} should admit the Brownian map as
a scaling limit. In particular, in his work mentioned at the end of
the Introduction, Le Gall proves Theorem \ref{sec:introduction-4} for
$2p$-angulations for any given $p\geq 2$, that is maps with faces of
degree $2p$. Up to a deterministic, $p$-dependent multiplicative
factor, the limit is still the Brownian map.

We believe that our method is robust enough to tackle this more
general problem. In particular, the results of Sections
\ref{sec:coding-labeled-maps} could be adapted to more general maps
using variants of the Bouttier-Di Francesco-Guitter bijection
\cite{BdFGmobiles} that incorporate multiple points. The labeled maps
that would intervene would be more complicated than the ones of the
present paper, but we expect that the scaling limits of Section
\ref{sec:scal-limits-label} remain valid. Indeed, such generalizations
hold for the basic case of trees, by the invariance principles
developed in \cite{jfmgm05,miergwmulti}. This program seems reasonable
to carry out at least in the case where quadrangulations are replaced
by random {\em bipartite} maps. Since the currently available
invariance principles for trees coding non-bipartite maps are
considerably weaker than those for bipartite maps, the case of
non-bipartite maps (including random triangulations) would probably be
notably harder to analyze.

Once results analogous to those of Sections
\ref{sec:coding-labeled-maps} and \ref{sec:scal-limits-label} are
obtained, our strategy of proof would be valid without much change as
soon as we have the prior knowledge that (loosely speaking)
\begin{enumerate}
\item any subsequential scaling limit $(S,D)$ of the family of maps
  under consideration is homeomorphic to $(S,D^*)$
\item the estimates for the volumes of balls in $(S,D)$ of Section
  \ref{sec:rough-comp-betw} hold
\item typical geodesics in $(S,D)$ are a.s.\ unique. 
\end{enumerate}
All these points are also quite robust, at least in the bipartite
case. They are known for $2p$-angulations, and likely extend to more
general cases using the approaches of
\cite{legall06,legall08,miertess}.

Last but not least, we believe that similar methods extend to higher
genera, probably to the cost of technical complication. In particular,
the analogs of points 1., 2., 3.\ above have been derived for maps on
orientable compact surfaces in
\cite{miertess,bettinelli10,bettinelli11}.



\paragraph{Stable maps. }A one-parameter family of scaling limits of
maps, different from the Brownian map, can be obtained by considering
Boltzmann distributions on maps for which the degrees of faces have
heavy tails \cite{LGMi09}. The problem of the uniqueness of the
scaling limit is still open for these maps. We do not know if the
methods of the present paper can be adapted to tackle this problem.

\paragraph{Geodesic stars in the Brownian map. }The methods of
Sections \ref{sec:coding-labeled-maps}, \ref{sec:scal-limits-label}
and \ref{sec:proof-key-lemmas} allow to give estimates on the
probability of the event $\mathcal{G}(\eps,k)$ that, if
$x_0,x_1,\ldots,x_k$ are uniformly chosen points in $S$, the geodesics
from $x_0$ to $x_1,\ldots,x_k$ are pairwise disjoint outside of
$B_{D_*}(x_0,\eps)$. Using similar arguments as in Section
\ref{sec:relat-label-maps}, one finds that this probability is of the
same order as
$$\frac{1}{\eps}\CLM_1^{(k+1)}\Big(\min_{e\in \vec{E}(f_0)}\un{W}_e\geq
-\eps\Big)\, .$$ By using bounds of the type \eqref{eq:16}, this is
bounded by $C\eps^{k-1}$, and we think that the exponent is
sharp. Since one needs about $\eps^{-4}$ balls of radius $\eps$ to
cover $(S,D^*)$, possibly up to slowly varying terms, this estimate
seems to indicate that there is an order of $\eps^{k-5}$ such
$\eps$-geodesic stars with $k$ arms in the Brownian map. Letting
$\eps\to 0$, and making a leap of faith, this suggests that the
probability that there exist points $x_1,\ldots,x_k$ satisfying the
geodesic star event $\mathcal{G}(S;x_1,\ldots,x_k)$ of Definition
\ref{sec:rough-comp-betw-1} is $1$ as long as $k\in \{1,2,3,4\}$, and
$0$ when $k\geq 6$, the case $k=5$ being critical and harder to
settle. We hope to study more detailed aspects of geodesic stars in
the Brownian map in future work.

\def\polhk#1{\setbox0=\hbox{#1}{\ooalign{\hidewidth
  \lower1.5ex\hbox{`}\hidewidth\crcr\unhbox0}}}

\end{document}